\documentclass[12pt, reqno]{amsart}
\usepackage[margin=1in]{geometry}
\usepackage{amssymb,latexsym,amsmath,amscd,amsfonts}
\usepackage{latexsym}
\usepackage[mathscr]{eucal}
\usepackage{bm}
\usepackage{mathptmx}
\usepackage{amssymb}
\usepackage{amsthm}
\usepackage{dcolumn}
\usepackage[all]{xy}
\usepackage{enumitem}
\usepackage{tikz}
\usepackage[utf8]{inputenc}

\def \qed {\hfill \vrule height6pt width 6pt depth 0pt}
\def\textmatrix#1&#2\\#3&#4\\{\bigl({#1 \atop #3}\ {#2 \atop #4}\bigr)}
\def\dispmatrix#1&#2\\#3&#4\\{\left({#1 \atop #3}\ {#2 \atop #4}\right)}
\newcommand{\beg}{\begin{equation}}
	\newcommand{\eeg}{\end{equation}}
\newcommand{\ben}{\begin{eqnarray*}}
	\newcommand{\een}{\end{eqnarray*}}

\newtheorem{theorem}{Theorem}[section]
\newtheorem{lemma}[theorem]{Lemma}
\newtheorem{corollary}[theorem]{Corollary}

\newtheorem{proposition}[theorem]{Proposition}

\newtheorem{definition}[theorem]{Definition}

\newtheorem{remark}[theorem]{Remark}

\numberwithin{equation}{section}

\newcommand{\HS}{\mathcal H}

\begin{document}
	\title{Hardy Spaces and Toeplitz operators on the Hexablock}
	
\author{Sourav Pal, Guicong Su and Nitin Tomar}

\address[Sourav Pal]{Mathematics Department, Indian Institute of Technology Bombay,
	Powai, Mumbai - 400076, India.} \email{sourav@math.iitb.ac.in, souravmaths@gmail.com}

\address[Guicong Su]{School of Mathematics and Statistics, Wuhan University of Technology, Wuhan, Hubei 430070, P.R.China}\email{suguicong@whut.edu.cn}

\address[Nitin Tomar]{Statistics and Mathematics Unit, Indian Statistical Institute, Bangalore, Karnataka - 560059, India.} \email{tomarnitin414@gmail.com}

	%\author{Sourav Pal}
	%\address{Department of Mathematics, Indian Institute of Technology Bombay, Powai, Mumbai,
	%	Maharashtra – 400076, India.}
	%\email{sourav@math.iitb.ac.in, souravmaths@gmail.com}
	
	%\thanks{* Corresponding author.}
	%\thanks{ }

	%\author{Guicong Su}
	%\address{School of Mathematics and Statistics, Wuhan University of Technology, Wuhan, Hubei 430070, P.
	%	R. China}
	%\email{suguicong@whut.edu.cn}
	
	%\author{Nitin Tomar}
	%\address{Statistics and Mathematics Unit, Indian Statistical Institute, Bangalore,
	%	Karnataka – 560059, India.}
	%\email{tomarnitin414@gmail.com, nitin\_pd@isibang.ac.in}

	
	
\begin{abstract} 
	The hexablock is a domain in $\mathbb{C}^4$ arising from the $\mu$-synthesis problem and is given by
	\[
	\mathbb{H}=\left\{(a,x_1,x_2,x_3)\in \mathbb{C}\times\mathbb{E}:
	\sup_{z_1,z_2\in\mathbb{D}}
	\left|
	\frac{a\sqrt{(1-|z_1|^2)(1-|z_2|^2)}}
	{1-x_1z_1-x_2z_2+x_3z_1z_2}
	\right|<1
	\right\},
	\]
	where
	$
	\mathbb{E}
	=
	\left\{
	(x_1,x_2,x_3)\in\mathbb{C}^{3}:
	1-x_1z_1-x_2z_2+x_3z_1z_2\neq 0
	\ \text{for all}\ z_1,z_2\in\overline{\mathbb D}
	\right\}
	$
	is the tetrablock. Two sets closely related to $\mathbb{H}$ are the $\mu$-hexablock $\mathbb{H}_{\mu}$ and the normed hexablock $\mathbb{H}_N$ that arise in connection with the $\mu_E$-unit ball and the norm unit ball of $M_2(\mathbb{C})$, respectively. We introduce the Hardy space $H^2(\mathbb H)$ associated with a
	canonical measure $\mu_{\mathbb H}$ on the distinguished boundary $b\mathbb H$ of the hexablock. We identify $H^2(\mathbb H)$ with a closed subspace of
	$L^2(b\mathbb H,\mu_{\mathbb H})$ and show that it is a reproducing kernel Hilbert space. We establish a new description of $\mathbb H$ and obtain a subset of it that becomes a uniqueness set for the algebra $\text{Hol}(\mathbb H)$. In fact,  we prove that $\mathbb H$ is the interior of the polynomial convex hull of $\Psi_{\mathbb H}(\overline{\mathbb D}^4)$, where
\[
	\Psi_{\mathbb H}(\xi_1, \xi_2,\lambda_1,\lambda_2)
	=
	\left(\frac{\lambda_1(\xi_1-1)}{2},
	\frac{\xi_2(\xi_1+1)}{2}, \ 
	\frac{\lambda_1\lambda_2(\xi_1+1)}{2}, \
	\xi_1 \xi_2\lambda_1\lambda_2
	\right).
\]
	We also show that $H^2(\mathbb H)$ embeds isometrically as a closed subspace of
	$H^2(\mathbb D^4)$. We then obtain a Brown-Halmos type
	characterization of Toeplitz operators on $H^2(\mathbb H)$ and a
	commutant lifting theorem for $\mathbb  H$-isometries. As consequences, we determine the
	commutant of the coordinate multipliers, characterize the compact
	Toeplitz operators and show that the coordinate multipliers have no
	nontrivial joint reducing subspaces. Finally, we characterize $\mathbb H$-unitaries as the images of commuting quadruples of unitaries under $\Psi_{\mathbb H}$ and present a class of $\mathbb H$-contractions from commuting quadruples having
	$\overline{\mathbb D}^4$ as a spectral set. We prove an analogous characterization of $\mathbb E$-unitaries and present a class of $\mathbb E$-contractions arising from commuting triples of contractions having $\overline{\mathbb D}^3$ as a spectral set. En route, we also study a weighted Hardy space on the hexablock.
\end{abstract}

\maketitle

\tableofcontents

\section{Introduction}

\noindent
The class of Toeplitz operators occupies a central position in operator theory owing to their close connections with function theory, harmonic analysis and the structure of multiplication operators on Hilbert spaces of holomorphic functions. On the classical Hardy space $H^2(\mathbb D)$, a bounded operator $T$ is Toeplitz  if and only if it satisfies the Brown-Halmos identity
\[
M_z^*TM_z=T,
\]
where $M_z$ denotes multiplication by the coordinate function on $H^2(\mathbb D)$. This characterization was established by Brown and Halmos \cite{BrownHalmos1964}. Since then, it has given rise to the study of Hardy spaces over various domains in $\mathbb C^n$ and to the characterization of the associated Toeplitz operators, typically referred to as Brown-Halmos-type theorems. For example, the analogues of the Brown-Halmos theorem have been obtained for the polydisc \cite{Ding2003, DingSunZheng2012, MajiSarkarSarkar2018}, the Euclidean unit ball \cite{DavieJewell1977} and the Hartogs triangle \cite{JainPramanick2025}. The Hardy space  of the symmetrized bidisc was introduced in \cite{MisraShyamRoyZhang2013} and the associated Toeplitz operators were studied in \cite{BhattacharyyaDasSau2021}. Moreover, the class of Hardy spaces and Toeplitz operators on the tetrablock and the pentablock were introduced in  \cite{JainKumarMalPramanick2026} and \cite{SuTomar}, respectively. Also, the study of Toeplitz operators on bounded symmetric domains and Cartan domains can be found in \cite{Upmeier1983, KumarMalPramanick2025, GhoshShyamRoy2024}. The existing works show that the study of Hardy type spaces on a domain is closely intertwined with its geometric structure, its distinguished boundary and the multiplication operators associated with its coordinate functions.

\smallskip 

In this  direction, there are three notable domains, namely, the symmetrized bidisc $\mathbb G_2$, the tetrablock $\mathbb E$ and the pentablock $\mathbb P$ that are defined as follows:
\begin{align*}
	\mathbb G_2&=\{(\text{tr}(A), \det(A)): A \in M_2(\mathbb C), \|A\|<1\},\\
	\mathbb E&=\{(a_{11}, a_{22}, \det(A)): A \in M_2(\mathbb C), \|A\|<1\},\\
	\mathbb P&=\{(a_{21}, \text{tr}(A), \det(A)): A \in M_2(\mathbb C), \|A\|<1\}.
\end{align*}
Here $\text{tr}(A)$ and $\det(A)$ denotes the trace and determinant of $A \in M_2(\mathbb C)$, respectively. Each of these three domains originated in the context of $\mu$-synthesis problem. Given a linear subspace $E$ of $M_n(\mathbb C)$, the structured singular value associated with $E$ is the map $\mu_E: M_n(\mathbb C)\to[0,\infty)$ defined by
\[
\mu_E(A)=1\slash \inf\{\|X\|\in E,\ \det(I-AX)=0\}
\]
with the understanding that $\mu_E(A)=0$ if there does not exist $X\in E$ with $\det(I-AX)=0$. Here, $\|\cdot\|$ denotes the operator norm on $M_n(\mathbb{C})$ and $I$ is the identity matrix. Given a linear subspace $E\subseteq M_n(\mathbb C)$, the corresponding $\mu$-synthesis problem asks whether, for given distinct points $\lambda_1,\ldots,\lambda_m\in\mathbb D$ and matrices $A_1,\ldots,A_m\in M_n(\mathbb C)$, there exists a holomorphic map $F:\mathbb D\to M_n(\mathbb C)$ satisfying
$F(\lambda_j)=A_j$ for $1\leq j\leq m$ such that $\mu_E(F(\lambda))\le1$ for $\lambda\in\mathbb D$. For further background on control-theoretic origins of the structured singular value, we refer the reader to \cite{DoyleII, Francis}. Let $E_1, E_2$ and $E_3$ be linear subspaces of $M_2(\mathbb C)$ consisting of all scalar matrices, diagonal matrices and upper triangular matrices with same diagonal entries, respectively. The domains $\mathbb G_2$, $\mathbb E$ and $\mathbb P$ are associated with the structured singular values corresponding to $E_1$, $E_2$ and $E_3$, respectively, as follows:
\begin{align*}
	\mathbb G_2&=\{(\text{tr}(A), \det(A)): A \in M_2(\mathbb C), \mu_{E_1}(A)<1\},\\
	\mathbb E&=\{(a_{11}, a_{22}, \det(A)): A \in M_2(\mathbb C), \mu_{E_2}(A)<1\},\\
	\mathbb P&=\{(a_{21}, \text{tr}(A), \det(A)): A \in M_2(\mathbb C), \mu_{E_3}(A)<1\},
\end{align*}
for $A=(a_{ij})_{i, j=1}^2 \in M_2(\mathbb C)$. The domains $\mathbb{G}_2$, $\mathbb{E}$ and $\mathbb{P}$ were introduced in \cite{AglerYoung1999}, \cite{AbouhajarWhiteYoung2007} and \cite{AglerLykovaYoung2015}, respectively. Since the linear subspaces $E_1, E_2$ and $E_3$ are all contained in the subspace $E_4$ of all upper triangular matrices in $M_2(\mathbb C)$, it is natural to consider the $\mu$-synthesis problem corresponding to $E_4$. Motivated by this, the authors of \cite{BiswasPalI} introduced the hexablock $\mathbb H$ defined as
\[
\mathbb{H}=\left\{(a,x_1,x_2,x_3)\in \mathbb{C}\times\mathbb{E}:
\sup_{z_1,z_2\in\mathbb{D}}
\left|
\frac{a\sqrt{(1-|z_1|^2)(1-|z_2|^2)}}
{1-x_1z_1-x_2z_2+x_3z_1z_2}
\right|<1
\right\}.
\]
The set $\mathbb H$ is a non-convex, polynomially convex and $(1,1, 1,2)$-quasi-balanced domain in $\mathbb C^4$. To explain the connection between the hexablock $\mathbb H$ and the $\mu$-synthesis problem, we recall from \cite{BiswasPalI} the sets
\[
\mathbb{H}_{\mu}=\{\pi_{\mathbb H}(A): A\in M_2(\mathbb C),\ \mu_{E_4}(A)<1\},
\ \ 
\mathbb{H}_N=\{\pi_{\mathbb H}(A): A\in M_2(\mathbb C),\ \|A\|<1\},
\]
where $\pi_{\mathbb H}(A)=(a_{21},a_{11},a_{22},\det(A))$ for $A=(a_{ij})_{i,j=1}^2 \in M_2(\mathbb C)$. These sets are known as the $\mu$-hexablock and the normed hexablock, respectively. Unlike the domains $\mathbb G$, $\mathbb E$ and $\mathbb P$, which admit equivalent descriptions using the corresponding $\mu_E$-unit ball and the operator norm unit ball of $M_2(\mathbb C)$, the sets $\mathbb{H}_\mu$ and $\mathbb{H}_N$ are different. Moreover, neither of them is a domain in $\mathbb C^4$. However, the hexablock $\mathbb H$ can be recovered from their polynomial convex hulls through the identities
\[
\mathbb H
=
\operatorname{int}(\overline{\mathbb{H}}_\mu)
=
\operatorname{int}\left(\widehat{\overline{\mathbb H}}_N\right),
\]
where $\widehat{\overline{\mathbb{H}}}_N$ denotes the polynomial convex hull of $\mathbb{H}_{N}$ and $\operatorname{int}(X)$ denotes the interior of $X$. Further details on the origin and geometric properties of $\mathbb H$ can be found in \cite{BiswasPalI, BiswasPalII, Su_Hexa}. 

\smallskip 

The purpose of this paper is to construct a Hardy space on the hexablock $\mathbb H$ with particular emphasis on the associated Toeplitz operators and the operator-theoretic consequences of a particularly interesting natural parametrization map. Indeed, a central role in our approach is played by the parametrization map $\Psi_{\mathbb H}: \mathbb C^4 \to \mathbb C^4$ defined as
\[
\Psi_{\mathbb H}(\xi_1,\xi_2,\lambda_1,\lambda_2)
=
\left(
\frac{\lambda_1(\xi_1-1)}{2},
\frac{\xi_2(\xi_1+1)}{2},
\frac{\lambda_1\lambda_2(\xi_1+1)}{2},
\xi_1\xi_2\lambda_1\lambda_2
\right).
\]
We prove in Section \ref{sec002_hex} that $\Psi_{\mathbb H}(\overline{\mathbb D} \times \mathbb D^2 \times \overline{\mathbb D})\subseteq\mathbb H$ and $\Psi_{\mathbb H}(\mathbb T^4)=b\mathbb H$. Here, $b \mathbb H$ is the distinguished boundary of $\mathbb H$. Furthermore,
\[
\widehat{\Psi_{\mathbb H}(\overline{\mathbb D}^{\,4})}=\overline{\mathbb H}
\quad \text{and} \quad \text{int}\left(\widehat{\Psi_{\mathbb H}(\overline{\mathbb D}^{\,4})}\right)=\mathbb H.
\]
We further prove in Section \ref{sec002_hex} that $\Psi_{\mathbb H}(\mathbb D^4)$ is a uniqueness set for the algebra $\text{Hol}(\mathbb H)$ of holomorphic functions on $\mathbb H$ in the sense that $f=0$ on $\mathbb H$ whenever $f \in \text{Hol}(\mathbb H)$ with $f|_{\Psi_{\mathbb H}(\mathbb D^4)}=0$. We present an analogous map for the tetrablock $\mathbb E$, which is given by
\[
\Psi_{\mathbb E}(\lambda_1, \lambda_2, \xi)=\left(\frac{\xi(\lambda_1+\lambda_2)}{2}, \frac{\lambda_1+\lambda_2}{2}, \xi \lambda_1 \lambda_2\right), \quad  (\lambda_1, \lambda_2, \xi) \in \mathbb C^3.
\]
The map $\Psi_{\mathbb E}$ parametrizes the tetrablock in the sense that $\widehat{\Psi_{\mathbb E}(\overline{\mathbb D}^{3})}=\overline{\mathbb E}$ and $\text{int}(\widehat{\Psi_{\mathbb E}(\overline{\mathbb D}^{3})})=\mathbb E$. We also prove that $\Psi_{\mathbb E}(\mathbb D^3)$ is a uniqueness set for
$\operatorname{Hol}(\mathbb E)$.

\smallskip  

The parametrization map $\Psi_{\mathbb H}$ and the fact that $\Psi_{\mathbb H}(\mathbb T^4)=b\mathbb H$ allows us to transfer normalized Haar measure on $\mathbb T^4$ to a canonical probability measure $\mu_{\mathbb H}$ on $b\mathbb H$. Consequently, we define in Section \ref{sec003_hex} the Hardy space on the hexablock $\mathbb H$ as
\[
H^2(\mathbb H)
=
\overline{\mathbb C[a,x_1,x_2,x_3]}
^{L^2(b\mathbb H,\mu_{\mathbb H})}.
\]
The motivation behind considering the measure $\mu_{\mathbb H}$ is discussed in Section \ref{sec003_hex}. The construction of Hardy space $H^2(\mathbb H)$ capitalizes on the geometry of $\mathbb H$ and provides the appropriate setting for studying its coordinate multiplication operators. We present two decomposition structures for $H^2(\mathbb H)$ and then show that the space $H^2(\mathbb H)$ is a reproducing kernel Hilbert space. We further identify it with a closed subspace of $L^2(b\mathbb H, \mu_{\mathbb H})$ via radial boundary maps
\[
\delta_r(a, x_1, x_2, x_3)=(ra, rx_1, rx_2, r^2x_3) \quad (0 \leq r<1).
\]
We study in Section \ref{sec004_hex} the relationship between $H^2(\mathbb H)$ and the classical Hardy space on the polydisc $\mathbb D^4$ through composition with $\Psi_{\mathbb H}$. This provides a way to establish the basic function-theoretic framework required for the subsequent operator-theoretic developments on $H^2(\mathbb H)$. In particular, we identify the multiplier algebra of $H^2(\mathbb H)$ and show that it is isometrically isomorphic to $H^\infty(\mathbb H)$. It is also proved in Section \ref{sec004_hex} that there is an isometric algebra embedding of $H^\infty(\mathbb H)$ into $L^\infty(b\mathbb H, \mu_{\mathbb H})$. Next, we study in Section \ref{sec005_hex} the Toeplitz operators on $H^2(\mathbb H)$. Given $\varphi\in L^\infty(b\mathbb H,\mu_{\mathbb H})$, we define
\[
T_\varphi
=
P_{H^2(\mathbb H)}
M_\varphi\big|_{H^2(\mathbb H)},
\]
where $M_\varphi$ denotes multiplication by $\varphi$ on $L^2(b\mathbb H,\mu_{\mathbb H})$ and $P_{H^2(\mathbb H)}$ is the orthogonal projection onto $H^2(\mathbb H)$. The map $T_\varphi$ is referred to as the Toeplitz operator with symbol $\varphi$. We show that the coordinate multiplication operators $(T_a,T_{x_1},T_{x_2},T_{x_3})$
form an $\mathbb H$-isometry on $H^2(\mathbb H)$ and the corresponding multiplication operators $(M_a,M_{x_1},M_{x_2},M_{x_3})$ on $L^2(b\mathbb H,\mu_{\mathbb H})$ is its minimal normal extension that itself is an $\mathbb H$-unitary. The fundamental operator identities associated with $\mathbb H$-isometries from \cite{PalTomarI} suggest the appropriate analogues of the classical Brown-Halmos relations. In particular, the $\mathbb H$-unitary $(M_a, M_{x_1}, M_{x_2}, M_{x_3})$ satisfies (see Theorem 3.5 in \cite{PalTomarI}) the following identities:
\[
M_a^*M_a+M_{x_1}^*M_{x_1}=I, \quad M_a^*M_a+M_{x_2}^*M_{x_2}=I \quad 
M_{x_1}=M_{x_2}^*M_{x_3},\quad 
M_{x_3}^*M_{x_3}=I.
\]
Capitalizing on this, we prove the following relations for a Toeplitz operator $T$:
\begin{align*} 
	T_{x_3}^*T T_{x_3}
	=T,\quad
	T_{x_1}^*T T_{x_3}
	=T T_{x_2},\quad
	T_a^*T T_a+
	T_{x_2}^*T T_{x_2}
	=T,\quad
	T_a^*T T_a+
	T_{x_1}^*T T_{x_1}
	=T.
\end{align*}
We refer to these identities as the Brown-Halmos identities in the present setting. One of our main results is the following Brown-Halmos type theorem for the Toeplitz operators on $H^2(\mathbb H)$.

\begin{theorem}
	A bounded linear operator $T$ acting on $H^2(\mathbb H)$ is Toeplitz if and only if $T$ satisfies the Brown-Halmos identities. 
\end{theorem} 

The above result is presented as Theorem \ref{thm_BH} in Section \ref{sec005_hex}. Its proof relies on a commutant lifting result presented as Lemma \ref{lem:lifting-hex}. Several corollaries follow directly from this Brown-Halmos characterization. First, we show that every bounded operator on $H^2(\mathbb H)$ commuting with the four coordinate multiplication operators is necessarily Toeplitz. This provides a direct connection between the commutant of the canonical $\mathbb H$-isometry and the Toeplitz algebra associated with the hexablock. We further establish that the coordinate multipliers have no nontrivial joint reducing subspaces. In fact, we prove a stronger conclusion that every bounded operator that doubly commutes with the coordinate multipliers is a scalar multiple of the identity. Another consequence of the Brown-Halmos characterization is that the zero operator is the only compact Toeplitz operator on $H^2(\mathbb H)$. Finally, we formulate and prove in Theorem \ref{thm_CLT} a commutant lifting theorem for an arbitrary $\mathbb H$-isometry and its minimal $\mathbb H$-unitary extension. 

\begin{theorem}[Commutant lifting theorem]
	Let $(T_1,T_2,T_3,T_4)$ be an $\mathbb H$-isometry acting on a Hilbert space $\mathcal H$ and let $(N_1,N_2,N_3,N_4)$ be its minimal $\mathbb H$-unitary extension on a Hilbert space $\mathcal K$ containing $\mathcal H$. An operator $T\in\mathcal B(\mathcal H)$ satisfies the Brown-Halmos relations with respect to $(T_1,T_2,T_3,T_4)$ if and only if there exists a unique operator $N\in\mathcal B(\mathcal K)$ such that
	\[
	T=P_{\mathcal H}N|_{\mathcal H},\qquad
	\|N\|=\|T\|,\qquad
	NN_j=N_jN,\quad 1\leq j\leq 4.
	\]
\end{theorem}

\smallskip 

We introduce in Section \ref{sec006_hex} a weighted Hardy-type space on $\mathbb H$. To do so, we prove that the restriction of $\pi_{\mathbb H}$ from the space $U(2)$ of unitary matrices in $M_2(\mathbb C)$ to $b\mathbb H$ is in fact a homeomorphism. We then define a probability measure $\sigma_{\mathbb H}$ on $b\mathbb H$ as the push forward measure $\sigma_{\mathbb H}=(\pi_{\mathbb H})_*m_{U(2)}$, where $m_{U(2)}$ is the normalized Haar measure on $U(2)$. We prove that the two measures $\sigma_{\mathbb H}$ and $\mu_{\mathbb H}$ are mutually absolutely continuous. A weighted Hardy space is defined as 
\[
H^2_\sigma(\mathbb H)
=
\overline{\mathbb C[a,x_1,x_2,x_3]}
^{\,L^2(b\mathbb H,\sigma_{\mathbb H})}.
\]
Furthermore, we study the connections between the space $H^2(\mathbb H)$ and $H^2_\sigma(\mathbb H)$. En route, we present an alternative description of the measure $\sigma_{\mathbb H}$ as the pushforward measure induced by a homeomorphism between $b\mathbb H$ and $\mathbb T \times \partial \mathbb B_2$, where $\mathbb B_2$ is the Euclidean unit ball in $\mathbb C^2$.

\smallskip 

The parametrization $\Psi_{\mathbb H}$ also provides a useful mechanism for constructing operator tuples associated with the hexablock. In Section \ref{sec007_hex}, we investigate the relationship between polydisc $\mathbb D^4$ operator tuples and the operators associated with $\mathbb H$. We first establish an operator-theoretic counterpart of the boundary parametrization as follows: a commuting quadruple is an $\mathbb H$-unitary if and only if it admits a representation through $\Psi_{\mathbb H}$ evaluated at a commuting quadruple of unitaries. This characterization gives a simple way for producing $\mathbb H$-isometries from commuting quadruples of isometries. We also show that if a commuting quadruple of contractions has $\overline{\mathbb D}^{4}$ as a spectral set, then its image under $\Psi_{\mathbb H}$ is an $\mathbb H$-contraction. Moreover, whenever the original quadruple admits a commuting unitary dilation, the resulting $\mathbb H$-contraction admits a dilation to an $\mathbb H$-unitary. These results illustrate how the parametrization of the hexablock can be used to transfer spectral-set and dilation properties from $\mathbb D^4$ to $\mathbb H$. We prove analogous results in Section \ref{sec008_hex} for the operators associated with $\mathbb E$.

\smallskip

\noindent \textbf{Notations.} Throughout the paper, all operators are bounded linear maps acting on complex Hilbert spaces and $\mathcal B(\mathcal H)$ denotes the algebra of operators acting on a Hilbert space $\mathcal H$. We denote by $\delta_r(a, x_1, x_2, x_3)=(ra, rx_1, rx_2, r^2x_3)$ for $0 \leq r<1$, which we refer to as the radial boundary map. For a commuting tuple of operators $\underline{T}=(T_1, \dotsc, T_n)$ acting on a Hilbert space $\HS$, we denote by $\sigma_T(T_1, \dotsc, T_n)$ the Taylor joint spectrum of $(T_1, \dotsc, T_n)$.

\section{A uniqueness set of the hexablock and tetrablock}\label{sec002_hex}

\noindent In this section, we present alternative characterizations of the domains $\mathbb H$ and $\mathbb E$ via the parametrization maps denoted by $\Psi_{\mathbb H}$ and $\Psi_{\mathbb E}$, respectively. We show that the image of $\mathbb D^4$ under the map $\Psi_{\mathbb H}$ is a uniqueness set for $\text{Hol}(\mathbb H)$, that is, $f \in \text{Hol}(\mathbb H)$ with $f|_{\Psi_{\mathbb H}}=0$ implies $f=0$ on $\mathbb H$. Moreover, we show that an analogous result holds for $\mathbb E$ and the associated map $\Psi_{\mathbb E}$. 

\smallskip 

To do so, we recall from the literature that the distinguished boundary $b\Omega$ of a domain $\Omega \subseteq \mathbb C^n$ is the smallest closed subset of $\overline{\Omega}$ on which every function that is continuous on $\overline{\Omega}$ and holomorphic on $\Omega$ attains its maximum modulus. For example, the authors of \cite{AbouhajarWhiteYoung2007} and \cite{BiswasPalI} proved that 
\begin{align*}
	b\mathbb E&=\{(u_{11}, u_{22}, \det(U)): U=(u_{ij}) \in U(2)\} \quad \text{and} \\
	b\mathbb H&=\{(u_{21}, u_{11}, u_{22}, \det(U)): U=(u_{ij}) \in U(2) \},
\end{align*}
respectively, where $U(2)$ is the space of all unitary matrices in $M_2(\mathbb C)$. We study below the parametrization map $\Psi_{\mathbb H}$ associated with $\mathbb H$. 

\begin{theorem}\label{prop:Psi-interior}
	The polynomial map $\Psi_{\mathbb H}:\mathbb C^4\to\mathbb C^4$ given by
	\[
	\Psi_{\mathbb H}(\xi_1,\xi_2,\lambda_1,\lambda_2)
	=
	\left(
	\frac{\lambda_1(\xi_1-1)}{2},
	\frac{\xi_2(\xi_1+1)}{2},
	\frac{\lambda_1\lambda_2(\xi_1+1)}{2},
	\xi_1\xi_2\lambda_1\lambda_2
	\right)
	\]
	satisfies  $\Psi_{\mathbb H}(\overline{\mathbb D} \times \mathbb D^2 \times \overline{\mathbb D})\subseteq\mathbb H$ and $\Psi_{\mathbb H}(\mathbb T^4)=b\mathbb H$. Furthermore,
	\[
	\widehat{\Psi_{\mathbb H}(\overline{\mathbb D}^{\,4})}=\overline{\mathbb H}
	\quad \text{and} \quad \text{int}\left(\widehat{\Psi_{\mathbb H}(\overline{\mathbb D}^{\,4})}\right)=\mathbb H.
	\]
\end{theorem}

\begin{proof}
	Let $\pi_{\mathbb H}:M_2(\mathbb C)\to\mathbb C^4$ be given by
	$\pi_{\mathbb H}(A)=(a_{21},a_{11},a_{22},\det A)$, where
	$A=(a_{ij})_{i,j=1}^2 \in M_2(\mathbb C)$. We use the fact from \cite{BiswasPalI} that the normed hexablock $\mathbb H_N$ is contained in $\mathbb H$, that is, 
	\[
	\mathbb H_N=\pi_{\mathbb H}\big(\{A\in M_2(\mathbb C):\|A\|<1\}\big)
	\subseteq\mathbb H
	\]
	and the distinguished-boundary of $\mathbb H$ is
	$b\mathbb H=\pi_{\mathbb H}(U(2))$. Let $(\xi_1,\xi_2,\lambda_1,\lambda_2)\in \overline{\mathbb D} \times \mathbb D^2 \times \overline{\mathbb D}$.
	Consider
	\[
	A=
	\frac12
	\begin{pmatrix}
		\xi_2(\xi_1+1) & \xi_2\lambda_2(\xi_1-1)\\
		\lambda_1(\xi_1-1) & \lambda_1\lambda_2(\xi_1+1)
	\end{pmatrix} \quad \text{and} \quad 
	H_0=\frac1{\sqrt2}
	\begin{pmatrix}
		1&1\\
		1&-1
	\end{pmatrix}.
	\]
	It is not difficult to see that 
	\begin{equation}\label{eqn_factor_A}
		A=
		\begin{pmatrix}
			\xi_2&0\\
			0&\lambda_1
		\end{pmatrix}
		H_0
		\begin{pmatrix}
			\xi_1&0\\
			0&1
		\end{pmatrix}
		H_0
		\begin{pmatrix}
			1&0\\
			0&\lambda_2
		\end{pmatrix}.
	\end{equation}
	Also, $\det A=\xi_1\xi_2\lambda_1\lambda_2$ and $\pi_{\mathbb H}(A)=\Psi_{\mathbb H}(\xi_1,\xi_2,\lambda_1,\lambda_2)$. As $H_0$ is unitary, we have 
	\[
	\|A\|
	\leq
	\max\{|\xi_2|,|\lambda_1|\} \cdot 
	\max\{|\xi_1|,1\} \cdot 
	\max\{1,|\lambda_2|\}
	=
	\max\{|\xi_2|,|\lambda_1|\}<1.
	\]
	Consequently, $\Psi_{\mathbb H}(\overline{\mathbb D} \times \mathbb D^2 \times \overline{\mathbb D})\subseteq\mathbb H$. Next, assume that $(\xi_1,\xi_2,\lambda_1,\lambda_2)\in\mathbb T^4$. Then
	every factor in the factorization \eqref{eqn_factor_A} of $A$ is unitary.
	Thus, $A$ is a unitary matrix and so, $\Psi_{\mathbb H}(\mathbb T^4)\subseteq b\mathbb H$. To prove the reverse inclusion, let
	$(a,x_1,x_2,x_3)\in b\mathbb H$. Then there exists a unitary matrix
	\[
	U=
	\begin{pmatrix}
		\alpha&\beta\\
		\gamma&\delta
	\end{pmatrix}
	\]
	such that $(a,x_1,x_2,x_3)
	=
	(\gamma,\alpha,\delta,\det U).
	$
	Since $U$ is unitary, we have that 
	\[
	|\alpha|=|\delta|,
	\quad 
	|\beta|=|\gamma|,
	\quad |\alpha|^2+|\gamma|^2=1
	\quad \text{and} \quad \delta=\overline{\alpha}\det U.
	\]
	First consider the case when $0<|\alpha|<1$. Choose
	$\xi_1\in\mathbb T\setminus\{-1,1\}$ such that
	\[
	\left|\frac{\xi_1+1}{2}\right|=|\alpha|.
	\]
	Since $|\xi_1|=1$, we have
	\[
	\left|\frac{\xi_1-1}{2}\right|
	=
	\sqrt{1-|\alpha|^2}
	=
	|\gamma|.
	\]
	Define
	\[
	\xi_2=\frac{2\alpha}{\xi_1+1},\qquad
	\lambda_1=\frac{2\gamma}{\xi_1-1},\qquad
	\lambda_2=\frac{2\delta}{\lambda_1(\xi_1+1)}.
	\]
	Then $\xi_2,\lambda_1,\lambda_2\in\mathbb T$ and
	\[
	\frac{\lambda_1(\xi_1-1)}2=\gamma,\qquad
	\frac{\xi_2(\xi_1+1)}2=\alpha,\qquad
	\frac{\lambda_1\lambda_2(\xi_1+1)}2=\delta.
	\]
	Using $|\xi_1|=1$ and $\delta=\overline{\alpha}\det U$, we have
	\begin{align*}
		\xi_1\xi_2\lambda_1\lambda_2
		=
		\frac{4\xi_1\alpha\delta}{(\xi_1+1)^2}
		=
		\frac{4\alpha\delta}{|\xi_1+1|^2}
		=
		\frac{\alpha\delta}{|\alpha|^2}
		=\det U.
	\end{align*}
	Therefore, $\Psi_{\mathbb H}(\xi_1,\xi_2,\lambda_1,\lambda_2)=(\gamma,\alpha,\delta,\det U)$. Next, assume that $|\alpha|=1$. Since $U$ is unitary, it follows that $\beta=\gamma=0$ and
	$\det U=\alpha\delta$. Taking $(\xi_1,\xi_2,\lambda_1,\lambda_2)
	=
	(1,\alpha,1,\delta)$ gives that
	\[
	\Psi_{\mathbb H}(\xi_1,\xi_2,\lambda_1,\lambda_2)
	=
	(0,\alpha,\delta,\alpha\delta)
	=
	(\gamma,\alpha,\delta,\det U).
	\]
	Finally, take $\alpha=0$. Then $\delta=0$
	and $|\gamma|=1$. Taking $(\xi_1,\xi_2,\lambda_1,\lambda_2)
	=
	\left(-1, \overline{\gamma}\det U,-\gamma,1\right)$,
	we have that
	$
	\Psi_{\mathbb H}(\xi_1,\xi_2,\lambda_1,\lambda_2)
	=
	(\gamma,0,0,\det U).
	$
	In all cases, $(a,x_1,x_2,x_3)\in\Psi_{\mathbb H}(\mathbb T^4)$ and thus
	$\Psi_{\mathbb H}(\mathbb T^4)=b\mathbb H$. Since $\overline{\mathbb H}$ is
	polynomially convex, it follows that $\widehat{\Psi_{\mathbb H}(\overline{\mathbb D}^{\,4})} \subseteq \overline{\mathbb H}$. Let $z \in \overline{\mathbb H}$ and let $q$ be a holomorphic polynomial in four variables. Then 
	\[
	|q(z)| \leq \|q\|_{\infty, \overline{\mathbb H}}=\|q\|_{\infty, b \mathbb H}=\|q\|_{\infty, \Psi_{\mathbb H}(\mathbb T^4)} \leq \|q\|_{\infty, \Psi_{\mathbb H}(\overline{\mathbb D}^4)} \quad \text{and so,}
	\quad z \in \widehat{\Psi_{\mathbb H}(\overline{\mathbb D}^4)}.
	\]
	Consequently, $\overline{\mathbb H}
	=
	\widehat{\Psi_{\mathbb H}(\overline{\mathbb D}^{4})}$ and the remaining  conclusion now follows.
\end{proof}

Next, we recall the notion of \emph{uniqueness sets} from \cite{Chee1978} and show that the image of $\mathbb D^4$ under the map $\Psi_{\mathbb H}$ is a uniqueness set for $\mathrm{Hol}(\mathbb H)$. An analogous result has been established for the class of holomorphic functions on the pentablock in \cite{SuTomar}.

\begin{definition}
	Let $\Omega\subseteq\mathbb C^n$ be a domain and let
	$\Lambda$ be a family of functions on $\Omega$. A set
	$\Omega' \subseteq\Omega$ is called a \emph{uniqueness set for $\Lambda$}
	if $f\in \Lambda$ and $f|_{\Omega'}=0$ implies $f\equiv0$ on $\Omega$. 
\end{definition}

\begin{theorem} \label{thm:Psi-uniqueness}
	The set $\Psi_{\mathbb H}(\mathbb D^4)$ is a uniqueness set for
	$\operatorname{Hol}(\mathbb H)$.
\end{theorem}

\begin{proof}
	Take $f \in \text{Hol}(\mathbb H)$ and let $f=0$ on $\Psi_{\mathbb H}(\mathbb D^4)$. We have by Theorem \ref{prop:Psi-interior} that $\Psi_{\mathbb H}(\overline{\mathbb{D}}\times\mathbb{D}^2 \times \overline{\mathbb D})\subseteq\mathbb{H}$. A simple continuity argument shows that $f=0$ on $\Psi_{\mathbb H}(\overline{\mathbb D} \times \mathbb D^2 \times \overline{\mathbb D})$. We now use the $(1, 1, 1, 2)$-quasi-balanced property of $\mathbb H$ and $\overline{\mathbb H}$ to show that $f=0$ on $\mathbb H$. For $0<s<1$, define $f_s=f\circ \delta_s$, where $\delta_s(a, x_1, x_2, x_3)=(sa, sx_1, sx_2, s^2x_3)$ for $(a, x_1, x_2, x_3) \in \overline{\mathbb H}$. Evidently, each $f_s$ is holomorphic on $\overline{\mathbb H}$. Fix $s \in (0, 1)$. Let $(a, x_1, x_2, x_3) \in b \mathbb H$. By Theorem \ref{prop:Psi-interior}, there exists $(\xi_1, \xi_2, \lambda_1, \lambda_2) \in \mathbb T^4$ such that $\Psi_{\mathbb H}(\xi_1, \xi_2, \lambda
	_1, \lambda_2)=(a, x_1, x_2, x_3)$. Since $f|_{\Psi_{\mathbb H}(\overline{\mathbb{D}}\times\mathbb{D}^2 \times \overline{\mathbb D})}=0$, we have 
	\[
	f_s(a, x_1, x_2, x_3)=f(sa, sx_1, sx_2, s^2x_3)=f\left(\Psi_{\mathbb H}(\xi_1, s\xi_2, s\lambda_1, \lambda_2,\right))=0
	\]
	and so, $f_s|_{b \mathbb H}=0$. Therefore, $f_s=0$ on $\overline{\mathbb{H}}$ for every $s \in (0, 1)$ and thus $f=0$ on $\mathbb H_s$ for $0<s<1$, where $\mathbb H_s=\{(sa, sx_1, sx_2, s^2x_3): (a, x_1, x_2, x_3) \in \overline{\mathbb{H}}\}$. It was proved in  \cite{BiswasPalI} that $\mathbb H=\cup_{0<s<1}\mathbb H_s$ and so, $f=0$ on $\mathbb H$. The proof is now complete.
\end{proof}

Theorem \ref{prop:Psi-interior} shows that the set $\Psi_{\mathbb H}(\mathbb D^4)$ completely determines the domain $\mathbb H$. Moreover, the
values of a holomorphic function on $\mathbb H$ on this set determine the
function on $\mathbb H$. In particular, we have by Theorem \ref{thm:Psi-uniqueness} that if
$f,g\in\operatorname{Hol}(\mathbb H)$ and
\[
f=g \quad \text{on } \Psi_{\mathbb H}(\mathbb D^4),
\]
then $f=g$ on $\mathbb H$. We shall in later sections that the parametrization map $\Psi_{\mathbb H}$ plays a pivotal role in the construction of the Hardy space $H^2(\mathbb H)$ associated with the hexablock. Next, we present an analogous result for the domain $\mathbb E$ and obtain a uniqueness set of $\mathbb E$, which is a proper subset of $\mathbb E$. 

\begin{theorem}
	The polynomial map $\Psi_{\mathbb E}: \mathbb C^3 \to \mathbb C^3$ defined as
	\[
	\Psi_{\mathbb E}(\lambda_1, \lambda_2, \xi)=\left(\frac{\xi(\lambda_1+\lambda_2)}{2}, \frac{\lambda_1+\lambda_2}{2}, \xi \lambda_1 \lambda_2\right) 
	\]
	satisfies $\Psi_{\mathbb E}(\mathbb D^2 \times \overline{\mathbb D})\subseteq\mathbb E$ and $\Psi_{\mathbb E}(\mathbb T^3)=b\mathbb E$. Also,
	\[
	\widehat{\Psi_{\mathbb E}(\overline{\mathbb D}^{3})}=\overline{\mathbb E}
	\quad \text{and} \quad \text{int}\left(\widehat{\Psi_{\mathbb E}(\overline{\mathbb D}^{3})}\right)=\mathbb E.
	\]
	Furthermore, $\Psi_{\mathbb E}(\mathbb D^3)$ is a uniqueness set for
	$\operatorname{Hol}(\mathbb E)$.
\end{theorem}

\begin{proof}
	Let $\pi_{\mathbb E}:M_2(\mathbb C)\to\mathbb C^3$ be given by
	$\pi_{\mathbb E}(A)=(a_{11},a_{22},\det A)$, where
	$A=(a_{ij})_{i,j=1}^2 \in M_2(\mathbb C)$. It was proved in \cite{AbouhajarWhiteYoung2007} that 
	$\mathbb E=\pi_{\mathbb E}\big(\{A\in M_2(\mathbb C):\|A\|<1\}\big)$. Let $(\lambda_1,\lambda_2, \xi)\in\mathbb D^2 \times \overline{\mathbb D}$.
	Consider
	\[
	A=
	\frac12
	\begin{pmatrix}
		\xi(\lambda_1+\lambda_2) & \xi(\lambda_1-\lambda_2)\\
		\lambda_1-\lambda_2 & \lambda_1+\lambda_2
	\end{pmatrix} \quad \text{and} \quad 
	H_0=\frac1{\sqrt2}
	\begin{pmatrix}
		1&1\\
		1&-1
	\end{pmatrix}.
	\]
	Note that $\det A=\xi\lambda_1\lambda_2$ and $\pi_{\mathbb E}(A)=\Psi_{\mathbb E}(\lambda_1,\lambda_2, \xi)$. Also, we have 
	\begin{equation}\label{eqn_factor_A_2}
		A=
		\begin{pmatrix}
			\xi&0\\
			0&1
		\end{pmatrix}
		H_0
		\begin{pmatrix}
			\lambda_1&0\\
			0&\lambda_2
		\end{pmatrix}
		H_0.
	\end{equation}
	Then $\|A\| \leq \max\{|\lambda_1|,|\lambda_2|\}<1$. Consequently, $\Psi_{\mathbb E}(\mathbb D^2 \times \overline{\mathbb D})\subseteq\mathbb E$. Next, assume that $(\lambda_1,\lambda_2, \xi)\in\mathbb T^3$. It follows from \eqref{eqn_factor_A_2} that $A$ is unitary and so, $\Psi_{\mathbb E}(\mathbb T^3)\subseteq b\mathbb E$. Next, suppose
	$(x_1,x_2,x_3)\in b\mathbb E$. Since $\overline{\mathbb D}=\{(z_1+z_2)\slash 2: z_1, z_2 \in \mathbb T\}$ and $x_2 \in \overline{\mathbb D}$, one can choose $\lambda_1, \lambda_2 \in \mathbb T$ such that $x_2=(\lambda_1+\lambda_2)\slash 2$. Set $\xi=\overline{\lambda}_1 \overline{\lambda}_2 x_3$. Since $(x_1, x_2, x_3) \in b\mathbb E$, we have that $x_1=\overline{x}_2x_3$ and thus 
	\[
	x_1=\frac{(\overline{\lambda}_1+\overline{\lambda}_2)\xi \lambda _1 \lambda_2}{2}=\frac{\xi(\lambda_1+\lambda_2)}{2}. 
	\]
	Hence, $(x_1, x_2, x_3)=\Psi_{\mathbb E}(\lambda_1, \lambda_2, \xi)$ and $\xi, \lambda_1, \lambda_2 \in \mathbb T$. Hence, $\Psi_{\mathbb E}(\mathbb T^3)=b\mathbb E$. Next, the polynomial convexity of $\overline{\mathbb E}$ gives that $\widehat{\Psi_{\mathbb E}(\overline{\mathbb D}^{3})} \subseteq \overline{\mathbb E}$. Let $x \in \overline{\mathbb E}$ and let $q$ be a holomorphic polynomial in three variables. Then 
	\[
	|q(x)| \leq \|q\|_{\infty, \overline{\mathbb E}}=\|q\|_{\infty, b \mathbb E}=\|q\|_{\infty, \Psi_{\mathbb E}(\mathbb T^3)} \leq \|q\|_{\infty, \Psi_{\mathbb E}(\overline{\mathbb D}^3)} \quad \text{and so,}
	\quad x \in \widehat{\Psi_{\mathbb E}(\overline{\mathbb D}^3)}.
	\]
	Consequently, $\overline{\mathbb E}
	=
	\widehat{\Psi_{\mathbb E}(\overline{\mathbb D}^{3})}$ and thus $\mathbb E=\mathrm{int}(\widehat{\Psi_{\mathbb E}(\overline{\mathbb D}^{3})})$.
	
	\smallskip 
	
	Next, we show that $\Psi_{\mathbb E}(\mathbb D^3)$ is a uniqueness set for $\operatorname{Hol}(\mathbb E)$. Take $f \in \text{Hol}(\mathbb E)$ such that $f=0$ on $\Psi_{\mathbb E}(\mathbb D^3)$. We have proved above that $\Psi_{\mathbb E}(\mathbb{D}^2 \times \overline{\mathbb D})\subseteq\mathbb{H}$. By continuity arguments, $f=0$ on $\Psi_{\mathbb E}(\mathbb D^2 \times \overline{\mathbb D})$. Since $\mathbb E$ and $\overline{\mathbb E}$ are $(1, 1, 2)$-quasi-balanced sets, one can define $f_s=f\circ \theta_s$, where $\theta_s(x_1, x_2, x_3)=(sx_1, sx_2, s^2x_3)$ for $(x_1, x_2, x_3) \in \overline{\mathbb E}$ and $s \in (0, 1)$. Clearly, each $f_s$ is holomorphic on $\overline{\mathbb E}$. Fix $s \in (0, 1)$. Let $(x_1, x_2, x_3) \in b \mathbb E$. Then there exists $(\lambda_1, \lambda_2, \xi) \in \mathbb T^3$ such that $\Psi_{\mathbb E}(\lambda
	_1, \lambda_2, \xi)=(x_1, x_2, x_3)$. Since $f|_{\Psi_{\mathbb E}(\mathbb{D}^2 \times \overline{\mathbb D})}=0$, we have 
	\[
	f_s(x_1, x_2, x_3)=f(sx_1, sx_2, s^2x_3)=f\left(\Psi_{\mathbb E}(s\lambda_1, s\lambda_2, \xi \right))=0
	\]
	and so, $f_s|_{b \mathbb E}=0$. Therefore, $f_s=0$ on $\overline{\mathbb{E}}$ for every $s \in (0, 1)$ and thus $f=0$ on $\mathbb E_s$ for $0<s<1$, where $\mathbb E_s=\{(sx_1, sx_2, s^2x_3): (x_1, x_2, x_3) \in \overline{\mathbb{E}}\}$. Since $\mathbb E=\cup_{0<s<1}\mathbb E_s$, we have that $f=0$ on $\mathbb E$, which completes the proof.
\end{proof}

\section{The Hardy space on the hexablock}\label{sec003_hex}

\noindent In this section, we introduce the Hardy space associated with the hexablock. We begin by recalling from literature the constructions of Hardy spaces on the symmetrized bidisc $\mathbb G_2$, the tetrablock $\mathbb E$ and the pentablock $\mathbb P$. Recall that the symmetrized bidisc is given by
$\mathbb G_2=\pi_{\mathrm{symm}}(\mathbb D^2)$, where
\[
\pi_{\mathrm{symm}}:\mathbb C^2\to\mathbb C^2, \quad \pi_{\mathrm{symm}}(z_1,z_2)=(z_1+z_2,z_1z_2).
\]
is the symmetrization map. The closure and
distinguished boundary of $\mathbb G_2$ are given by
$\Gamma_2=\pi_{\mathrm{symm}}(\overline{\mathbb D}^2)$ and
$b\Gamma_2=\pi_{\mathrm{symm}}(\mathbb T^2)$, respectively. The map $\pi_{\mathrm{symm}}$ is a proper holomorphic map of multiplicity two with $J_{\pi_{\mathrm{symm}}}(z_1,z_2)=z_1-z_2$. The authors of \cite{MisraShyamRoyZhang2013} introduced the Hardy space
$H^2(\mathbb G_2)$ as the space consisting of all holomorphic maps $g$ on $\mathbb G_2$ such that
\[
\|g\|_{H^2(\mathbb G_2)}^2
=
\frac{1}{2}
\sup_{0<r<1}
\int_{\mathbb T^2}
\left|g\circ\pi_{\mathrm{symm},r}(\lambda_1, \lambda_2)
\right|^2
\left|J_{\pi_{\mathrm{symm}}}(r\lambda_1,r\lambda_2)\right|^2
\,dm_2 < \infty,
\]
where $\pi_{\mathrm{symm},r}(\lambda_1,\lambda_2)
=
\left(r(\lambda_1+\lambda_2),r^2\lambda_1\lambda_2\right)$ and $m_2$ is the normalized Haar measure on $\mathbb T^2$. A similar construction from \cite{BeraChavanJain2025} applies to the tetrablock $\mathbb E$.  Let
\[
\mathcal R_{II}
=
\left\{
(z_1,z_2,z_3)\in\mathbb C^3:
\left\|
\begin{pmatrix}
	z_1&z_3\\
	z_3&z_2
\end{pmatrix}
\right\|<1
\right\}
\]
and consider $\varphi:\mathcal R_{II}\to \mathbb E$ defined as $
\varphi(z_1,z_2,z_3)
=
(z_1,z_2,z_1z_2-z_3^2)$. Moreover, $J_\varphi(z_1,z_2,z_3) =-2z_3$. The corresponding Hardy space is given by
\[
H^2(\mathbb E)
=
\left\{
f\in\operatorname{Hol}(\mathbb E):
\|J_\varphi f\circ\varphi\|_{H^2(\mathcal R_{II})} <\infty 
\right\}.
\]
Hence, in both cases, the Hardy space on the target domain is obtained by pulling back to a domain of the same complex dimension and incorporating the
complex Jacobian of the proper holomorphic map. However, it was shown in \cite{SuTomar} that the construction of the Hardy space on the domain $\mathbb P$ does not follow an analogous framework. The primary reason for such an anomaly is the fact that the natural matrix realization of $\mathbb P$ given by
\[
\pi_\mathbb P: \left\{A=[a_{ij}] \in M_2(\mathbb C): \|A\|<1  \right\} \to \mathbb C^3,
\qquad
\pi_{\mathbb P}(A)=\bigl(a_{21}, a_{11}+ a_{22},\det(A)\bigr)
\]
is neither proper nor does it preserve the distinguished boundaries in a surjective manner. To overcome this, the authors of \cite{SuTomar} consider the polynomial map 
\[
\Psi_{\mathbb P}:\mathbb D^3\to\mathbb C^3, \quad 
(\xi,\lambda_1,\lambda_2)
\mapsto
\left(
\frac{\xi(\lambda_1-\lambda_2)}{2},
\lambda_1+\lambda_2,
\lambda_1\lambda_2
\right).
\]
The map $\Psi_{\mathbb P}$ satisfies $\Psi_{\mathbb P}(\mathbb T^3)=b\mathbb P$. Let $m_3$ denote the normalized Haar measure on $\mathbb T^3$. Let $\mu_{\mathbb P}$ be the probability measure on $b\mathbb P$ given by $\mu_{\mathbb P}=(\Psi_{\mathbb P})_*m_3$. Then
\[
\int_{b\mathbb{P}}F(a,s,p)\,\mathrm{d}\mu_{\mathbb P}
=
\int_{\mathbb{T}^3}
F\!\left(
\frac{\xi(\lambda_1-\lambda_2)}{2},
\lambda_1+\lambda_2,
\lambda_1\lambda_2
\right)
\,\mathrm{d}m_3(\xi,\lambda_1,\lambda_2)
\]
for every $F\in C(b\mathbb{P})$. The Hardy space of the pentablock is defined as
\[
H^2(\mathbb{P})
=
\overline{\mathbb{C}[a,s,p]}^{\,L^2(b\mathbb{P},\mu)}.
\]
Our construction of the Hardy space $H^2(\mathbb H)$ follows a similar approach as that used for $H^2(\mathbb P)$.  Recall from Theorem \ref{prop:Psi-interior} that $b\mathbb H$ admits the parameterization given by
\[
b \mathbb H=
\left\{\left(
\frac{\lambda_1(\xi_1-1)}{2},
\frac{\xi_2(\xi_1+1)}{2},
\frac{\lambda_1\lambda_2(\xi_1+1)}{2},
\xi_1\xi_2\lambda_1\lambda_2
\right):
(\xi_1, \xi_2,\lambda_1,\lambda_2)\in\mathbb{T}^4\right\}.
\]
Let $m_k$ denote normalized Haar measure on the $k$-torus $\mathbb{T}^k$ and define the pushforward measure $\mu_{\mathbb H}=(\Psi_\mathbb H)_*m_4$.
Equivalently, for every \(F\in\mathcal{C}(b\mathbb{H})\),
\[
\int_{b\mathbb{H}}F(a, x_1, x_2, x_3)\,\mathrm{d}\mu_{\mathbb H}
=
\int_{\mathbb{T}^4}
F\left(\Psi_{\mathbb H}(\xi_1, \xi_2,\lambda_1,\lambda_2)\right)
\mathrm{d}m_4(\xi_1, \xi_2,\lambda_1,\lambda_2).
\]
Thus, the map $\Psi_{\mathbb H}$ gives a canonical boundary measure on $b\mathbb H$. One can further study the structure of the measure $\mu_{\mathbb H}$ more explicitly. To do so, we first describe the
fibres of the natural projection from $b\mathbb{H}$ onto $b\mathbb E$. Define
$ \operatorname{pr}:b\mathbb{H}\to b\mathbb E$ as
$\operatorname{pr}(a, x_1, x_2, x_3)=(x_1, x_2, x_3)$. For $(x_1, x_2, x_3)\in b\mathbb E$, put
\[
\rho(x_1, x_2, x_3)
=
\sqrt{1-|x_1|^2},
\]
which is same as $\sqrt{1-|x_2|^2}$ since $|x_1|=|x_2|$ for $(x_1, x_2, x_3) \in b\mathbb E$. As $|a|^2+|x_1|^2=1$
on $b\mathbb{H}$, the fibre of $\operatorname{pr}$ over $(x_1, x_2, x_3) \in b\mathbb H$ is given by
\begin{align*}
	\operatorname{pr}^{-1}(x_1, x_2, x_3)
	=
	\left\{
	(a,x_1, x_2, x_3)\in b\mathbb{H}: |a|=\rho(x_1, x_2, x_3)
	\right\}
	 =
	\left\{
	\bigl(\xi\rho(x_1, x_2, x_3), x_1, x_2, x_3\bigr):
	\xi\in\mathbb{T}
	\right\}.
\end{align*}
Let $\nu_{\mathbb H}=(\operatorname{pr})_*\mu_{\mathbb H}$ denote the pushforward of $\mu_{\mathbb H}$ under the projection
$\operatorname{pr}:b\mathbb H\to b\mathbb E$. Thus, $\nu_{\mathbb H}$ is a probability measure on $b\mathbb E$. By the above description of the fibres of $\operatorname{pr}$, for every $F\in C(b\mathbb H)$, we have
\[
\int_{b\mathbb H}F\,d\mu_{\mathbb H}
=
\int_{b\mathbb E}\int_{\mathbb T}
F\bigl(\xi\rho(x_1,x_2,x_3),x_1,x_2,x_3\bigr)
\,dm_1(\xi)\,d\nu_{\mathbb H}.
\]
On each fibre over $(x_1,x_2,x_3)\in b\mathbb E$, the measure induced by $\mu_{\mathbb H}$ is the normalized rotationally invariant measure obtained from Haar measure on $\mathbb T$ under the parametrization
\[
\xi\mapsto \bigl(\xi\rho(x_1,x_2,x_3),x_1,x_2,x_3\bigr).
\]
The above description gives the invariance of $\mu_{\mathbb H}$ under rotations of the $a$-fibres. For $\eta\in\mathbb T$, define
\[
R_\eta(a,x_1,x_2,x_3)
=(\eta a,x_1,x_2,x_3).
\]
For a point on the fibre over $(x_1,x_2,x_3)$, we have
\[
\begin{aligned}
	R_\eta\bigl(\xi\rho(x_1,x_2,x_3),x_1,x_2,x_3\bigr)
	=
	\bigl((\eta\xi)\rho(x_1,x_2,x_3),x_1,x_2,x_3\bigr).
\end{aligned}
\]
Thus, $R_\eta$ fixes the last three coordinates and rotates the first coordinate by $\eta$. Since $m_1$ is normalized Haar
measure on $\mathbb{T}$,
\[
\int_{\mathbb{T}} G(\eta\xi)\,dm_1(\xi)
=
\int_{\mathbb{T}} G(\xi)\,dm_1(\xi)
\]
for every $G\in L^1(\mathbb{T})$. Then
\begin{align*}
	\int_{b\mathbb H}F\circ R_\eta\,d\mu_{\mathbb H}
	&=
	\int_{b\mathbb E}\int_{\mathbb T}
	F\bigl(\eta\xi\rho(x_1,x_2,x_3),x_1,x_2,x_3\bigr)
	\,dm_1(\xi)\,d\nu_{\mathbb H}=
	\int_{b\mathbb H}F\,d\mu_{\mathbb H}.
\end{align*}
and thus $(R_\eta)_*\mu_{\mathbb H}=\mu_{\mathbb H}$ for every $\eta\in\mathbb T$. 

\begin{remark}
	A primary motivation for our choice of the measure $\mu_{\mathbb H}$
	comes from the construction of the Hardy space on the symmetrized
	bidisc \cite{MisraShyamRoyZhang2013}. Recall that the symmetrized bidisc is the image of the bidisc
	under the symmetrization map
	\[
	\pi_{\mathrm{symm}}:\mathbb D^2\to\mathbb G_2,
	\qquad
	\pi_{\mathrm{symm}}(z_1,z_2)=(z_1+z_2,z_1z_2).
	\]
	This function sends $\mathbb D^2$ onto $\mathbb G_2$ and it 
	maps the distinguished boundary $\mathbb T^2$ onto the distinguished boundary $b\Gamma_2$ of $\mathbb G_2$. Thus, the symmetrization map provides a natural parametrization of the symmetrized bidisc in terms of the bidisc, including its distinguished boundary. The polynomial map $\Psi_{\mathbb H}$ plays an analogous role for the hexablock.
	Indeed, we have
	\[
	\Psi_{\mathbb H}(\mathbb T^4)=b\mathbb H, \quad  \widehat{\Psi_{\mathbb H}(\overline{\mathbb D}^{\,4})}
	=\overline{\mathbb H}
	\quad \text{and} \quad 
	\operatorname{int}\left(
	\widehat{\Psi_{\mathbb H}(\overline{\mathbb D}^{\,4})}
	\right)=\mathbb H.
	\]
	Consequently, just as the symmetrization map relates the geometry
	of the bidisc to that of the symmetrized bidisc, the map
	$\Psi_{\mathbb H}$ relates the geometry of the polydisc $\mathbb D^4$ to that
	of the hexablock. This motivates us to use the pushforward
	of normalized Haar measure on $\mathbb T^4$, namely
	\[
	\mu_{\mathbb H}=(\Psi_{\mathbb H})_*m_4,
	\]
	as a natural probability measure on $b\mathbb H$. Another motivation comes from the construction of the Hardy
	space on the pentablock by the authors of \cite{SuTomar},
	where a boundary measure is obtained similarly as the pushforward
	of normalized Haar measure on $\mathbb T^3$ under a polynomial parametrization of the distinguished boundary. Following this approach, we
	define the Hardy space of the hexablock using the measure
	$\mu_{\mathbb H}$. We note that, unlike the construction on
	the symmetrized bidisc, our definition is based directly on the
	boundary parametrization and does not involve a Jacobian weight
	arising from a proper holomorphic mapping. \qed
\end{remark}

We next consider the $(1,1,1,2)$-quasi-circular action on the domain $\mathbb H$ given by
\[
q_\omega(a,x_1,x_2,x_3)
=
(\omega a,\omega x_1,\omega x_2,\omega^2x_3),
\qquad \omega\in\mathbb T.
\]
Since $\rho(\omega x_1,\omega x_2,\omega^2x_3)
=
\rho(x_1,x_2,x_3)$, we have that 
\[
q_\omega\bigl(\xi\rho(x_1,x_2,x_3),x_1,x_2,x_3\bigr)=
\bigl((\omega\xi)\rho(\omega x_1,\omega x_2,\omega^2x_3),
\omega x_1,\omega x_2,\omega^2x_3\bigr).
\] 
The invariance of $\mu_{\mathbb H}$ under $q_\omega$ can be observed directly from the parametrization map $\Psi_{\mathbb H}$. Indeed, we have
\begin{align*}
	\Psi_{\mathbb H}(\xi_1,\omega\xi_2,\omega\lambda_1,\lambda_2)
	&=
	\left(
	\frac{\omega\lambda_1(\xi_1-1)}{2},
	\frac{\omega\xi_2(\xi_1+1)}{2},
	\frac{\omega\lambda_1\lambda_2(\xi_1+1)}{2},
	\omega^2\xi_1\xi_2\lambda_1\lambda_2
	\right)\\
	&=
	q_\omega\bigl(\Psi_{\mathbb H}(\xi_1,\xi_2,\lambda_1,\lambda_2)\bigr).
\end{align*}
Since the map $(\xi_1,\xi_2,\lambda_1,\lambda_2)
\mapsto
(\xi_1,\omega\xi_2,\omega\lambda_1,\lambda_2)
$
preserves the normalized Haar measure on $\mathbb T^4$ and
$\mu_{\mathbb H}=(\Psi_{\mathbb H})_*m_4$, we have that
$(q_\omega)_*\mu_{\mathbb H}=\mu_{\mathbb H}$, that is,
\[
\int_{b\mathbb H}
F(\omega a,\omega x_1,\omega x_2,\omega^2x_3)
\,d\mu_{\mathbb H}
=
\int_{b\mathbb H}F(a,x_1,x_2,x_3)\,d\mu_{\mathbb H}
\]
for $F\in C(b\mathbb H)$. Thus, the measure $\mu_{\mathbb H}$ is invariant under two natural symmetries of $b\mathbb H$. The first rotates each circular $a$-fibre while leaving its base point in $b\mathbb E$ fixed. The second carries the fibre over $(x_1,x_2,x_3)$ onto the fibre over $(\omega x_1,\omega x_2,\omega^2x_3)$. These invariance properties follow from the rotational invariance of Haar measure on $\mathbb T$ and the invariance of Haar measure on $\mathbb T^4$ under the corresponding torus transformation. We now define the boundary Hardy space associated with this measure.

\begin{definition}
	The boundary Hardy space of $\mathbb H$ is defined by
	\[
	H^2(\mathbb H)
	=
	\overline{\mathbb C[a,x_1,x_2,x_3]}^{\,L^2(b\mathbb H,\mu_{\mathbb H})},
	\]
	equipped with the inner product
	\[
	\langle f,g\rangle_{H^2(\mathbb H)}
	=
	\int_{b\mathbb H}
	f(a,x_1,x_2,x_3)
	\overline{g(a,x_1,x_2,x_3)}
	\,d\mu_{\mathbb H}.
	\]
\end{definition}

The two invariance properties of $\mu_{\mathbb H}$ established above give rise to two natural orthogonal decompositions of $H^2(\mathbb H)$. The invariance under the $(1,1,1,2)$-quasi-circular action leads to a decomposition according to weighted homogeneous degree and the rotational invariance of the $a$-fibres yields a decomposition according to the powers of the first variable $a$. 

\smallskip

\noindent \textbf{Weighted homogeneous decomposition.} For $n\in\mathbb N_0$, let $\mathcal Q_n$ denote the space of weighted homogeneous polynomials of degree $n$, that is, $q(\omega a,\omega x_1,\omega x_2,\omega^2x_3)
=
\omega^n q(a,x_1,x_2,x_3)$ for every $\omega\in\mathbb T$. More precisely,
\[
\mathcal Q_n
=
\operatorname{span}
\left\{
a^j x_1^k x_2^\ell x_3^m:
j,k,\ell,m\geq0,\;
j+k+\ell+2m=n
\right\}.
\]
The invariance of $\mu_{\mathbb H}$ under the quasi-circular action implies that $\mathcal Q_n\perp\mathcal Q_m$ whenever $n\neq m$. Since polynomials are dense in $H^2(\mathbb H)$ by definition, we obtain the orthogonal decomposition
\begin{equation}\label{eq:orth-hexablock}
	H^2(\mathbb H)
	=
	\bigoplus_{n=0}^{\infty}\mathcal Q_n.
\end{equation}
Moreover, the dimension of $\mathcal Q_n$ is given by
\[
\dim\mathcal Q_n
=
\sum_{m=0}^{\lfloor n/2\rfloor}
\binom{n-2m+2}{2}.
\]

\noindent \textbf{Fibre decomposition.} The second decomposition is motivated by the circular fibre structure in the $a$-variable. For each fixed $(x_1,x_2,x_3)\in b\mathbb E$, the corresponding fibre is
\[
\left\{
\bigl(\xi\rho(x_1,x_2,x_3),x_1,x_2,x_3\bigr):
\xi\in\mathbb T
\right\}.
\]
For $j,k\in\mathbb Z_+$, we have
\begin{align*}
	\int_{\mathbb T}
	\bigl(\xi\rho(x_1,x_2,x_3)\bigr)^k
	\overline{\bigl(\xi\rho(x_1,x_2,x_3)\bigr)^j}
	\,dm_1(\xi)
	&=
	\rho(x_1,x_2,x_3)^{k+j}
	\int_{\mathbb T}\xi^{k-j}\,dm_1(\xi)\\
	&=
	\begin{cases}
		\rho(x_1,x_2,x_3)^{2k}, & j=k,\\[4pt]
		0, & j\neq k.
	\end{cases}
\end{align*}
Consequently, distinct powers of the fibre variable $a$ are orthogonal in $L^2(b\mathbb H,\mu_{\mathbb H})$. Next, if $f\in\mathbb C[x_1,x_2,x_3]$, then
\[
\begin{aligned}
	\|a^k f\|_{H^2(\mathbb H)}^2
	&=
	\int_{b\mathbb E}\int_{\mathbb T}
	\bigl|\xi\rho(x_1,x_2,x_3)\bigr|^{2k}
	|f(x_1,x_2,x_3)|^2
	\,dm_1(\xi)\,d\nu_{\mathbb H}\\
	&=
	\int_{b\mathbb E}
	|f(x_1,x_2,x_3)|^2
	\rho(x_1,x_2,x_3)^{2k}
	\,d\nu_{\mathbb H}.
\end{aligned}
\]
Thus, the natural coefficient space corresponding to the $k$-th power of the fibre variable is the weighted Hardy space
\[
H_k^2(\mathbb E)
=
\overline{\mathbb C[x_1,x_2,x_3]}^{
	\,L^2(b\mathbb E,\rho^{2k}\,d\nu_{\mathbb H})
}
\ \ \text{with norm} \ \ 
\|f\|_{H_k^2(\mathbb E)}^2
=
\int_{b\mathbb E}
|f(x_1,x_2,x_3)|^2
\rho(x_1,x_2,x_3)^{2k}
\,d\nu_{\mathbb H}.
\]
Since $\rho(x_1,x_2,x_3)>0$ for $\nu_{\mathbb H}$-almost every $(x_1,x_2,x_3)\in b\mathbb E$, this defines a norm on $\mathbb C[x_1,x_2,x_3]$. Every polynomial $q\in\mathbb C[a,x_1,x_2,x_3]$ admits a unique decomposition
\[
q(a,x_1,x_2,x_3)
=
\sum_{k=0}^N a^k q_k(x_1,x_2,x_3),
\]
where $q_k\in\mathbb C[x_1,x_2,x_3]$. By the orthogonality of distinct powers of $a$, we have
\[
\|q\|_{H^2(\mathbb H)}^2
=
\sum_{k=0}^N
\|a^kq_k\|_{H^2(\mathbb H)}^2
=
\sum_{k=0}^N
\|q_k\|_{H_k^2(\mathbb E)}^2.
\]
Motivated by the above discussion, we obtain the following result.

\begin{proposition}\label{prop:decompositions-hexablock}
	The Hardy space $H^2(\mathbb H)$ admits the two orthogonal
	decompositions
	\begin{equation}\label{eq:quasi-decomposition-hexablock}
		H^2(\mathbb H)
		=
		\bigoplus_{n=0}^{\infty}\mathcal Q_n
		=
		\bigoplus_{k=0}^{\infty}
		a^k H_k^2(\mathbb E),
	\end{equation}
	where $\mathcal Q_n$ is the space of weighted homogeneous polynomials
	of degree $n$. Moreover, every $f\in H^2(\mathbb H)$ has a unique expansion
	\[
	f(a,x_1,x_2,x_3)
	=
	\sum_{k=0}^{\infty}a^k f_k(x_1,x_2,x_3),
	\]
	where $f_k\in H_k^2(\mathbb E)$ and the series converges in $H^2(\mathbb H)$. Also,
	$\displaystyle
	\|f\|_{H^2(\mathbb H)}^2
	=
	\sum_{k=0}^{\infty}
	\|f_k\|_{H_k^2(\mathbb E)}^2.
	$
\end{proposition}

\begin{proof}
	By the above discussion, it suffices to establish the expansion of $f$
	as in the statement. Let $f\in H^2(\mathbb H)$. Since polynomials are
	dense in $H^2(\mathbb H)$, there exists a sequence of polynomials
	$\{q^{(n)}\}_{n\geq1}$ such that
	$q^{(n)}\to f$ in $H^2(\mathbb H)$. For each $n$, write
	\[
	q^{(n)}(a,x_1,x_2,x_3)
	=
	\sum_{k=0}^{N_n}a^kq_k^{(n)}(x_1,x_2,x_3),
	\qquad
	q_k^{(n)}\in\mathbb C[x_1,x_2,x_3],
	\]
	and set $q_k^{(n)}=0$ for $k>N_n$. By the orthogonality of distinct
	powers of $a$ and the definition of $H_k^2(\mathbb E)$, we have
	\[
	\left\|
	\sum_{k=0}^{\infty}a^kq_k^{(n)}
	\right\|_{H^2(\mathbb H)}^2
	=
	\sum_{k=0}^{\infty}
	\|q_k^{(n)}\|_{H_k^2(\mathbb E)}^2
	\ \ \text{and} \ \
	\|q^{(n)}-q^{(m)}\|_{H^2(\mathbb H)}^2
	=
	\sum_{k=0}^{\infty}
	\|q_k^{(n)}-q_k^{(m)}\|_{H_k^2(\mathbb E)}^2.
	\]
	Since $\{q^{(n)}\}$ converges in $H^2(\mathbb H)$, it is a Cauchy
	sequence. Thus, for each $k\geq0$, the sequence
	$\{q_k^{(n)}\}_{n\geq1}$ is Cauchy in $H_k^2(\mathbb E)$. By the
	completeness of $H_k^2(\mathbb E)$, there exists
	$f_k\in H_k^2(\mathbb E)$ such that $q_k^{(n)}\to f_k$ in $H_k^2(\mathbb E)$. Next, we show that the series
	\[
	\sum_{k=0}^{\infty}a^kf_k
	\]
	converges in $H^2(\mathbb H)$. Let $C>0$ be such that
	$\|q^{(n)}\|_{H^2(\mathbb H)}^2\leq C$ for all $n\geq1$.
	Fix $N\in\mathbb N$. Then
	\[
	\sum_{k=0}^{N}
	\|q_k^{(n)}\|_{H_k^2(\mathbb E)}^2
	\leq C \quad \text{and so,} \quad 
	\sum_{k=0}^{N}
	\|f_k\|_{H_k^2(\mathbb E)}^2
	\leq C.
	\]
	Since $N$ is arbitrary, it follows that
	\[
	\sum_{k=0}^{\infty}
	\|f_k\|_{H_k^2(\mathbb E)}^2<\infty \quad \text{and thus} \quad 
	\left\|
	\sum_{k=N+1}^{M}a^kf_k
	\right\|_{H^2(\mathbb H)}^2
	=
	\sum_{k=N+1}^{M}
	\|f_k\|_{H_k^2(\mathbb E)}^2
	\to 0
	\]
	as $M,N\to\infty$. Since $H^2(\mathbb H)$ is complete, there
	exists $g\in H^2(\mathbb H)$ such that
	$
	g=\sum_{k=0}^{\infty}a^kf_k
	$
	and the series converges in $H^2(\mathbb H)$. It remains to show that $g=f$. Fix $n\geq1$. Since
	$q_k^{(m)}\to f_k$ in $H_k^2(\mathbb E)$ for every $k$, Fatou's
	lemma gives
	\[
	\begin{aligned}
		\sum_{k=0}^{\infty}
		\|f_k-q_k^{(n)}\|_{H_k^2(\mathbb E)}^2
		\leq
		\liminf_{m\to\infty}
		\sum_{k=0}^{\infty}
		\|q_k^{(m)}-q_k^{(n)}\|_{H_k^2(\mathbb E)}^2
		&=
		\liminf_{m\to\infty}
		\|q^{(m)}-q^{(n)}\|_{H^2(\mathbb H)}^2\\
		&=
		\|f-q^{(n)}\|_{H^2(\mathbb H)}^2.
	\end{aligned}
	\]
	Also, the orthogonality of the summands in the
	representations of $g$ and $q^{(n)}$ gives
	\[
	\|g-q^{(n)}\|_{H^2(\mathbb H)}^2
	=
	\sum_{k=0}^{\infty}
	\|f_k-q_k^{(n)}\|_{H_k^2(\mathbb E)}^2
	\leq
	\|f-q^{(n)}\|_{H^2(\mathbb H)}^2.
	\]
	Thus,
	$
	\|f-g\|_{H^2(\mathbb H)}
	\leq
	\|f-q^{(n)}\|_{H^2(\mathbb H)}
	+
	\|q^{(n)}-g\|_{H^2(\mathbb H)}
	\leq
	2\|f-q^{(n)}\|_{H^2(\mathbb H)}
	\to 0
	$
	as $n\to\infty$. Hence $f=g$. The uniqueness of the expansion follows from the orthogonality of the subspaces $a^kH_k^2(\mathbb E)$. Indeed, if
	\[
	\sum_{k=0}^{\infty}a^kf_k=0,
	\quad \text{then} \quad 
	\|f_j\|_{H_j^2(\mathbb E)}^2
	=
	\left\langle
	\sum_{k=0}^{\infty}a^kf_k,a^jf_j
	\right\rangle_{H^2(\mathbb H)}
	=
	0
	\]
	and hence $f_j=0$ for every $j\geq0$. Finally, since the partial sums converge to $f$ and the summands
	are orthogonal, we have
	\[
	\left\|
	\sum_{k=0}^{N}a^kf_k
	\right\|_{H^2(\mathbb H)}^2
	=
	\sum_{k=0}^{N}
	\|f_k\|_{H_k^2(\mathbb E)}^2
	\]
	for every $N \geq 0$. Letting $N\to\infty$, it follows that
	$\displaystyle
	\|f\|_{H^2(\mathbb H)}^2
	=
	\sum_{k=0}^{\infty}
	\|f_k\|_{H_k^2(\mathbb E)}^2.
	$
\end{proof}

Next, we show that every element of $H^2(\mathbb H)$ admits a holomorphic realization on $\mathbb H$ and that $H^2(\mathbb H)$ is a reproducing kernel Hilbert space. To this end, we first establish an estimate that allows us to pass from the orthogonal decomposition of $H^2(\mathbb H)$ into weighted homogeneous polynomials to a locally uniformly convergent holomorphic series.

\begin{lemma}\label{lemma:hexablock-estimate}
	Let $K\Subset\mathbb H$. There exist constants $0<r_K<1$ and $C_K>0$ such that
	\[
	\sup_{z\in K}|q(z)|
	\leq
	C_K(n+1)^2r_K^n
	\|q\|_{L^2(b\mathbb H,\mu_{\mathbb H})}
	\]
	for every $q\in\mathcal Q_n$, where
	$
	\mathcal Q_n
	=
	\operatorname{span}
	\left\{
	a^j x_1^k x_2^\ell x_3^m:
	j,k,\ell,m\geq0,\;
	j+k+\ell+2m=n
	\right\}.
	$
\end{lemma}

\begin{proof}
	Since $K\Subset\mathbb H$ and $\mathbb H$ is $(1,1,1,2)$-quasi-balanced, there exists $R>1$ (sufficiently close to $1$) such that
	$\delta_R(K)\Subset\mathbb H$, where $\delta_r(a,x_1,x_2,x_3)=(ra,rx_1,rx_2,r^2x_3)$. Put $r_K=R^{-1}$. Since $q\in\mathcal Q_n$, its weighted homogeneity gives
	\[
	q(\lambda a,\lambda x_1,\lambda x_2,\lambda^2x_3)
	=
	\lambda^nq(a,x_1,x_2,x_3)
	\]
	for every $\lambda\in\mathbb C$. Then 
	\[
	|q(z)|
	=
	R^{-n}|q(w)|=r_K^{n}|q(w)|
	\leq
	r_K^n\|q\|_{L^\infty(b\mathbb H)}
	\]
	for $z\in K$ and $w=\delta_R(z)$. We now estimate the above inequality in terms of the $L^2$-norm. Set $Q=q\circ\Psi_{\mathbb H}$. Since $q\in\mathcal Q_n$, the function $Q$ is a trigonometric polynomial on $\mathbb T^4$ involving at most $(n+1)^4$ Fourier monomials. Write
	$Q=\sum_{\alpha\in A}c_\alpha\chi_\alpha$ with the cardinality $|A|\leq(n+1)^4$. By the Cauchy-Schwarz inequality and the orthogonality of the characters of $\mathbb T^4$, we have
	\begin{align*}
		\|q\|_{L^\infty(b\mathbb H)}
		=
		\|Q\|_{L^\infty(\mathbb T^4)}
		\leq
		\sum_{\alpha\in A}|c_\alpha|
		\leq
		|A|^{1/2}
		\left(\sum_{\alpha\in A}|c_\alpha|^2\right)^{1/2}
		\leq
		(n+1)^2\|Q\|_{L^2(\mathbb T^4)}.
	\end{align*}
	Since $\mu_{\mathbb H}=(\Psi_{\mathbb H})_*m_4$, we have that $\|Q\|_{L^2(\mathbb T^4)}
	=
	\|q\|_{L^2(b\mathbb H,\mu_{\mathbb H})}$. Combining the above estimates, the desired conclusion follows. 
\end{proof}

The estimate in Lemma \ref{lemma:hexablock-estimate} provides sufficient decay in the weighted degree to ensure local uniform convergence of the homogeneous components. We therefore obtain the following holomorphic realization theorem.

\begin{theorem}[Holomorphic realization]\label{thm:hex-hol-realization}
	Let $f\in H^2(\mathbb H)$ and let
	$
	f=\sum_{n=0}^{\infty}q_n
	$
	be its weighted homogeneous decomposition with each $q_n\in\mathcal Q_n$. Then, for every $z\in\mathbb H$, the series
	$
	\sum_{n=0}^{\infty}q_n(z)
	$
	converges absolutely and the convergence is uniform on compact subsets of $\mathbb H$. The resulting map
	\[
	T:H^2(\mathbb H) \to \operatorname{Hol}(\mathbb H),
	\qquad
	(Tf)(z)=\sum_{n=0}^{\infty}q_n(z),
	\]
	is injective and continuous when $\operatorname{Hol}(\mathbb H)$ is equipped with the compact-open topology. In particular, $H^2(\mathbb H)$ is a reproducing kernel Hilbert space.
\end{theorem}

\begin{proof}
	Let $K\Subset\mathbb H$ and $n\geq0$. By Lemma \ref{lemma:hexablock-estimate}, there exist constants $C_K>0$ and $0<r_K<1$ such that
	$\sup_{z\in K}|q_n(z)|
	\leq
	C_K(n+1)^2r_K^n
	\|q_n\|_{L^2(b\mathbb H,\mu_{\mathbb H})}$. Since the decomposition of $f$ is orthogonal in $H^2(\mathbb H)$, we have
	\[
	\sum_{n=0}^{\infty}
	\|q_n\|_{L^2(b\mathbb H,\mu_{\mathbb H})}^2
	=
	\|f\|_{H^2(\mathbb H)}^2<\infty.
	\]
	By the Cauchy-Schwarz inequality,
	\[
	\begin{aligned}
		&\sum_{n=0}^{\infty}
		C_K(n+1)^2r_K^n
		\|q_n\|_{L^2(b\mathbb H,\mu_{\mathbb H})} \leq
		C_K
		\left(
		\sum_{n=0}^{\infty}(n+1)^4r_K^{2n}
		\right)^{1/2}
		\left(
		\sum_{n=0}^{\infty}
		\|q_n\|_{L^2(b\mathbb H,\mu_{\mathbb H})}^2
		\right)^{1/2}.
	\end{aligned}
	\]
	The first series on the right is finite since $0<r_K<1$, while the second one equals $\|f\|_{H^2(\mathbb H)}^2$. By the Weierstrass $M$-test, the series
	\[
	\sum_{n=0}^{\infty}q_n(z)
	\]
	converges absolutely and uniformly on $K$. Since $K\Subset\mathbb H$ was arbitrary, the convergence is locally uniform on $\mathbb H$. Each $q_n$ is a polynomial and hence holomorphic on $\mathbb H$. Therefore, by the Weierstrass theorem, the function
	\[
	Tf(z)=\sum_{n=0}^{\infty}q_n(z)
	\]
	is holomorphic on $\mathbb H$. We next prove that $T$ is injective. Suppose that $Tf\equiv0$ on $\mathbb H$. Fix $z=(a,x_1,x_2,x_3)\in\mathbb H$ and let $\lambda\in\mathbb D$. Then 
	\[
	\begin{aligned}
		0
		=Tf(\lambda a,\lambda x_1,\lambda x_2,\lambda^2x_3)
		=\sum_{n=0}^{\infty}
		q_n(\lambda a,\lambda x_1,\lambda x_2,\lambda^2x_3)
		=\sum_{n=0}^{\infty}\lambda^nq_n(a,x_1,x_2,x_3).
	\end{aligned}
	\]
	The second equality follows from the locally uniform convergence established above and the third one follows from the weighted homogeneity of $q_n$. Consequently, the uniqueness of the Taylor coefficients gives $q_n(z)=0$ for every $n\geq0$. Since $z\in\mathbb H$ was arbitrary, we have that $q_n\equiv0$ for every $n$ and hence $f=0$ in $H^2(\mathbb H)$. Thus $T$ is injective. Next, for every $K\Subset\mathbb H$, the above estimates give a constant $A_K>0$ (depending only on $K$) such that
	\[
	\sup_{z\in K}|Tf(z)|
	\leq
	A_K\|f\|_{H^2(\mathbb H)}.
	\]
	This proves that $T$ is continuous from $H^2(\mathbb H)$ into $\operatorname{Hol}(\mathbb H)$ endowed with the compact-open topology. In particular, point evaluations are bounded linear functionals on $H^2(\mathbb H)$. Hence $H^2(\mathbb H)$ is a reproducing kernel Hilbert space.
\end{proof}

The above theorem identifies each element of the boundary Hardy
space $H^2(\mathbb H)$ with a holomorphic function on $\mathbb H$.
We next characterize this holomorphic realization in terms of radial boundary norms. Recall that $\mathbb H$ is
$(1,1,1,2)$-quasi-balanced. Thus, the radial maps
\[
\delta_r(a,x_1,x_2,x_3)
=
(ra,rx_1,rx_2,r^2x_3),
\qquad 0\leq r<1,
\]
satisfy
$\delta_r(\overline{\mathbb H})\subseteq\mathbb H$
for every $r\in(0,1)$. This leads naturally to the following
radial Hardy space:
\[
H^2_{\mathrm{rad}}(\mathbb H)
=
\left\{
f\in\operatorname{Hol}(\mathbb H):
\|f\|_{H^2_{\mathrm{rad}}(\mathbb H)}<\infty
\right\},
\]
where
\[
\|f\|_{H^2_{\mathrm{rad}}(\mathbb H)}
=
\sup_{0<r<1}
\left(
\int_{b\mathbb H}
|f\circ\delta_r|^2\,d\mu_{\mathbb H}
\right)^{1/2}.
\]
We first verify that the above quantity defines a norm. Let
$\|f\|_{H^2_{\mathrm{rad}}(\mathbb H)}=0$.
Then $f\circ\delta_r=0$ $\mu_{\mathbb H}$-almost everywhere on
$b\mathbb H$ for every $r\in(0,1)$. Since
$\mu_{\mathbb H}=(\Psi_{\mathbb H})_*m_4$,
$\Psi_{\mathbb H}(\mathbb T^4)=b\mathbb H$ and the normalized
Haar measure $m_4$ has full support on $\mathbb T^4$, it follows that the
measure $\mu_{\mathbb H}$ has full support on $b\mathbb H$.
Moreover, $f\circ\delta_r$ is continuous on
$\overline{\mathbb H}$. Consequently,
$f\circ\delta_r=0$ on $b\mathbb H$. For $r\in(0,1)$, define
\[
\mathbb H_r
=
\delta_r(\overline{\mathbb H})
=
\left\{
(a,x_1,x_2,x_3):
(r^{-1}a,r^{-1}x_1,r^{-1}x_2,r^{-2}x_3)
\in\overline{\mathbb H}
\right\}.
\]
It was proved in \cite{BiswasPalI} that
\[
\mathbb H=\bigcup_{0<r<1}\mathbb H_r.
\]
Fix $r\in(0,1)$ and put $f_r=f\circ\delta_r$.
Since $\delta_r(\overline{\mathbb H})\subset\mathbb H$,
the function $f_r$ is holomorphic on $\mathbb H$ and continuous
on $\overline{\mathbb H}$. We have already shown that
$f_r=0$ on $b\mathbb H$. Since $b\mathbb H$ is the distinguished
boundary of $\mathbb H$, it follows that
$f_r=0$ on $\overline{\mathbb H}$. Consequently,
$f=0$ on $\mathbb H_r$ for every $r\in(0,1)$ and hence
$f=0$ on $\mathbb H$. Thus,
$\|f\|_{H^2_{\mathrm{rad}}(\mathbb H)}=0$
if and only if $f=0$. The absolute homogeneity and the
triangle inequality follow from the corresponding properties
of the $L^2(b\mathbb H,\mu_{\mathbb H})$ norm.
Therefore, $H^2_{\mathrm{rad}}(\mathbb H)$ is a normed space.
We shall establish its completeness after identifying it
isometrically with the boundary Hardy space $H^2(\mathbb H)$.
The next result provides this identification and shows that
the radial norm coincides with the boundary Hardy-space norm.

\begin{theorem}[Radial realization]
	\label{thm:radial-realization-hex}
	The holomorphic realization identifies
	$H^2(\mathbb H)$ isometrically with
	$H^2_{\mathrm{rad}}(\mathbb H)$. More precisely, assume that
	$f\in H^2(\mathbb H)$ and
	$
	f=\sum_{n=0}^{\infty}q_n
	$
	is its weighted homogeneous decomposition, where
	$q_n\in\mathcal Q_n$. Then
	\[
	\|f\|_{H^2(\mathbb H)}^2
	=
	\sup_{0<r<1}
	\int_{b\mathbb H}|f\circ\delta_r|^2\,d\mu_{\mathbb H} \ \ \text{and} \ \ \int_{b\mathbb H}|f\circ\delta_r|^2\,d\mu_{\mathbb H}
	=
	\sum_{n=0}^{\infty}
	r^{2n}\|q_n\|_{L^2(b\mathbb H,\mu_{\mathbb H})}^2
	\]
	for every $0<r<1$. Conversely, every
	$f\in H^2_{\mathrm{rad}}(\mathbb H)$ admits a unique expansion $f=\sum_{n=0}^{\infty}q_n$ with $q_n\in\mathcal Q_n$, where the series converges locally uniformly on $\mathbb H$
	and
	\[
	\sum_{n=0}^{\infty}
	\|q_n\|_{L^2(b\mathbb H,\mu_{\mathbb H})}^2<\infty.
	\]
\end{theorem}

\begin{proof}
	Let $f\in H^2(\mathbb H)$ and consider its weighted homogeneous decomposition
	\[
	f=\sum_{n=0}^{\infty}q_n \quad (q_n \in \mathcal Q_n).
	\]
	By Theorem
	\ref{thm:hex-hol-realization}, the series converges locally
	uniformly on $\mathbb H$ to the holomorphic realization of
	$f$. Since
	$\delta_r(\overline{\mathbb H})\Subset\mathbb H$
	for every $r\in(0,1)$, it follows that
	$
	f\circ\delta_r
	=
	\sum_{n=0}^{\infty}r^nq_n
	$
	pointwise on $b\mathbb H$. By the orthogonality of the spaces $\mathcal Q_n$, we have
	\begin{align*}
		\sum_{n=0}^{\infty}
		\|r^nq_n\|_{L^2(b\mathbb H,\mu_{\mathbb H})}^2
		=
		\sum_{n=0}^{\infty}
		r^{2n}\|q_n\|_{L^2(b\mathbb H,\mu_{\mathbb H})}^2
		\leq
		\sum_{n=0}^{\infty}
		\|q_n\|_{L^2(b\mathbb H,\mu_{\mathbb H})}^2
		=
		\|f\|_{H^2(\mathbb H)}^2.
	\end{align*}
	Hence, the series $\sum_{n=0}^{\infty}r^nq_n$ converges in $L^2(b\mathbb H,\mu_{\mathbb H})$. Consequently,
	\begin{align*}
		\int_{b\mathbb H}|f\circ\delta_r|^2\,d\mu_{\mathbb H}
		=
		\left\|
		\sum_{n=0}^{\infty}r^nq_n
		\right\|_{L^2(b\mathbb H,\mu_{\mathbb H})}^2
		=
		\sum_{n=0}^{\infty}
		r^{2n}\|q_n\|_{L^2(b\mathbb H,\mu_{\mathbb H})}^2.
	\end{align*}
	Since the terms on the right-hand side are nonnegative
	and increase with $r$, the monotone convergence theorem gives
	\begin{align*}
		\sup_{0<r<1}
		\int_{b\mathbb H}|f\circ\delta_r|^2\,d\mu_{\mathbb H}
		=
		\lim_{r\uparrow1}
		\sum_{n=0}^{\infty}
		r^{2n}\|q_n\|_{L^2(b\mathbb H,\mu_{\mathbb H})}^2
		=
		\sum_{n=0}^{\infty}
		\|q_n\|_{L^2(b\mathbb H,\mu_{\mathbb H})}^2
		=
		\|f\|_{H^2(\mathbb H)}^2.
	\end{align*}
	Thus, $f\in H^2_{\mathrm{rad}}(\mathbb H)$ and $\|f\|_{H^2(\mathbb H)}=\|f\|_{H^2_{\mathrm{rad}}(\mathbb H)}$.
	
	\smallskip 
	
	Conversely, let $f\in H^2_{\mathrm{rad}}(\mathbb H)$. Since $\mathcal H$ is a $(1,1,1,2)$-quasi-balanced domain containing the origin, for every $z=(a,x_1,x_2,x_3)\in \mathbb H$, the function
	\[
	\lambda\longmapsto
	f(\lambda a,\lambda x_1,\lambda x_2,\lambda^2x_3)
	\]
	is holomorphic on $\mathbb D$. Hence, it admits a Taylor expansion
	\[
	f(\lambda a,\lambda x_1,\lambda x_2,\lambda^2x_3)
	=
	\sum_{n=0}^{\infty}\lambda^n q_n(a,x_1,x_2,x_3),
	\qquad \lambda\in\mathbb D.
	\]
	We now show that each $q_n$ is a weighted homogeneous polynomial of degree $n$. Since $f$ is holomorphic on $\mathbb H$ and $0\in\mathbb H$, its Taylor expansion in a neighborhood of the origin is given by
	\[
	f(a,x_1,x_2,x_3)
	=
	\sum_{j,k,\ell,m\geq 0}
	c_{jk\ell m}a^jx_1^kx_2^\ell x_3^m.
	\]
	Substituting $(\lambda a,\lambda x_1,\lambda x_2,\lambda^2x_3)$
	into this expansion, it follows that
	\[
	f(\lambda a,\lambda x_1,\lambda x_2,\lambda^2x_3)
	=
	\sum_{j,k,\ell,m\geq 0}
	c_{jk\ell m}
	\lambda^{j+k+\ell+2m}
	a^jx_1^kx_2^\ell x_3^m.
	\]
	Collecting the terms corresponding to the same power of $\lambda$,
	we deduce that
	\[
	q_n(a,x_1,x_2,x_3)
	=
	\sum_{\substack{j,k,\ell,m\geq 0\\
			j+k+\ell+2m=n}}
	c_{jk\ell m}a^jx_1^kx_2^\ell x_3^m.
	\]
	In particular, each $q_n$ is a polynomial belonging to the
	weighted homogeneous polynomial space  $\mathcal Q_n$. The identity theorem therefore yields
	\[
	f(\lambda a,\lambda x_1,\lambda x_2,\lambda^2x_3)
	=
	\sum_{n=0}^{\infty}\lambda^nq_n(a,x_1,x_2,x_3),
	\qquad \lambda\in\mathbb D.
	\]
	In particular, letting $\lambda \to 1$ gives
	\[
	f(a,x_1,x_2,x_3)
	=
	\sum_{n=0}^{\infty}q_n(a,x_1,x_2,x_3)
	\]  
	for every  $(a,x_1,x_2,x_3)\in\mathbb H$. Thus, $f$ admits the desired weighted homogeneous decomposition $f=\sum_{n=0}^{\infty}q_n$ with  $q_n\in\mathcal Q_n$. We next establish the square-summability of these
	homogeneous components. Fix $0<r<1$ and define
	\[
	F_r(a,x_1,x_2,x_3)
	=
	f(ra,rx_1,rx_2,r^2x_3), \quad (a, x_1, x_2, x_3) \in b\mathbb H.
	\]
	By the above Taylor expansion of $f$, we have for every $(a, x_1, x_2, x_3) \in b\mathbb H$ that
	$F_r(a,x_1,x_2,x_3)=\sum_{n=0}^{\infty}
	r^nq_n(a,x_1,x_2,x_3)$. Let $(a,x_1,x_2,x_3)\in b\mathbb H$. Then 
	\[
	F_r(\omega a, \omega x_1, \omega x_2, \omega^2 x_3)
	=
	f(r\omega a,r\omega x_1,
	r\omega x_2,r^2\omega^2x_3)
	=
	\sum_{n=0}^{\infty}
	r^n\omega^nq_n(a,x_1,x_2,x_3).
	\]
	By Parseval's identity, 
	\[
	\int_{\mathbb T}
	|F_r(\omega a, \omega x_1, \omega x_2, \omega^2 x_3)|^2\,dm_1(\omega)
	=
	\sum_{n=0}^{\infty}
	r^{2n}|q_n(a, x_1, x_2, x_3)|^2.
	\]
	Integrating both sides over $b\mathbb H$ with respect
	to $\mu_{\mathbb H}$, and using Fubini's theorem and
	the invariance of $\mu_{\mathbb H}$ under the
	$(1,1,1,2)$-quasi-circular action, we have that
	\[
	\int_{b\mathbb H}|F_r|^2\,d\mu_{\mathbb H}
	=
	\int_{b\mathbb H}
	\int_{\mathbb T}
	|F_r(q_\omega(z))|^2
	\,dm_1(\omega)\,d\mu_{\mathbb H}(z)
	=
	\sum_{n=0}^{\infty}
	r^{2n}
	\|q_n\|_{L^2(b\mathbb H,\mu_{\mathbb H})}^2.
	\]
	Since $F_r=f\circ\delta_r$ on $b\mathbb H$, it follows
	that
	\[
	\sup_{0<r<1}
	\sum_{n=0}^{\infty}
	r^{2n}
	\|q_n\|_{L^2(b\mathbb H,\mu_{\mathbb H})}^2
	=
	\|f\|_{H^2_{\mathrm{rad}}(\mathbb H)}^2
	<\infty.
	\]
	Letting $r\uparrow1$ and applying the monotone
	convergence theorem, we have
	\[
	\sum_{n=0}^{\infty}
	\|q_n\|_{L^2(b\mathbb H,\mu_{\mathbb H})}^2
	<\infty
	\]
	Thus, the series $g=\sum_{n=0}^{\infty}q_n$ defines an element of $H^2(\mathbb H)$. By Theorem \ref{thm:hex-hol-realization}, the series of $g$ is a locally uniformly convergent series on $\mathbb H$. It remains to show that $g=f$. For every $z=(a,x_1,x_2,x_3)\in\mathbb H$, we have by
	the one-variable Taylor expansion that
	\[
	f(\lambda a,\lambda x_1,\lambda x_2,\lambda^2x_3)
	=
	\sum_{n=0}^{\infty}
	\lambda^nq_n(a,x_1,x_2,x_3),
	\qquad |\lambda|<1.
	\]
	Since $\mathbb H$ is open and $z\in\mathbb H$, there
	exists $\rho>1$ sufficiently close to $1$ such that
	$(\lambda a,\lambda x_1,\lambda x_2,\lambda^2x_3)
	\in\mathbb H$ for every $|\lambda|<\rho$. In particular, we may take $\lambda=1$ in the
	corresponding Taylor expansion. Hence,
	\[
	f(a,x_1,x_2,x_3)
	=
	\sum_{n=0}^{\infty}
	q_n(a,x_1,x_2,x_3).
	\]
	Therefore, $f=g$ on $\mathbb H$, which completes the proof.
\end{proof}

As a consequence of Theorem \ref{thm:radial-realization-hex}, we have the following corollary.

\begin{corollary}[Radial approximation]
	\label{cor:radial-approx-H}
	For $f\in H^2(\mathbb H)$, let
	$f_r=f\circ\delta_r$, where
	$
	\delta_r(a,x_1,x_2,x_3)
	=(ra,rx_1,rx_2,r^2x_3).
	$
	Then $\lim_{r\uparrow1}
	\|f_r-f\|_{H^2(\mathbb H)}=0$.
	
\end{corollary}

\begin{proof}
	Let $f=\sum_{n=0}^{\infty}q_n$ be the weighted homogeneous decomposition of $f$ in
	$H^2(\mathbb H)$, where
	$q_n\in\mathcal Q_n$. By the weighted homogeneity of $q_n$, we have that
	$f_r=\sum_{n=0}^{\infty}r^nq_n$. By the orthogonality of the spaces
	$\mathcal Q_n$ in $L^2(b\mathbb H,\mu_{\mathbb H})$, we have
	\[
	\|f_r-f\|_{H^2(\mathbb H)}^2
	=\sum_{n=0}^{\infty}
	|r^n-1|^2
	\|q_n\|_{L^2(b\mathbb H,\mu_{\mathbb H})}^2.
	\]
	Since $\sum_{n=0}^{\infty}
	\|q_n\|_{L^2(b\mathbb H,\mu_{\mathbb H})}^2
	=
	\|f\|_{H^2(\mathbb H)}^2<\infty$, the dominated convergence theorem gives the desired conclusion.
\end{proof}

We next identify $H^2(\mathbb H)$ with a closed subspace of 
$L^2(b\mathbb H, \mu_{\mathbb H})$ via radial boundary values.

\begin{corollary}[Radial boundary values]
	\label{cor:radial-boundary-H}
	For every $f\in H^2(\mathbb H)$, there exists a unique
	$f^*\in L^2(b\mathbb H,\mu_{\mathbb H})$ such that
	\[
	\lim_{r\uparrow1}
	\|f\circ\delta_r-f^*\|_
	{L^2(b\mathbb H,\mu_{\mathbb H})}=0 \quad \text{and} \quad 
	\|f\|_{H^2(\mathbb H)}
	=
	\|f^*\|_{L^2(b\mathbb H,\mu_{\mathbb H})}.
	\]
\end{corollary}

\begin{proof}
	Let $f=\sum_{n=0}^{\infty}q_n$ be the weighted homogeneous decomposition of $f$, where $q_n\in\mathcal Q_n$. We have by Theorem \ref{thm:radial-realization-hex} that
	\[
	\sum_{n=0}^{\infty}
	\|q_n\|_{L^2(b\mathbb H,\mu_{\mathbb H})}^2
	=
	\|f\|_{H^2(\mathbb H)}^2<\infty.
	\]
	Since $\mathcal Q_n$'s are
	mutually orthogonal in
	$L^2(b\mathbb H,\mu_{\mathbb H})$, the partial sums $s_N=\sum_{n=0}^Nq_n$ form a Cauchy sequence in
	$L^2(b\mathbb H,\mu_{\mathbb H})$.
	Hence, there exists a unique
	$f^*\in L^2(b\mathbb H,\mu_{\mathbb H})$
	such that $\lim_{N\to\infty}
	\|s_N-f^*\|_{L^2(b\mathbb H,\mu_{\mathbb H})}=0$. Note that $f\circ\delta_r
	=
	\sum_{n=0}^{\infty}r^nq_n$ belongs to
	in $L^2(b\mathbb H,\mu_{\mathbb H})$ for $0<r<1$. Thus, we have
	\[
	\|f\circ\delta_r-f^*\|_
	{L^2(b\mathbb H,\mu_{\mathbb H})}^2
	=
	\sum_{n=0}^{\infty}
	|r^n-1|^2
	\|q_n\|_{L^2(b\mathbb H,\mu_{\mathbb H})}^2 \to 0,
	\]
	where we have used the dominated convergence theorem. Then
	\begin{align*}
		\|f^*\|_{L^2(b\mathbb H,\mu_{\mathbb H})}^2
		=\lim_{r\uparrow1}
		\|f\circ\delta_r\|_
		{L^2(b\mathbb H,\mu_{\mathbb H})}^2
		=\lim_{r\uparrow1}
		\sum_{n=0}^{\infty}
		r^{2n}
		\|q_n\|_{L^2(b\mathbb H,\mu_{\mathbb H})}^2
		=\|f\|_{H^2(\mathbb H)}^2,
	\end{align*}
	which completes the proof. 
\end{proof}

\section{The multiplier algebra of $H^2(\mathbb H)$}\label{sec004_hex}

\noindent In this section, we prove that the Hardy space $H^2(\mathbb H)$ can be isometrically embedded in the classical Hardy space $H^2(\mathbb D^4)$ on the four dimensional polydisc $\mathbb D^4$. Next, we study the multiplier algebra of $H^2(\mathbb H)$ and show that it is isometrically isomorphic to the Hilbert space $H^\infty(\mathbb H)$. Moreover, we present a boundary realization theorem for $H^\infty(\mathbb H)$. 

\begin{proposition}\label{prop:isometric-embedding}
	The map $J_{\mathbb H}:H^2(\mathbb{H})\to H^2(\mathbb{D}^4)$ defined as
	\[
	J_{\mathbb H}f=f\circ\Psi_{\mathbb H}.
	\]
	is an isometric embedding. 
\end{proposition}

\begin{proof}
	Let $f\in H^2(\mathbb{H})$ and denote by $F=J_{\mathbb H}f$. Clearly, the map $F$ is holomorphic on $\mathbb D^4$. Fix $0<r<1$. It follows from Theorem \ref{prop:Psi-interior} that $\Psi_{\mathbb H}(\overline{\mathbb{D}}\times (r\overline{\mathbb{D}})^2 \times \overline{\mathbb D})\subset\mathbb{H}$. Thus, the map $(\xi_1, \xi_2,\lambda_1,\lambda_2) \mapsto
	F(\xi_1, r\xi_2, r\lambda_1, \lambda_2)$ extends continuously to $\mathbb{T}^4$. It is easy to see that the classical Hardy norm on $\mathbb{D}^4$ coincides with
	\begin{equation}\label{mixed-Hardy-norm}
		\|F\|_{H^2(\mathbb{D}^4)}^2
		=
		\sup_{0<r<1}
		\int_{\mathbb{T}^4}
		|F(\xi_1,r\xi_2, r\lambda_1, \lambda_2)|^2\,\mathrm{d}m_4.
	\end{equation}
	Suppose $(\xi_1, \xi_2, \lambda_1, \lambda_2) \in \mathbb T^4$. Let $(a, x_1, x_2, x_3)=\Psi_{\mathbb H}(\xi_1, \xi_2, \lambda_1, \lambda_2)$.  Then
	\begin{align*}
		\Psi_{\mathbb H}(\xi_1, r\xi_2, r\lambda_1, \lambda_2)
		&=\left(
		\frac{r\lambda_1(\xi_1-1)}{2},
		\frac{r\xi_2(\xi_1+1)}{2},
		\frac{r\lambda_1\lambda_2(\xi_1+1)}{2},
		r^2\xi_1\xi_2\lambda_1\lambda_2
		\right)\\
		&=(ra, rx_1, rx_2, r^2x_3)\\
		&=\delta_r(a, x_1, x_2, x_3),
	\end{align*}
	where $\delta_r$ is the radial map $\delta_r(a, x_1, x_2, x_3)=(ra, rx_1, rx_2, r^2x_3)$ for $0 \leq r <1$. Thus,
	\begin{align*}
		\|J_{\mathbb H}f\|_{H^2(\mathbb{D}^4)}^2
		&=
		\sup_{0<r<1}
		\int_{\mathbb{T}^4}
		|f\circ\Psi_{\mathbb H}(\xi_1, r\xi_2, r\lambda_1, \lambda_2)|^2\,\mathrm{d}m_4(\xi_1, \xi_2, \lambda_1, \lambda_2)\\
		&=
		\sup_{0<r<1}
		\int_{b\mathbb{H}}|f\circ\delta_r(a, x_1, x_2, x_3)|^2\,\mathrm{d}\mu_{\mathbb H}(a, x_1, x_2, x_3) \quad [\text{as $\mu_{\mathbb H}=(\Psi_{\mathbb H})_*m_4$}]\\
		&=
		\|f\|_{H^2(\mathbb{H})}^2,
	\end{align*}
	where the last equality holds due to Theorem \ref{thm:radial-realization-hex}. Thus, $J$ is an isometry.
\end{proof}

A holomorphic map $\varphi$ on $\mathbb H$ is called a
\textit{multiplier} of $H^2(\mathbb H)$ if $\varphi f\in H^2(\mathbb H)$
for every $f\in H^2(\mathbb H)$. The set of all multipliers is denoted
by $\operatorname{Mult}(H^2(\mathbb H))$. Equivalently,
\[
\operatorname{Mult}(H^2(\mathbb H))
=
\{\varphi\in\operatorname{Hol}(\mathbb H):
M_\varphi\in\mathcal B(H^2(\mathbb H))\},
\]
where $M_\varphi$ denotes the multiplication operator defined by
$M_\varphi f=\varphi f$. Equipped with the norm
\[
\|\varphi\|_{\operatorname{Mult}}
=
\|M_\varphi\|,
\]
the space $\operatorname{Mult}(H^2(\mathbb H))$ is called the
\textit{multiplier algebra} of $H^2(\mathbb H)$. We denote by
$H^\infty(\mathbb H)$ the space of bounded holomorphic functions on
$\mathbb H$ equipped with the supremum norm $\|.\|_{\infty,\mathbb H}$. Below, we present a complete description of the multiplier algebra of $H^2(\mathbb H)$.

\begin{theorem}\label{thm:multipliers-hex}
	The multiplier algebra of $H^2(\mathbb H)$ is isometrically
	isomorphic to $H^\infty(\mathbb H)$.
\end{theorem}

\begin{proof}
	Let $\varphi\in H^\infty(\mathbb H)$. For any
	$f\in H^2(\mathbb H)$, we have by Theorem \ref{thm:radial-realization-hex} that 
	\begin{align*}
		\|\varphi f\|_{H^2(\mathbb H)}^2
		&=
		\sup_{0<r<1}
		\int_{b\mathbb H}
		|\varphi\circ\delta_r|^2
		|f\circ\delta_r|^2\,d\mu_{\mathbb H}\\
		&\leq
		\|\varphi\|_{\infty,\mathbb H}^2
		\sup_{0<r<1}
		\int_{b\mathbb H}
		|f\circ\delta_r|^2\,d\mu_{\mathbb H}\\
		&=
		\|\varphi\|_{\infty,\mathbb H}^2
		\|f\|_{H^2(\mathbb H)}^2.
	\end{align*}
	Thus, $M_\varphi$ is bounded and $\|M_\varphi\|
	\leq
	\|\varphi\|_{\infty,\mathbb H}$. Conversely, let $\varphi\in\operatorname{Mult}(H^2(\mathbb H))$ and take $z \in \mathbb H$.
	Since $H^2(\mathbb H)$ is a reproducing kernel Hilbert space, the corresponding reproducing kernel $k_z$ satisfies $M_\varphi^*k_z=
	\overline{\varphi(z)}k_z$. Consequently,
	\[
	|\varphi(z)|
	=
	\frac{\|M_\varphi^*k_z\|_{H^2(\mathbb H)}}
	{\|k_z\|_{H^2(\mathbb H)}}
	\leq
	\|M_\varphi\| \quad \text{for every $z \in \mathbb H$ and so,} \quad \|\varphi\|_{\infty,\mathbb H}
	\leq
	\|M_\varphi\|.
	\]
	Hence, $\varphi\in H^\infty(\mathbb H)$. Combining the two inequalities, it follows that $\|M_\varphi\|_{\mathcal B(H^2(\mathbb H))}
	=
	\|\varphi\|_{\infty,\mathbb H}$, which completes the proof.
\end{proof}

\begin{corollary}\label{cor:multiplier-algebra-hex}
	The map
	\[
	\pi_{\infty, \mathbb H}:H^\infty(\mathbb H)
	\to
	\mathcal B(H^2(\mathbb H)),
	\qquad
	\pi_{\infty, \mathbb H}(\varphi)=M_\varphi
	\]
	is an isometric algebra isomorphism onto
	$\operatorname{Mult}(H^2(\mathbb H))$. Consequently,
	$H^\infty(\mathbb H)$ is a norm-closed unital operator algebra on
	$H^2(\mathbb H)$.
\end{corollary}

\begin{proof}
	By Theorem \ref{thm:multipliers-hex}, $\pi_{\infty, \mathbb H}$ is an
	isometric algebra isomorphism onto its range. Since
	$H^\infty(\mathbb H)$ is complete and $\|M_\varphi\|=\|\varphi\|_{\infty,\mathbb H}$, we have that the range of $\pi_{\infty, \mathbb H}$ is complete and hence is a norm-closed
	subalgebra of $\mathcal B(H^2(\mathbb H))$.
\end{proof}

We next establish a boundary realization theorem for
$H^\infty(\mathbb H)$.

\begin{theorem}\label{thm:boundary-Hinfty-hex}
	Let $f\in H^\infty(\mathbb H)$ and let
	$f^*\in L^2(b\mathbb H,\mu_{\mathbb H})$ be the unique boundary
	function defined as in Corollary \ref{cor:radial-boundary-H}. Then $f^*\in L^\infty(b\mathbb H,\mu_{\mathbb H})$ and $\|f^*\|_{\infty,b\mathbb H}
	=
	\|f\|_{\infty,\mathbb H}$. In particular, the boundary value map
	\[
	\tau_{\infty, \mathbb H}:
	H^\infty(\mathbb H)
	\to
	L^\infty(b\mathbb H,\mu_{\mathbb H}),
	\qquad
	f\mapsto f^*,
	\]
	is an isometric algebra embedding.
\end{theorem}

\begin{proof}
	Define $f_r(a,x_1,x_2,x_3)
	=
	f(ra,rx_1,rx_2,r^2x_3)$ for $(a,x_1,x_2,x_3)\in b\mathbb H$ and $0<r<1$.
	Since $\delta_r(\overline{\mathbb H})\subseteq\mathbb H$, the
	function $f_r$ is continuous on $\overline{\mathbb H}$. Moreover, $\|f_r\|_{L^\infty(b\mathbb H,\mu_{\mathbb H})}
	\leq \|f\|_{\infty,\mathbb H}$ for $0<r<1$. In particular, $f_r$ is uniformly bounded in $L^\infty(b\mathbb H,\mu_{\mathbb H})$. Since $H^\infty(\mathbb H)\subseteq H^2(\mathbb H)$, it follows
	from Corollary \ref{cor:radial-boundary-H} that there exists a
	unique $f^*\in L^2(b\mathbb H,\mu_{\mathbb H})$ such that
	\[
	\lim_{r\uparrow1}
	\|f_r-f^*\|_{L^2(b\mathbb H,\mu_{\mathbb H})}=0.
	\]
	We first show that $f^*\in L^\infty(b\mathbb H,\mu_{\mathbb H})$.
	By the $L^2$-convergence, there exists a sequence
	$r_k\uparrow1$ such that $f_{r_k}\to f^*$ a.e. on $b\mathbb H$. Since $|f_{r_k}(a,x_1,x_2,x_3)| \leq \|f\|_{\infty,\mathbb H}$, we have that
	$|f^*(a,x_1,x_2,x_3)|
	\leq \|f\|_{\infty,\mathbb H}$ for a.e. $(a,x_1,x_2,x_3)\in b\mathbb H$. Consequently,
	\[
	f^*\in L^\infty(b\mathbb H,\mu_{\mathbb H})
	\quad\text{and}\quad
	\|f^*\|_{\infty,b\mathbb H}
	\leq \|f\|_{\infty,\mathbb H}.
	\]
	We next show the reverse inequality. Let
	$h\in H^2(\mathbb H)$ and let
	$h^*\in L^2(b\mathbb H,\mu_{\mathbb H})$ be its radial
	boundary function, as given by Corollary
	\ref{cor:radial-boundary-H}. We claim that
	\[
	(fh)^*=f^*h^*.
	\]
	It follows from Theorem
	\ref{thm:multipliers-hex} that $fh\in H^2(\mathbb H)$. Moreover, $(fh)_r=f_rh_r$, where $h_r=h\circ\delta_r$ for $0<r<1$. Since $\mu_{\mathbb H}(b\mathbb H)=1$, we have
	\[
	\|g\|_{L^1(b\mathbb H,\mu_{\mathbb H})}
	\leq
	\|g\|_{L^2(b\mathbb H,\mu_{\mathbb H})}
	\]
	for every $g\in L^2(b\mathbb H,\mu_{\mathbb H})$. 
	Using this and the Cauchy-Schwarz inequality, we have
	\begin{align*}
		\|f_rh_r-f^*h^*\|_{L^1(b\mathbb H,\mu_{\mathbb H})}
		&\leq
		\|f_r(h_r-h^*)\|_{L^1}
		+
		\|(f_r-f^*)h^*\|_{L^1}\\
		&\leq
		\|f_r\|_{L^\infty}\|h_r-h^*\|_{L^1}
		+
		\|f_r-f^*\|_{L^2}\|h^*\|_{L^2}\\
		&\leq
		\|f_r\|_{L^\infty}\|h_r-h^*\|_{L^2}
		+
		\|f_r-f^*\|_{L^2}\|h^*\|_{L^2}\\
		&\leq
		\|f\|_{\infty,\mathbb H}
		\|h_r-h^*\|_{L^2}
		+
		\|f_r-f^*\|_{L^2}\|h^*\|_{L^2}.
	\end{align*}
	By Corollary \ref{cor:radial-boundary-H}, both
	$h_r\to h^*$ and $f_r\to f^*$ in $L^2(b\mathbb H,\mu_{\mathbb H})$
	as $r\uparrow1$. Therefore,
	\[
	\lim_{r\uparrow1}
	\|f_rh_r-f^*h^*\|_{L^1(b\mathbb H,\mu_{\mathbb H})}
	=0.
	\]
	Since $fh\in H^2(\mathbb H)$, another
	application of Corollary \ref{cor:radial-boundary-H} gives
	\[
	\lim_{r\uparrow1}
	\|(fh)_r-(fh)^*\|_{L^2(b\mathbb H,\mu_{\mathbb H})}
	=0 \quad \text{and so,}\quad 
	\lim_{r\uparrow1}
	\|(fh)_r-(fh)^*\|_{L^1(b\mathbb H,\mu_{\mathbb H})}
	=0.
	\]
	As $(fh)_r=f_rh_r$, we have by uniqueness of the limit in
	$L^1(b\mathbb H,\mu_{\mathbb H})$ that $(fh)^*=f^*h^*$. In particular, taking $h\in H^\infty(\mathbb H)$, we have
	\[
	\tau_{\infty,\mathbb H}(fh)
	=
	\tau_{\infty,\mathbb H}(f)
	\tau_{\infty,\mathbb H}(h).
	\]
	Thus, the boundary value map
	$\tau_{\infty,\mathbb H}$ is multiplicative. Its linearity
	follows directly from the uniqueness of radial boundary
	functions. Hence, $\tau_{\infty,\mathbb H}$ is an algebra
	homomorphism. It remains to prove that $\tau_{\infty,\mathbb H}$ is an
	isometry. Let $h\in H^2(\mathbb H)$. By Corollary \ref{cor:radial-boundary-H} and the identity $(fh)^*=f^*h^*$ established
	above, we have
	\begin{align*}
		\|M_fh\|_{H^2(\mathbb H)}
		=\|fh\|_{H^2(\mathbb H)}
		=\|(fh)^*\|_{L^2(b\mathbb H,\mu_{\mathbb H})}
		&=\|f^*h^*\|_{L^2(b\mathbb H,\mu_{\mathbb H})}\\
		&\leq
		\|f^*\|_{\infty,b\mathbb H}
		\|h^*\|_{L^2(b\mathbb H,\mu_{\mathbb H})}\\
		&=
		\|f^*\|_{\infty,b\mathbb H}
		\|h\|_{H^2(\mathbb H)}.
	\end{align*}
	Taking the supremum over $h\in H^2(\mathbb H)$ with
	$\|h\|_{H^2(\mathbb H)}\leq1$, it follows that $\|M_f\|_{\mathcal B(H^2(\mathbb H))}
	\leq \|f^*\|_{\infty,b\mathbb H}$. It now follows from  Theorem \ref{thm:multipliers-hex} that
	\[
	\|f\|_{\infty,\mathbb H}
	=
	\|M_f\|_{\mathcal B(H^2(\mathbb H))}
	\leq
	\|f^*\|_{\infty,b\mathbb H}.
	\]
	Combining this inequality with the reverse inequality
	proved above, we conclude that $\|f^*\|_{\infty,b\mathbb H}
	=
	\|f\|_{\infty,\mathbb H}$. Therefore, $\tau_{\infty,\mathbb H}$ is an isometric algebra homomorphism. In particular, it is injective and
	hence an isometric algebra embedding.
\end{proof}

The following commutative diagram summarizes the results in this section. 

\begin{center}
	\begin{tikzpicture}[>=stealth]
		
		% Top level
		\node (Hinf) at (0,4) {$H^\infty(\mathbb H)$};
		\node (Linf) at (7,4) {$L^\infty(b\mathbb H,\mu_{\mathbb H})$};
		
		% Middle level
		\node (H2P) at (0,2) {$H^2(\mathbb H)$};
		\node (L2P) at (7,2) {$L^2(b\mathbb H,\mu_{\mathbb H})$};
		
		% Bottom level
		\node (H2D) at (0,0) {$H^2(\mathbb D^4)$};
		\node (L2T) at (7,0) {$L^2(\mathbb T^4, m_4)$};
		
		% Top horizontal boundary map
		\draw[->]
		(Hinf) -- node[above] {$f\mapsto f^*$} (Linf);
		
		% Middle horizontal boundary map
		\draw[->]
		(H2P) -- node[above] {$f\mapsto f^*$} (L2P);
		
		% Bottom horizontal boundary map
		\draw[->]
		(H2D) -- node[above] {\small{radial boundary map}} (L2T);
		
		% Left vertical inclusions / embeddings
		\draw[->, right=20]
		(Hinf) to
		node[left] {$f\mapsto f$} (H2P);
		
		\draw[->]
		(H2P) -- node[left] {$f\mapsto f\circ\Psi_{\mathbb H}$} (H2D);
		
		% Right vertical inclusions / embeddings
		\draw[->, left=20]
		(Linf) to
		node[right] {$g\mapsto g$} (L2P);
		
		\draw[->]
		(L2P) -- node[right] {$g\mapsto g\circ\Psi$} (L2T);
		
	\end{tikzpicture}
\end{center}

\section{Brown-Halmos characterization and commutant lifting for $\mathbb H$}\label{sec005_hex}

\noindent In this section, we study Toeplitz operators on $H^2(\mathbb H)$ and the coordinate multiplication operators $(T_a,T_{x_1},T_{x_2},T_{x_3})$ on $H^2(\mathbb H)$. We first establish the Brown--Halmos identities characterizing Toeplitz operators. We then study the commutant of the coordinate multiplication operators and their joint reducing subspaces. Finally, we establish a commutant lifting theorem for $\mathbb H$-isometries. Let
\[
P_{H^2(\mathbb H)}:
L^2(b\mathbb H,\mu_{\mathbb H})
\to H^2(\mathbb H)
\]
denote the orthogonal projection. For $\varphi\in L^\infty(b\mathbb H,\mu_{\mathbb H})$, let
$
M_\varphi:
L^2(b\mathbb H,\mu_{\mathbb H})
\to L^2(b\mathbb H,\mu_{\mathbb H})
$
be the multiplication operator defined by
\[
M_\varphi(f)(a,x_1,x_2,x_3)
=
\varphi(a,x_1,x_2,x_3)f(a,x_1,x_2,x_3)
\]
for every $f\in L^2(b\mathbb H,\mu_{\mathbb H})$.
Evidently, $M_\varphi$ is a bounded normal operator satisfying
$\|M_\varphi\|=\|\varphi\|_{\infty,b\mathbb H}$ and $M_\varphi^*=M_{\overline{\varphi}}$. For $\varphi\in L^\infty(b\mathbb H,\mu_{\mathbb H})$, the Toeplitz operator with symbol $\varphi$ is the bounded operator $T_\varphi:H^2(\mathbb H)\to H^2(\mathbb H)$ defined by
\[
T_\varphi
=
P_{H^2(\mathbb H)}M_\varphi|_{H^2(\mathbb H)}.
\]
The coordinate functions $a,x_1,x_2,x_3$ are holomorphic polynomials on $\mathbb H$. Since the polynomial algebra is invariant under multiplication by each of the coordinate functions and is dense in $H^2(\mathbb H)$, it follows that $H^2(\mathbb H)$ is jointly invariant under $M_a, M_{x_1}, M_{x_2}$ and $M_{x_3}$. Consequently, the corresponding Toeplitz operators
$T_a,T_{x_1},T_{x_2}$ and $T_{x_3}$ are simply the restrictions of the multiplication operators $M_a,M_{x_1},M_{x_2}$ and $M_{x_3}$, respectively. The multiplication operators associated with $H^2(\mathbb H)$ and
$L^2(b\mathbb H,\mu_{\mathbb H})$ provide natural examples of isometries and unitaries, respectively, associated with $\mathbb H$. This naturally leads to the definitions of $\mathbb H$-isometries and $\mathbb H$-unitaries, which we recall next from \cite{PalTomarI}.

\begin{definition}
	A commuting quadruple $(A, X_1, X_2, X_3)$ of operators acting on a Hilbert space $\mathcal H$ is said to be 
	\begin{enumerate}[leftmargin=*]
		\item an \textit{$\mathbb H$-unitary} if its joint spectrum $\sigma_T(A, X_1, X_2, X_3) \subseteq b \mathbb H$;
		\item an \textit{$\mathbb H$-isometry} if it is the restriction of a $\mathbb H$-unitary to a joint invariant subspace.
	\end{enumerate}
\end{definition}

To begin with, we prove that the coordinate multipliers acting on $L^2(b\mathbb H, \mu_{\mathbb H})$ and $H^2(\mathbb H)$ are the canonical examples of $\mathbb H$-unitaries and $\mathbb H$-isometries, respectively. 

\begin{lemma}\label{lem_mini}
	The quadruple $(M_a,M_{x_1},M_{x_2},M_{x_3})$ acting on
	$L^2(b\mathbb H,\mu_{\mathbb H})$ is an $\mathbb H$-unitary with its joint spectrum
	$\sigma_T(M_a,M_{x_1},M_{x_2},M_{x_3})=b\mathbb H$. Moreover, $(T_a,T_{x_1},T_{x_2},T_{x_3})$ acting on $H^2(\mathbb H)$ is an $\mathbb H$-isometry such that $T_{x_3}$ is a pure isometry and $\sigma_T(T_a,T_{x_1},T_{x_2},T_{x_3})
	=\overline{\mathbb H}$. Also, $(M_a,M_{x_1},M_{x_2},M_{x_3})$ acting on
	$L^2(b\mathbb H,\mu_{\mathbb H})$ is the minimal normal extension
	of $(T_a,T_{x_1},T_{x_2},T_{x_3})$ acting on $H^2(\mathbb H)$.
\end{lemma}

\begin{proof}
	Since $(M_a,M_{x_1},M_{x_2},M_{x_3})$ is a commuting tuple
	of multiplication operators on $L^2(b\mathbb H,\mu_{\mathbb H})$,
	its Taylor joint spectrum coincides with the joint essential range
	of the coordinate maps $(a,x_1,x_2,x_3)$. As $\mu_{\mathbb H}$
	is supported on $b\mathbb H$, we have
	\[
	\sigma_T(M_a,M_{x_1},M_{x_2},M_{x_3})
	\subseteq b\mathbb H.
	\]
	Hence, $(M_a,M_{x_1},M_{x_2},M_{x_3})$ is an
	$\mathbb H$-unitary. Let $z=(a_0,x_{1,0},x_{2,0},x_{3,0})\in b\mathbb H$ and let
	$U$ be any open neighbourhood of $z$ in $\mathbb C^4$.
	Since the coordinate map $(a,x_1,x_2,x_3)$ is the identity
	on $b\mathbb H$, we have that $(a,x_1,x_2,x_3)^{-1}(U)=U\cap b\mathbb H$. We claim that
	\[
	\mu_{\mathbb H}(U\cap b\mathbb H)>0.
	\]
	Recall that $\mu_{\mathbb H}=(\Psi_{\mathbb H})_*m_4$,
	where $\Psi_{\mathbb H}:\mathbb T^4\to b\mathbb H$ is
	continuous and surjective. Therefore, $\mu_{\mathbb H}(U\cap b\mathbb H)
	=m_4\bigl(\Psi_{\mathbb H}^{-1}(U\cap b\mathbb H)\bigr)$. Since $\Psi_{\mathbb H}$ is continuous,
	$\Psi_{\mathbb H}^{-1}(U\cap b\mathbb H)$ is an open subset
	of $\mathbb T^4$. Also, it is nonempty because
	$\Psi_{\mathbb H}$ is surjective and
	$U\cap b\mathbb H\neq\varnothing$.
	Every nonempty open subset of $\mathbb T^4$ contains a
	nonempty basic open set of the form
	$
	I_1\times I_2\times I_3\times I_4,
	$
	where each $I_j$ is a nonempty open arc of $\mathbb T$.
	Since the normalized Haar measure on $\mathbb T^4$ is the
	product of the normalized arc-length measures on the four
	copies of $\mathbb T$, we have that $m_4(I_1\times I_2\times I_3\times I_4)
	=
	\prod_{j=1}^4 m_1(I_j)>0$. Consequently,
	\[
	m_4\bigl(\Psi_{\mathbb H}^{-1}(U\cap b\mathbb H)\bigr)>0
	\quad \text{and hence}
	\quad 
	\mu_{\mathbb H}(U\cap b\mathbb H)>0.
	\]
	Thus, $\mu_{\mathbb H}\bigl((a,x_1,x_2,x_3)^{-1}(U)\bigr)>0$ for every open neighbourhood $U$ of $z$. Therefore, we have that $z\in\operatorname{ess\,ran}(a,x_1,x_2,x_3)$ and so,
	$
	\sigma_T(M_a,M_{x_1},M_{x_2},M_{x_3})=b\mathbb H.
	$
	Since $(T_a,T_{x_1},T_{x_2},T_{x_3})$ is the restriction of
	$(M_a,M_{x_1},M_{x_2},M_{x_3})$ to a joint invariant subspace,
	the quadruple $(T_a,T_{x_1},T_{x_2},T_{x_3})$ is an
	$\mathbb H$-isometry. In particular,
	$
	\sigma_T(T_a,T_{x_1},T_{x_2},T_{x_3})
	\subseteq\overline{\mathbb H}.
	$
	Moreover, since $|x_3|=1$ on $b\mathbb H$, the operator
	$M_{x_3}$ is unitary on $L^2(b\mathbb H,\mu_{\mathbb H})$.
	It follows that $T_{x_3}$ is an isometry on $H^2(\mathbb H)$. Next, we show that
	\[
	\sigma_T(T_a,T_{x_1},T_{x_2},T_{x_3})
	=\overline{\mathbb H}.
	\]
	Since $H^2(\mathbb H)$ is a reproducing kernel Hilbert space,
	for every $z=(a,x_1,x_2,x_3)\in\mathbb H$, there exists a
	reproducing kernel $k_z\in H^2(\mathbb H)$ that satisfies $f(z)=\langle f,k_z\rangle$ for $f\in H^2(\mathbb H))$. Consequently,
	$
	\langle f,T_a^*k_z\rangle
	=\langle T_af,k_z\rangle
	=(af)(z)
	=a f(z)
	=\langle f,\overline a\,k_z\rangle
	$
	for every $f\in H^2(\mathbb H)$. Hence, $T_a^*k_z=\overline a\,k_z$. Similarly, $T_{x_1}^*k_z=\overline{x}_1\,k_z,
	T_{x_2}^*k_z=\overline{x}_2\,k_z$ and $T_{x_3}^*k_z=\overline{x}_3\,k_z$. Thus, $k_z$ is a joint eigenvector of the commuting quadruple
	$(T_a^*,T_{x_1}^*,T_{x_2}^*,T_{x_3}^*)$ corresponding to
	the joint eigenvalue $(\overline a, \overline{x}_1,\overline{x}_2,\overline{x}_3)$. It now follows that $\overline{\mathbb H}
	\subseteq \sigma_T(T_a^*,T_{x_1}^*,T_{x_2}^*,T_{x_3}^*)$ and thus $\sigma_T(T_a,T_{x_1},T_{x_2},T_{x_3})
	=\overline{\mathbb H}$. We now show that $T_{x_3}$ is pure. For $k\geq0$, define
	\[
	H_{\geq k}
	=
	\bigoplus_{n\geq k}\mathcal Q_n,
	\]
	where $\mathcal Q_n=
	\operatorname{span}\left\{
	a^jx_1^kx_2^\ell x_3^m:
	j+k+\ell+2m=n
	\right\}$. Since multiplication by $x_3$ increases the weighted degree
	by $2$, we have that $T_{x_3}^n(H^2(\mathbb H))
	\subseteq H_{\geq 2n}$ for  $n\geq0$. As $T_{x_3}$ is an isometry, it follows for every $f\in H^2(\mathbb H)$ that
	\[
	\|T_{x_3}^{*n}f\|
	=\|P_{\operatorname{Ran}T_{x_3}^n}f\|
	\leq \|P_{H_{\geq 2n}}f\|
	\to 0
	\]
	as $n\to\infty$. Hence, $T_{x_3}$ is a pure isometry. It is only left to show the minimality of the normal extension.
	Let
	\[
	\mathcal M=
	\overline{\operatorname{span}}
	\left\{
	M_a^{*i}M_{x_1}^{*j}M_{x_2}^{*k}M_{x_3}^{*\ell}f:
	f\in H^2(\mathbb H),\
	i,j,k,\ell\geq0
	\right\},
	\]
	where the closure is taken in
	$L^2(b\mathbb H,\mu_{\mathbb H})$.
	Since the polynomial algebra is dense in $H^2(\mathbb H)$,
	we have that
	$\mathcal M=
	\overline{\operatorname{span}}
	\left\{
	M_a^{*i}M_{x_1}^{*j}M_{x_2}^{*k}M_{x_3}^{*\ell}q:
	q\in\mathbb C[a,x_1,x_2,x_3],\
	i,j,k,\ell\geq0
	\right\}$. Note that $M_a^{*i}M_{x_1}^{*j}M_{x_2}^{*k}M_{x_3}^{*\ell}q=
	\overline a^{\,i}\,
	\overline{x_1}^{\,j}\,
	\overline{x_2}^{\,k}\,
	\overline{x_3}^{\,\ell}\,q$ for every $q \in \mathbb C[a, x_1, x_2, x_3]$. Therefore, $
	\mathcal M
	=
	\overline{\mathcal A}^{\,L^2(b\mathbb H,\mu_{\mathbb H})},
	$
	where
	\[
	\mathcal A=
	\operatorname{span}
	\left\{
	a^i x_1^j x_2^k x_3^\ell
	\overline a^{\,m}
	\overline{x_1}^{\,n}
	\overline{x_2}^{\,p}
	\overline{x_3}^{\,q}:
	i,j,k,\ell,m,n,p,q\geq0
	\right\}.
	\]
	The algebra $\mathcal A$ contains the constants, is closed
	under complex conjugation and separates the points of
	$b\mathbb H$. Hence, by the complex Stone-Weierstrass theorem,
	$
	\overline{\mathcal A}^{\,\|\cdot\|_\infty}
	=
	C(b\mathbb H).
	$
	Since $\mu_{\mathbb H}$ is a finite regular Borel measure on
	$b\mathbb H$, the space $C(b\mathbb H)$ is dense in
	$L^2(b\mathbb H,\mu_{\mathbb H})$. It follows that
	$
	\mathcal M
	=
	\overline{\mathcal A}^{\,L^2(b\mathbb H,\mu_{\mathbb H})}
	=
	L^2(b\mathbb H,\mu_{\mathbb H})$. Thus, $(M_a,M_{x_1},M_{x_2},M_{x_3})$ is the minimal normal extension of $(T_a,T_{x_1},T_{x_2},T_{x_3})$. The proof is now complete.
\end{proof}

It is well-known that any bounded operator $M$ on $L^2(\mathbb T)$ that belongs to the commutant of the operator $M_z$ on $L^2(\mathbb T)$ is of the form $M=M_\varphi$ for some $\varphi \in L^\infty(\mathbb T)$. Similar results have been established for the symmetrized bidisc \cite{BhattacharyyaDasSau2021}, the tetrablock \cite{JainKumarMalPramanick2026} (also see \cite{KumarMalPramanick2025}) and the pentablock \cite{SuTomar}. We prove an analogue of this result for the domain $\mathbb H$. 
\begin{lemma}\label{lem_402}
	Let $M$ be a bounded operator on $L^2(b \mathbb H, \mu_{\mathbb H})$. If $M$ commutes with $M_a, M_{x_1}, M_{x_2}$ and $M_{x_3}$, then there exists $\varphi \in L^\infty(b \mathbb H, \mu_{\mathbb H})$ such that $M=M_{\varphi}$.
\end{lemma}

\begin{proof}
	We have by  Lemma \ref{lem_mini} that $(M_a, M_{x_1}, M_{x_2}, M_{x_3})$ is a commuting quadruple consisting of normal operators such that $\sigma_T(M_a, M_{x_1}, M_{x_2}, M_{x_3})=b \mathbb{H}$. By the spectral theorem, the von Neumann algebra $W^*(M_a, M_{x_1}, M_{x_2}, M_{x_3})$ generated by $M_a, M_{x_1}, M_{x_2}, M_{x_3}$ is the maximal abelian von Neumann algebra such that
	$
	W^*(M_a, M_{x_1}, M_{x_2}, M_{x_3})=\{M_\varphi: \varphi \in L^\infty(b \mathbb H, \mu_{\mathbb H})\}.
	$
\end{proof}

We now explore the Brown-Halmos type characterization of Toeplitz operators on $H^2(\mathbb H)$. To do so, we recall from Theorem 3.5 in \cite{PalTomarI} that a commuting quadruple $(M_a, M_{x_1}, M_{x_2}, M_{x_3})$ is an $\mathbb H$-unitary if and only if 
\begin{align}\label{eqn_601}
	M_a^*M_a+M_{x_1}^*M_{x_1}=I, \ \ M_{x_1}=M_{x_2}^*M_{x_3}, \ \ \|M_{x_1}\| \leq 1,  \ \ \text{and} \ \ M_{x_3}^*M_{x_3}=I.
\end{align}
Let $\varphi \in L^\infty(b\mathbb H,\mu_{\mathbb H})$ and let $T = T_\varphi$ be the associated Toeplitz operator. For
$f,g\in H^2(\mathbb H)$, we have 
\begin{align*}
	\langle T_{x_3}^*T_\varphi T_{x_3}f,g\rangle
	=\langle T_\varphi T_{x_3}f,T_{x_3}g\rangle
	&=\langle P_{H^2(\mathbb H)}M_\varphi T_{x_3}f,T_{x_3}g\rangle\\
	&=\langle M_\varphi M_{x_3}f,M_{x_3}g\rangle\\
	&=\langle M_{x_3}M_\varphi f,M_{x_3}g\rangle\\
	&=\langle M_\varphi f,g\rangle
	\qquad [\text{since }M_{x_3}^*M_{x_3}=I]\\
	&=\langle P_{H^2(\mathbb H)}M_\varphi f,g\rangle\\
	&=\langle T_\varphi f,g\rangle.
\end{align*}
Similarly,
\begin{align*}
	\langle T_{x_1}^*T_\varphi T_{x_3}f,g\rangle
	=\langle T_\varphi T_{x_3}f,T_{x_1}g\rangle
	&=\langle P_{H^2(\mathbb H)}M_\varphi T_{x_3}f,T_{x_1}g\rangle\\
	&=\langle M_\varphi M_{x_3}f,M_{x_1}g\rangle\\
	&=\langle M_{x_3}M_\varphi f,M_{x_1}g\rangle\\
	&=\langle M_{x_1}^*M_{x_3}M_\varphi f,g\rangle\\
	&=\langle M_{x_2}M_\varphi f,g\rangle
	\qquad [\text{since }M_{x_1}=M_{x_2}^*M_{x_3}]\\
	&=\langle M_\varphi M_{x_2}f,g\rangle\\
	&=\langle P_{H^2(\mathbb H)}M_\varphi M_{x_2}f,g\rangle\\
	&=\langle T_\varphi T_{x_2}f,g\rangle.
\end{align*}
Finally,
\begin{align*}
	\langle T_a^*T_\varphi T_af,g\rangle
	=\langle T_\varphi T_af,T_ag\rangle
	&=\langle P_{H^2(\mathbb H)}M_\varphi T_af,T_ag\rangle\\
	&=\langle M_\varphi M_af,M_ag\rangle\\
	&=\langle M_aM_\varphi f,M_ag\rangle\\
	&=\langle M_a^*M_aM_\varphi f,g\rangle\\
	&=\langle M_\varphi f,g\rangle
	-\langle M_{x_1}^*M_{x_1}M_\varphi f,g\rangle\\
	&=\langle M_\varphi f,g\rangle
	-\langle M_{x_1}M_\varphi f,M_{x_1}g\rangle\\
	&=\langle M_\varphi f,g\rangle
	-\langle M_\varphi M_{x_1}f,M_{x_1}g\rangle \  [\text{as} \ M_a^*M_a+M_{x_1}^*M_{x_1}=I]\\
	&=\langle P_{H^2(\mathbb H)}M_\varphi f,g\rangle
	-\langle P_{H^2(\mathbb H)}M_\varphi M_{x_1}f,M_{x_1}g\rangle\\
	&=\langle T_\varphi f,g\rangle
	-\langle T_\varphi T_{x_1}f,T_{x_1}g\rangle\\
	&=\left\langle
	\left(T_\varphi-T_{x_1}^*T_\varphi T_{x_1}\right)f,g
	\right\rangle.
\end{align*}
Similarly, one can prove that $\langle T_a^*T_\varphi T_af,g\rangle
=\left\langle
\left(T_\varphi-T_{x_2}^*T_\varphi T_{x_2}\right)f,g
\right\rangle$. The latter holds since $M_a^*M_a+M_{x_2}^*M_{x_2}=I$. Therefore, we have proved that
\begin{align*}
	T_{x_3}^*T_\varphi T_{x_3}
	=T_\varphi, \quad  
	T_{x_1}^*T_\varphi T_{x_3}=T_\varphi T_{x_2}, \quad
	T_a^*T_\varphi T_a+T_{x_2}^*T_\varphi T_{x_2}
	=T_\varphi, \quad
	T_a^*T_\varphi T_a+T_{x_1}^*T_\varphi T_{x_1}
	=T_\varphi.
\end{align*}
The above identities are the natural Brown-Halmos type relations for Toeplitz operators on $H^2(\mathbb H)$. This motivates the following definition.

\begin{definition}
	We say that a bounded operator $T$ acting on $H^2(\mathbb H)$ satisfies the Brown-Halmos identities if the following holds:
	\begin{align}\label{eqn_BHR}
		T_{x_3}^*T T_{x_3}=T,\quad  
		T_{x_1}^*T T_{x_3}=T T_{x_2}, \quad
		T_{a}^*T T_{a}+T_{x_2}^*T T_{x_2}=T, \quad  
		T_a^*T T_a+T_{x_1}^*T T_{x_1}
		=T.
	\end{align}
\end{definition}

We now establish a commutant lifting lemma for operators acting on $H^2(\mathbb H)$ that satisfy the Brown-Halmos identities. The analogous statements appear in earlier works including the symmetrized bidisc \cite{BhattacharyyaDasSau2021}, the tetrablock \cite{JainKumarMalPramanick2026} and the pentablock \cite{SuTomar}.

\begin{lemma}\label{lem:lifting-hex}
	If $T\in\mathcal B(H^2(\mathbb H))$ satisfies the
	Brown-Halmos identities as in \eqref{eqn_BHR}, then there exists an operator
	$M\in\mathcal B(L^2(b\mathbb H,\mu_{\mathbb H}))$ commuting with the coordinate multiplication operators
	$M_a,M_{x_1},M_{x_2},M_{x_3}$ such that $\|M\|=\|T\|$ and
	$
	T=P_{H^2(\mathbb H)}M|_{H^2(\mathbb H)}.
	$
\end{lemma}

\begin{proof}
	It follows from \eqref{eqn_601} that
	\[
	\frac23M_a^*M_a+\frac13M_{x_1}^*M_{x_1}+\frac13M_{x_2}^*M_{x_2}+\frac13M_{x_3}^*M_{x_3}=I.
	\]
	Also, we have by \eqref{eqn_BHR} that 
	\begin{align}\label{eqn_4003_BH}
		\frac23T_a^*TT_a+\frac13T_{x_1}^*TT_{x_1}+\frac13T_{x_2}^*TT_{x_2}+\frac13T_{x_3}^*TT_{x_3}=T.
	\end{align}
	Define the commuting quadruple consisting of normal operators on $L^2(b \mathbb H, \mu_{\mathbb H})$ as
	\[
	(N_1, N_2, N_3, N_4)=\left(\sqrt{\frac 23} M_a, \frac1{\sqrt3}M_{x_1}, \frac1{\sqrt3}M_{x_2}, \frac1{\sqrt3}M_{x_3}\right).
	\]
	Clearly, $N_1^*N_1+N_2^*N_2+N_3^*N_3+N_4^*N_4=I$. Then the linear map
	\[
	\Phi_{\mathbb H}:\mathcal B(L^2(b\mathbb H,\mu_{\mathbb H}))
	\to
	\mathcal B(L^2(b\mathbb H,\mu_{\mathbb H})), \ 
	\Phi_{\mathbb H}(R)=N_1^*RN_1+N_2^*RN_2+N_3^*RN_3+N_4^*RN_4
	\]
	defines a contraction. Consider the operator $F=\begin{bmatrix}
		T & 0\\
		0 & 0
	\end{bmatrix}$ acting on  $L^2(b\mathbb H, \mu_{\mathbb H})=H^2(\mathbb H) \oplus H^2(\mathbb H)^\perp$.  Then $F$ is a norm-preserving extension of $T$. Define
	\[
	\Theta_m=\frac1m\sum_{k=1}^m\Phi_{\mathbb H}^k(F).
	\]
	Since $\Phi_{\mathbb H}$ is contractive, $\|\Theta_m\| \leq \|F\|
	=\|T\|$. Hence, $\{\Theta_m\}$ is bounded and therefore admits a weak operator topology limit point, say $M$. We claim that each $\Phi_{\mathbb H}^k(F)$ is a compression of $T$, that is,  $P_{H^2(\mathbb H)}\Phi_{\mathbb H}^k(F)|_{H^2(\mathbb H)} = T$. Let $R$ be a compression of $T$. Then for $f,g\in H^2(\mathbb H)$,
	\begin{align*}
		\langle P_{H^2(\mathbb H)}\Phi_{\mathbb H}(R)f,g\rangle
		=
		\langle \Phi_{\mathbb H}(R)f,g\rangle
		&=
		\sum_{j=1}^4
		\langle RN_jf,N_jg\rangle\\
		&=
		\sum_{j=1}^4
		\langle P_{H^2(\mathbb H)}RN_jf,N_jg\rangle \quad [\text{as each $N_jf, N_jg \in H^2(\mathbb H)$}]\\
		&=
		\sum_{j=1}^4
		\langle TN_jf, N_jg\rangle\\
		&=
		\frac{2}{3}\langle TM_af,M_ag\rangle + \frac{1}{3}\sum_{j=1}^3\langle TM_{x_j}f, M_{x_j}g\rangle\\
		&=
		\frac{2}{3}\langle TT_af, T_ag\rangle + \frac{1}{3}\sum_{j=1}^3\langle TT_{x_j}f, T_{x_j}g\rangle\\
		&=
		\langle Tf,g\rangle,
	\end{align*}
	where the last equality follows from \eqref{eqn_4003_BH}. Hence, $\Phi_{\mathbb H}(R)$ is again a compression of $T$. Since $F$ is a compression of $T$, one can show inductively that every $\Phi^k(F)$ is a compression of $T$ and thus, $P_{H^2(\mathbb H)}\Theta_m|_{H^2(\mathbb H)}=T$ for every $m \geq 1$. Letting $m \to \infty$, we have that 
	\[
	P_{H^2(\mathbb H)}M|_{H^2(\mathbb H)}=T.
	\]
	A routine computation shows that
	\[
	\Phi_{\mathbb H}(\Theta_m)-\Theta_m
	=
	\frac1m\bigl(\Phi_{\mathbb H}^{m+1}(F)-\Phi(F)\bigr).
	\]
	Since $\|\Phi_{\mathbb H}\| \leq 1$, we have that $\|\Phi^{m+1}(F)-\Phi(F)\|
	\leq
	2\|F\|$. Then
	\[
	\|\Phi_{\mathbb H}(\Theta_m)-\Theta_m\|
	\leq
	\frac{2\|F\|}{m}=\frac{2\|T\|}{m}
	\to 0.
	\]
	Since $\Phi_{\mathbb H}$ is weak operator continuous, we have that $\Phi_{\mathbb H}(M)=M$. Since $\Phi_{\mathbb H}(M)=M$ and each $N_j$ is normal, it follows from Proposition 2.4 in \cite{DavieJewell1977} that $MN_j=N_jM$ for $1 \leq j \leq 4$. Consequently, $M$ commutes with $M_a, M_{x_1}, M_{x_2}$ and $M_{x_3}$. Since $M$ is a limit point of the sequence $\{\Theta_m\}$ in the weak operator topology, $\|M\| \leq \liminf \|\Theta_m\| \leq \|T\|$ and thus, $\|M\|=\|T\|$. The proof is complete.
\end{proof}

Next, we present a characterization of Toeplitz operators on $H^2(\mathcal P)$.

\begin{theorem}\label{thm_BH}
	A bounded linear operator $T$ acting on $H^2(\mathbb H)$ is Toeplitz if and only if $T$ satisfies the Brown-Halmos identities as in \eqref{eqn_BHR}. 
\end{theorem}
\begin{proof}
	The discussion preceding \eqref{eqn_BHR} shows that the Brown-Halmos identities hold for every Toeplitz operator on $H^2(\mathbb H)$. We prove the converse. Let $T\in\mathcal B(H^2(\mathbb H))$ satisfy the Brown-Halmos identities. It follows from Lemma \ref{lem:lifting-hex} that there exists a
	bounded operator $M$ on $L^2(b\mathbb H,\mu_{\mathbb H})$ such that $\|M\|=\|T\|$, $M$ commutes with the coordinate multipliers $M_a, M_{x_1}, M_{x_2}, M_{x_3}$ and $T=P_{H^2(\mathbb H}M|_{H^2(\mathbb H)}$. Next, Lemma \ref{lem_402} guarantees the existence of a map
	$\varphi\in L^\infty(b\mathbb H,\mu_{\mathbb H})$ such that $M=M_\varphi$. Hence, $T
	=
	P_{H^2(\mathbb H)}M|_{H^2(\mathbb H)}
	=
	P_{H^2(\mathbb H)}M_\varphi|_{H^2(\mathbb H)}
	=
	T_\varphi$.
\end{proof}

We present some interesting consequences of the Brown-Halmos type characterization theorem.

\begin{corollary}\label{cor_406}
	If $T \in \mathcal{B}(H^2(\mathbb H))$ commutes with the coordinate multipication operators $T_a, T_{x_1}, T_{x_2}$ and $T_{x_3}$, then $T$ is a Toeplitz operator. 
\end{corollary}

\begin{proof}
	We have by Lemma \ref{lem_mini} that $(T_a, T_{x_1}, T_{x_2}, T_{x_3})$ is an $\mathbb H$-isometry. It follows from Theorem 4.4 in \cite{PalTomarI} that 
	\[
	T_a^*T_a+T_{x_1}^*T_{x_1}=I, \quad  T_a^*T_a+T_{x_2}^*T_{x_2}=I, \quad T_{x_1}=T_{x_2}^*T_{x_3} \quad \text{and} \quad T_{x_3}^*T_{x_3}=I.
	\]
	The commutativity hypothesis of $T$ with $T_a, T_{x_1}, T_{x_2}, T_{x_3}$ gives that $T$ satisfies the Brown-Halmos identities. The desired conclusion follows from Theorem \ref{thm_BH}.
\end{proof}

\begin{corollary}\label{cor_compact_toeplitz_hex}
	The only compact Toeplitz operator on $H^2(\mathbb H)$
	is the zero operator.
\end{corollary}

\begin{proof}
	Let $T_\varphi$ be a compact Toeplitz operator on
	$H^2(\mathbb H)$. We have by Theorem \ref{thm_BH} that
	$
	(T_{x_3}^*)^nT_\varphi T_{x_3}^n=T_\varphi
	$
	for every $n\geq 1$. Since $T_{x_3}$ is a pure isometry,
	we have
	\[
	\langle T_{x_3}^n f,g\rangle
	=
	\langle f,(T_{x_3}^*)^n g\rangle \to 0
	\]
	for all $f,g\in H^2(\mathbb H)$.
	Thus, $T_{x_3}^n f$ converges weakly to zero for every
	$f\in H^2(\mathbb H)$. By the compactness of $T_\varphi$, it follows that
	$\|T_\varphi T_{x_3}^n f\|\to 0$ for every $f\in H^2(\mathbb H)$.
	Now, for arbitrary $f,g\in H^2(\mathbb H)$, we have
	\begin{align*}
		\langle T_\varphi f,g\rangle
		=\langle (T_{x_3}^*)^nT_\varphi T_{x_3}^n f,g\rangle
		=\langle T_\varphi T_{x_3}^n f,T_{x_3}^n g\rangle
		\leq
		\|T_\varphi T_{x_3}^n f\|
		\|T_{x_3}^n g\|
		=
		\|T_\varphi T_{x_3}^n f\|\|g\|.
	\end{align*}
	Letting $n\to\infty$, $\langle T_\varphi f,g\rangle=0$
	for all $f,g\in H^2(\mathbb H)$ and thus, $T_\varphi=0$.
\end{proof}

Next, we investigate the joint reducing subspaces of the coordinate
multiplication operators $T_a,T_{x_1},T_{x_2}$ and $T_{x_3}$.
Recall that a closed subspace $\mathcal M$ of a Hilbert space
$\mathcal H$ is called a \textit{joint reducing subspace} for a
family of operators $\{T_\alpha\}_{\alpha\in I}\subseteq
\mathcal B(\mathcal H)$ if $\mathcal M$ reduces each $T_\alpha$,
that is,
\[
T_\alpha\mathcal M\subseteq\mathcal M
\quad\text{and}\quad
T_\alpha^*\mathcal M\subseteq\mathcal M
\]
for every $\alpha\in I$.

\begin{lemma}\label{lem:irreducibility-hex}
	The coordinate multiplication operators
	$T_a,T_{x_1},T_{x_2}$ and $T_{x_3}$ acting on
	$H^2(\mathbb H)$ have no nontrivial joint reducing subspaces.
\end{lemma}

\begin{proof}
	Let $\mathcal M$ be a joint reducing subspace for
	$(T_a,T_{x_1},T_{x_2},T_{x_3})$ and let $P_{\mathcal M}$
	denote the orthogonal projection of $H^2(\mathbb H)$ onto
	$\mathcal M$. Since $\mathcal M$ reduces each of the four
	coordinate multiplication operators, we have
	\[
	P_{\mathcal M}T_\nu=T_\nu P_{\mathcal M}
	\quad\text{and}\quad
	P_{\mathcal M}T_\nu^*=T_\nu^*P_{\mathcal M}
	\]
	for every $\nu\in\{a,x_1,x_2,x_3\}$. Recall that the constant function $1$ is the homogeneous component of degree zero in $H^2(\mathbb H)$. Since
	multiplication by each coordinate function raises the
	$(1,1,1,2)$-weighted degree, the orthogonality of the
	weighted homogeneous decomposition gives $T_a^*1=T_{x_1}^*1=T_{x_2}^*1=T_{x_3}^*1=0$. We claim that
	\[
	\ker T_a^*\cap\ker T_{x_1}^*
	\cap\ker T_{x_2}^*\cap\ker T_{x_3}^*
	=\mathbb C1.
	\]
	Let $g\in\ker T_a^*\cap\ker T_{x_1}^*
	\cap\ker T_{x_2}^*\cap\ker T_{x_3}^*$. Since $g$ belongs to the kernel of the adjoint of each
	coordinate multiplication operator, it follows that
	\[
	\langle g,af\rangle
	=\langle g,x_1f\rangle
	=\langle g,x_2f\rangle
	=\langle g,x_3f\rangle=0
	\]
	for every $f\in H^2(\mathbb H)$.
	Thus, $g$ is orthogonal to every non-constant
	monomial in $a,x_1,x_2,x_3$ since each such monomial
	is divisible by at least one of these coordinate functions. Let 
	$g=\sum_{n=0}^{\infty}g_n$ be the weighted homogeneous decomposition of
	$g\in H^2(\mathbb H)$, where $g_n\in\mathcal Q_n$. Since $g$ is orthogonal to every nonconstant monomial,
	we have $g\perp\mathcal Q_n$ for every $n\geq1$.
	By the orthogonality of the weighted homogeneous
	decomposition,
	\[
	0=\langle g,g_n\rangle=\|g_n\|^2
	\]
	for every $n\geq1$. Thus $g_n=0$ for every $n\geq1$
	and $g=g_0\in\mathcal Q_0=\mathbb C1$. Since $T_a^*1=T_{x_1}^*1=T_{x_2}^*1=T_{x_3}^*1=0$, we have that
	$
	\ker T_a^*\cap\ker T_{x_1}^*
	\cap\ker T_{x_2}^*\cap\ker T_{x_3}^*
	=\mathbb C1$. Set $h=P_{\mathcal M}1$. Using the commutativity of
	$P_{\mathcal M}$ with the adjoints of the coordinate
	multiplication operators, we have
	\[
	T_\nu^*h
	=T_\nu^*P_{\mathcal M}1
	=P_{\mathcal M}T_\nu^*1
	=0
	\]
	for every $\nu\in\{a,x_1,x_2,x_3\}$.
	Therefore, $h \in \ker T_a^*\cap\ker T_{x_1}^*
	\cap\ker T_{x_2}^*\cap\ker T_{x_3}^*$ and so, $P_{\mathcal M}1=c1$for some $c\in\mathbb C$. Since $P_{\mathcal M}$ commutes with all four coordinate
	multiplication operators, it commutes with every
	polynomial in these operators. Hence, 
	\[
	P_{\mathcal M}q
	=P_{\mathcal M}q(T_a,T_{x_1},T_{x_2},T_{x_3})1
	=q(T_a,T_{x_1},T_{x_2},T_{x_3})P_{\mathcal M}1
	=cq
	\]
	for every $q\in\mathbb C[a,x_1,x_2,x_3]$. Since the polynomial algebra
	$\mathbb C[a,x_1,x_2,x_3]$ is dense in
	$H^2(\mathbb H)$ and $P_{\mathcal M}$ is bounded,
	it follows that $P_{\mathcal M}=cI$ on $H^2(\mathbb H)$.
	As $P_{\mathcal M}$ is an orthogonal projection,
	we must have $c\in\{0,1\}$. Thus,
	\[
	\mathcal M=\{0\}
	\quad\text{or}\quad
	\mathcal M=H^2(\mathbb H).
	\]
	This completes the proof.
\end{proof}

The argument used in Lemma \ref{lem:irreducibility-hex} can also be used to characterize the operators that doubly commute
with the coordinate multiplication operators. To see this, let
$X\in\mathcal B(H^2(\mathbb H))$ satisfy $XT_\nu=T_\nu X$ and $XT_\nu^*=T_\nu^*X$ for $\nu\in\{a,x_1,x_2,x_3\}$. Taking $h=X1 \in H^2(\mathbb H)$, we get that $T_\nu^*h
=T_\nu^*X1
=XT_\nu^*1
=0$ for $\nu\in\{a,x_1,x_2,x_3\}$. The description of the joint kernel in the proof of
Lemma~\ref{lem:irreducibility-hex} therefore implies that
$h=c1$ for some $c\in\mathbb C$. Since $X$ commutes with
each coordinate multiplication operator, it follows that
$Xq
=Xq(T_a,T_{x_1},T_{x_2},T_{x_3})1
=q(T_a,T_{x_1},T_{x_2},T_{x_3})X1
=cq
$
for every $q\in\mathbb C[a,x_1,x_2,x_3]$.
The density of the polynomial algebra in $H^2(\mathbb H)$
and the boundedness of $X$ now give $X=cI$ on
$H^2(\mathbb H)$. We have thus arrive at the following
result.

\begin{theorem}\label{thm:double_commutant_hex}
	If $X\in\mathcal B(H^2(\mathbb H))$ doubly commute with
	$T_a,T_{x_1},T_{x_2}$ and $T_{x_3}$, then $X=cI$ for some $c\in\mathbb C$.
\end{theorem}

The Brown-Halmos characterization from Theorem \ref{thm_BH} suggests considering operators satisfying analogous relations with respect to arbitrary $\mathbb H$-isometries. This provides a framework for extending the results beyond the canonical coordinate multiplication operators on $H^2(\mathbb H)$. We begin with the following definition.

\begin{definition}
	Let $(A,X_1,X_2,X_3)$ be an $\mathbb H$-isometry acting on a Hilbert space $\mathcal H$. An operator $T\in\mathcal B(\mathcal H)$ is said to satisfy the \textit{Brown--Halmos relations} with respect to $(A,X_1,X_2,X_3)$ if
	\begin{align*}
		X_3^*TX_3=T, \quad 
		X_1^*TX_3=TX_2, \quad
		A^*TA+X_1^*TX_1=T,\quad
		A^*TA+X_2^*TX_2=T.
	\end{align*}
\end{definition}

We have by Theorem \ref{thm_BH} that an operator on $H^2(\mathbb H)$ is Toeplitz if and only if it satisfies the Brown-Halmos relations relative to the canonical $\mathbb H$-isometry $(T_a,T_{x_1},T_{x_2},T_{x_3})$. The lifting lemma, that is, Lemma \ref{lem:lifting-hex} used in its proof also suggests that these relations admit a formulation for general $\mathbb H$-isometries, which we prove in the following theorem.

\begin{theorem}[Commutant lifting theorem]\label{thm_CLT}
	Let $(T_1,T_2,T_3,T_4)$ be an $\mathbb H$-isometry acting on a Hilbert space $\mathcal H$ and let $(N_1,N_2,N_3,N_4)$ be its minimal $\mathbb H$-unitary extension on a Hilbert space $\mathcal K$ containing $\mathcal H$. An operator $T\in\mathcal B(\mathcal H)$ satisfies the Brown-Halmos relations with respect to $(T_1,T_2,T_3,T_4)$ if and only if there exists a unique operator $N\in\mathcal B(\mathcal K)$ such that
	\[
	T=P_{\mathcal H}N|_{\mathcal H},\qquad
	\|N\|=\|T\|,\qquad
	NN_j=N_jN,\quad 1\leq j\leq 4.
	\]
\end{theorem}

\begin{proof}
	We first establish the forward implication. Suppose that $T\in\mathcal B(\mathcal H)$ satisfies the Brown-Halmos relations relative to $(T_1,T_2,T_3,T_4)$. Then
	\begin{align}\label{eqn_hex_weighted_BH}
		\frac{2}{3}T_1^*TT_1+
		\frac{1}{3}T_2^*TT_2+
		\frac{1}{3}T_3^*TT_3+
		\frac{1}{3}T_4^*TT_4=T.
	\end{align}
	Since $(N_1, N_2, N_3, N_4)$ is an $\mathbb H$-unitary, it follows from Theorem 3.5 in \cite{PalTomarI} that
	\[
	\frac{2}{3}N_1^*N_1+
	\frac{1}{3}N_2^*N_2+
	\frac{1}{3}N_3^*N_3+
	\frac{1}{3}N_4^*N_4=I_{\mathcal K}.
	\]
	Following the proof of Lemma \ref{lem:lifting-hex} and applying Proposition 2.4 of \cite{DavieJewell1977} as in the proof, one can easily obtain an operator $N\in \mathcal B(\mathcal K)$ with the desired properties. Next,  assume that $N, N' \in\mathcal B(\mathcal K)$ both satisfy the conclusions of this theorem. Set $N''=N-N'$. Then $N''$ commutes with each $N_j$ and
	$P_{\mathcal H}N''|_{\mathcal H}=0$. Since each $N_j$ is normal, Fuglede's theorem implies that $N''$ also commutes with $N_j^*$ for $1\leq j \leq 4$.
	We claim that $N''$ vanishes identically. Let $h,k\in\mathcal H$ and let
	$\alpha=(\alpha_1,\alpha_2,\alpha_3,\alpha_4)\in\mathbb Z_+^4$. Then
	\begin{align*}
		\left\langle N''h,
		N_1^{*\alpha_1}N_2^{*\alpha_2}
		N_3^{*\alpha_3}N_4^{*\alpha_4}k
		\right\rangle
		=
		\left\langle
		N_1^{\alpha_1}N_2^{\alpha_2}
		N_3^{\alpha_3}N_4^{\alpha_4}N''h,k
		\right\rangle
			&=
		\left\langle
		N''N_1^{\alpha_1}N_2^{\alpha_2}
		N_3^{\alpha_3}N_4^{\alpha_4}h,k
		\right\rangle\\
		&=
		\left\langle
		P_{\mathcal H}N''N_1^{\alpha_1}N_2^{\alpha_2}
		N_3^{\alpha_3}N_4^{\alpha_4}h,k
		\right\rangle\\
		&=0,
	\end{align*}
	because $\mathcal H$ is jointly invariant under $N_j$ and the compression of $N''$ to $\mathcal H$ is zero. By the minimality of the $\mathbb H$-unitary extension, the space $\mathcal K$ is given by
	\[
	\mathcal K=
	\overline{\operatorname{span}}
	\left\{
	N_1^{*\alpha_1}N_2^{*\alpha_2}
	N_3^{*\alpha_3}N_4^{*\alpha_4}h:
	h\in\mathcal H,\ 
	\alpha\in\mathbb Z_+^4
	\right\}
	\]
	and so,  $N''N_1^{*\alpha_1}N_2^{*\alpha_2}
	N_3^{*\alpha_3}N_4^{*\alpha_4}h
	=
	N_1^{*\alpha_1}N_2^{*\alpha_2}
	N_3^{*\alpha_3}N_4^{*\alpha_4}N''h
	=0$ for every $h \in \mathcal H$. Hence, $N=N'$. For the converse, assume that there exists an operator $N\in\mathcal B(\mathcal K)$ satisfying
	$T=P_{\mathcal H}N|_{\mathcal H}, \|N\|=\|T\|$ and $NN_j=N_jN$ for $1\leq j\leq4$. Since $(N_1,N_2,N_3,N_4)$ is an $\mathbb H$-unitary, it follows from Theorem 3.5 in \cite{PalTomarI} that
	$N_1^*N_1+N_2^*N_2=I_{\mathcal K}, N_1^*N_1+N_3^*N_3=I_{\mathcal K},
	N_3^*N_4=N_2$ and $N_4^*N_4=I_{\mathcal K}$. Since $N$ commutes with each $N_j$, it satisfies the Brown-Halmos relations with respect to $(N_1, N_2, N_3, N_4)$. Next, one can write 
	\[
	N_j=
	\begin{bmatrix}
		T_j&*\\
		0&*
	\end{bmatrix}
	\quad (1\leq j\leq4) \quad \text{and} 
	\quad
	N=
	\begin{bmatrix}
		T&*\\
		*&*
	\end{bmatrix}
	\]
	with respect to the orthogonal decomposition $\mathcal K=\mathcal H\oplus\mathcal H^\perp$.
	Since $N$ satisfies the Brown-Halmos relations with respect to $(N_1, N_2, N_3, N_4)$ and $\mathcal H$ is jointly invariant under $N_1,N_2,N_3,N_4$, it follows that $T$ satisfies all the Brown-Halmos relations relative to $(T_1,T_2,T_3,T_4)$. The proof is now complete.
\end{proof}

\section{A weighted Hardy space on the hexablock}\label{sec006_hex}

\noindent In Section \ref{sec003_hex}, we constructed the Hardy space $H^2(\mathbb H)$ using the probability measure $\mu_{\mathbb H}$ on the distinguished boundary $b\mathbb H$, which was obtained as the pushforward of the normalized Haar measure on $\mathbb T^4$ under the parametrization $\Psi_{\mathbb H}$. We now introduce another natural probability measure on $b\mathbb H$ arising from its matrix representation. Recall from \cite{BiswasPalII} that
\[
b\mathbb H
=
\left\{
(u_{21},u_{11},u_{22},\det U):
U=(u_{ij})\in U(2)
\right\}.
\] 
Consider the map
\[
\pi_{\mathbb H}:U(2)\to b\mathbb H,
\qquad
\pi_{\mathbb H}(U)
=(u_{21},u_{11},u_{22},\det U).
\]
We show that $\pi_{\mathbb H}$ is a bijection. To see this,
let $U=(u_{ij}),V=(v_{ij})\in U(2)$ be such that $\pi_{\mathbb H}(U)=\pi_{\mathbb H}(V)$. Then
\[
u_{21}=v_{21},\quad
u_{11}=v_{11},\quad
u_{22}=v_{22},\quad
\det U=\det V.
\]
First, assume that $u_{21}\neq 0$. Then $u_{11}u_{22}-u_{21}u_{12}
=v_{11}v_{22}-v_{21}v_{12}$. Since $u_{11}=v_{11}$, $u_{21}=v_{21}$ and
$u_{22}=v_{22}$, it follows that $u_{21}u_{12}=u_{21}v_{12}$ and thus $u_{12}=v_{12}$. Next, let $u_{21}=0$. Then $v_{21}=0$ as well.
Since $U$ and $V$ are unitary, their first columns have
norm one and hence $|u_{11}|=|v_{11}|=1$. The orthogonality of the columns of $U$ and $V$ then gives $u_{12}=v_{12}=0$. In either case, we obtain $U=V$, proving that $\pi_{\mathbb H}$ is injective. Since $b\mathbb H$
is precisely the image of $U(2)$ under $\pi_{\mathbb H}$, the map $\pi_{\mathbb H}$ is also surjective. Therefore, $\pi_{\mathbb H}$ is a bijection. Since a continuous bijection from a compact space onto a Hausdorff space is a homeomorphism, it follows that $\pi_{\mathbb H}: U(2) \to b\mathbb H$ is a homeomorphism. Moreover, the inverse map can be explicitly described as
\[
\pi_{\mathbb H}^{-1}(a,x_1,x_2,x_3)
=
\begin{cases}
	\begin{pmatrix}
		x_1 & \dfrac{x_1x_2-x_3}{a}\\[5pt]
		a & x_2
	\end{pmatrix}, & a\neq 0,\\[14pt]
	\begin{pmatrix}
		x_1&0\\
		0&x_2
	\end{pmatrix}, & a=0.
\end{cases}
\]
Consequently, we arrive at the following result.

\begin{proposition}
	The map $\pi_{\mathbb H}: U(2) \to b\mathbb H$ defined as 
	\[
	\pi_{\mathbb H}(U)=(u_{21}, u_{11}, u_{22}, \det(U)), \qquad U=(u_{ij}) \in U(2)
	\]
	is a homeomorphism.
\end{proposition}

\smallskip

Let $m_{U(2)}$ denote the normalized Haar measure on $U(2)$. We define a probability measure $\sigma_{\mathbb H}$ on $b\mathbb H$ by $\sigma_{\mathbb H}=(\pi_{\mathbb H})_*m_{U(2)}$. Equivalently, for every $F\in C(b\mathbb H)$,
\[
\begin{aligned}
	\int_{b\mathbb H}F(a,x_1,x_2,x_3)\,
	d\sigma_{\mathbb H}
	=
	\int_{U(2)}
	F\circ \pi_{\mathbb H}(u_{21},u_{11},u_{22},\det U)\,
	dm_{U(2)}(U).
\end{aligned}
\]
We first establish the relationship between the two boundary measures.
\begin{proposition}
	The measures $\mu_{\mathbb H}$ and $\sigma_{\mathbb H}$ are mutually
	absolutely continuous. More precisely,
	\[
	\frac{d\sigma_{\mathbb H}}{d\mu_{\mathbb H}}
	(a,x_1,x_2,x_3)
	=
	\pi |a|\sqrt{1-|a|^2}
	\]
	almost everywhere on $b\mathbb H$.
\end{proposition}

\begin{proof}
	Consider the map $\Phi=\pi_{\mathbb H}^{-1}\circ\Psi_{\mathbb H}:
	\mathbb T^4\to U(2)$. We first recall the Haar integration formula
	on $U(2)$. Every matrix $U\in U(2)$ can be written as
	\[
	U=
	\begin{pmatrix}
		e^{i\alpha}\cos r & e^{i\beta}\sin r\\
		e^{i\gamma}\sin r & e^{i\delta}\cos r
	\end{pmatrix},
	\qquad 0\leq r\leq\frac{\pi}{2},
	\]
	where $\alpha+\delta-\beta-\gamma\equiv\pi\pmod{2\pi}$. In these coordinates, the normalized Haar measure on $U(2)$
	is given by
	\begin{align*}
		 \int_{U(2)}G(U)\,dm_{U(2)}(U)
		=
		\frac{1}{(2\pi)^3}
		\int_0^{\pi/2}\int_{[0,2\pi]^3}
		G\left(
		\begin{pmatrix}
			e^{i\alpha}\cos r & e^{i\beta}\sin r\\
			e^{i\gamma}\sin r & e^{i\delta}\cos r
		\end{pmatrix}
		\right)
		2\sin r\cos r\,
		d\alpha\,d\beta\,d\gamma\,dr
	\end{align*}
	for every $G\in C(U(2))$. We now show that the weighted pushforward of $m_4$ under
	$\Phi$ coincides with $m_{U(2)}$. Writing
	$\xi_1=e^{it}$, $\xi_2=e^{iu}$,
	$\lambda_1=e^{iv}$ and $\lambda_2=e^{iw}$, we have
	\[
	\Phi(\xi_1,\xi_2,\lambda_1,\lambda_2)
	=
	\begin{pmatrix}
		\displaystyle \frac{\xi_2(\xi_1+1)}{2} &  \displaystyle \frac{\xi_2\lambda_2(\xi_1-1)}{2}  \smallskip \\
		\displaystyle    \frac{\lambda_1(\xi_1-1)}{2} &  \displaystyle \frac{\lambda_1\lambda_2(\xi_1+1)}{2}
	\end{pmatrix}
	=
	\begin{pmatrix}
		\displaystyle e^{iu}\frac{e^{it}+1}{2}
		&
		\displaystyle e^{i(u+w)}\frac{e^{it}-1}{2}
		\\[6pt]
		\displaystyle e^{iv}\frac{e^{it}-1}{2}
		&
		\displaystyle e^{i(v+w)}\frac{e^{it}+1}{2}
	\end{pmatrix}.
	\]
	Putting $\theta=t/2$, this becomes
	$
	\Phi(\xi_1,\xi_2,\lambda_1,\lambda_2)
	=
	\begin{pmatrix}
		e^{i(u+\theta)}\cos\theta
		&
		ie^{i(u+w+\theta)}\sin\theta
		\\
		ie^{i(v+\theta)}\sin\theta
		&
		e^{i(v+w+\theta)}\cos\theta
	\end{pmatrix}.
	$
	Furthermore,
	\[
	\left|\frac{\lambda_1(\xi_1-1)}{2}\right|
	=\sin\theta,
	\qquad
	\sqrt{1-\left|\frac{\lambda_1(\xi_1-1)}{2}\right|^2}
	=|\cos\theta|.
	\]
	We next describe the changes of variables needed to express
	$\Phi$ in the above coordinates on $U(2)$. For
	$0\leq\theta\leq\pi/2$, put $r=\theta$ and define
	\[
	\alpha=u+\theta,\quad 
	\beta=u+w+\theta+\frac{\pi}{2},\quad  
	\gamma=v+\theta+\frac{\pi}{2},\quad 
	\delta=v+w+\theta.
	\]
	For $\pi/2\leq\theta\leq\pi$, put $r=\pi-\theta$ and define
	\[
	\alpha=u+\theta+\pi,\quad
	\beta=u+w+\theta+\frac{\pi}{2},\quad
	\gamma=v+\theta+\frac{\pi}{2},\quad 
	\delta=v+w+\theta+\pi.
	\]
	In each case, $\alpha+\delta-\beta-\gamma\equiv\pi\pmod{2\pi}$ and
	$
	\Phi(\xi_1,\xi_2,\lambda_1,\lambda_2)
	=
	U(\alpha,\beta,\gamma,r)$, where
	\[
	U(\alpha,\beta,\gamma,r)
	=
	\begin{pmatrix}
		e^{i\alpha}\cos r&e^{i\beta}\sin r\\
		e^{i\gamma}\sin r&
		-e^{i(\beta+\gamma-\alpha)}\cos r
	\end{pmatrix}.
	\]
	For each fixed $\theta$, both changes of variables from
	$(u,v,w)$ to $(\alpha,\beta,\gamma)$ admits the Jacobian
	matrix 
	\[
	\begin{pmatrix}
		1&0&0\\
		1&0&1\\
		0&1&0
	\end{pmatrix},
	\]
	which has determinant $-1$. Hence, the normalized
	Haar measure on $\mathbb T^3$ is preserved. Using
	$dm_4=(2\pi)^{-4}\,dt\,du\,dv\,dw$ and $dt=2\,d\theta$,
	we have
	\begin{align*}
		&\int_{\mathbb T^4}
		G(\Phi(\xi_1,\xi_2,\lambda_1,\lambda_2))
		\pi\left|\frac{\lambda_1(\xi_1-1)}{2}\right|
		\sqrt{1-\left|\frac{\lambda_1(\xi_1-1)}{2}\right|^2}
		\,dm_4\\
		&=\frac{1}{(2\pi)^3}
		\int_0^\pi\int_{[0,2\pi]^3}
		G\bigl(\Phi(\xi_1,\xi_2,\lambda_1,\lambda_2)\bigr)
		\sin\theta\,|\cos\theta|\,
		du\,dv\,dw\,d\theta\\
		&=\frac{1}{(2\pi)^3}
		\int_0^{\pi/2}\int_{[0,2\pi]^3}
		G\bigl(U(\alpha,\beta,\gamma,r)\bigr)
		\sin r\cos r\,
		d\alpha\,d\beta\,d\gamma\,dr\\
		&\quad+
		\frac{1}{(2\pi)^3}
		\int_0^{\pi/2}\int_{[0,2\pi]^3}
		G\bigl(U(\alpha,\beta,\gamma,r)\bigr)
		\sin r\cos r\,
		d\alpha\,d\beta\,d\gamma\,dr\\
		&=\frac{1}{(2\pi)^3}
		\int_0^{\pi/2}\int_{[0,2\pi]^3}
		G\bigl(U(\alpha,\beta,\gamma,r)\bigr)
		2\sin r\cos r\,
		d\alpha\,d\beta\,d\gamma\,dr\\
		&=\int_{U(2)}G(U)\,dm_{U(2)}(U).
	\end{align*}
	Thus, $m_{U(2)}=
	\Phi_*\left(
	\pi|a\circ\Psi_{\mathbb H}|
	\sqrt{1-|a\circ\Psi_{\mathbb H}|^2}\,m_4
	\right)$, where $a \circ \Psi_{\mathbb H}$ denotes the projection map onto the first coordinate of $\Psi_{\mathbb H}$. Then
	
	\begin{align*}
		\int_{b\mathbb H}F\,d\sigma_{\mathbb H}
		&=\int_{U(2)}F(\pi_{\mathbb H}(U))\,dm_{U(2)}(U)
		\qquad [\text{since }\sigma_{\mathbb H}=(\pi_{\mathbb H})_*m_{U(2)}]\\
		&=\int_{\mathbb T^4}
		F(\pi_{\mathbb H}(\Phi(\xi_1,\xi_2,\lambda_1,\lambda_2)))
		\,\pi|a\circ\Psi_{\mathbb H}| \times \sqrt{1-|a\circ\Psi_{\mathbb H}|^2}\,dm_4\\
		&\qquad \qquad \qquad \qquad \qquad [\text{by the Haar integration formula under }\Phi]\\
		&=\int_{\mathbb T^4}
		F(\Psi_{\mathbb H}(\xi_1,\xi_2,\lambda_1,\lambda_2))
		\,\pi|a\circ\Psi_{\mathbb H}|\times
		\sqrt{1-|a\circ\Psi_{\mathbb H}|^2}\,dm_4\\
		&\qquad\qquad \qquad \qquad \qquad \qquad \qquad [\text{since }\pi_{\mathbb H}\circ\Phi=\Psi_{\mathbb H}]\\
		&=\int_{b\mathbb H}
		F(a,x_1,x_2,x_3)\,
		\pi|a|\sqrt{1-|a|^2}\,d\mu_{\mathbb H}\\
		&\qquad [\text{since }\mu_{\mathbb H}=(\Psi_{\mathbb H})_*m_4
		\text{ and the pushforward integration formula holds}]
	\end{align*}
	for every $F\in C(b\mathbb H)$. Therefore,
	\[
	\frac{d\sigma_{\mathbb H}}{d\mu_{\mathbb H}}
	(a,x_1,x_2,x_3)
	=
	\pi|a|\sqrt{1-|a|^2}
	\]
	almost everywhere on $b\mathbb H$. The proof is now complete.
\end{proof}

We now define another Hardy type space on $\mathbb H$ associated with the measure $\sigma_{\mathbb H}$.

\begin{definition}
	The weighted Hardy space for $\mathbb H$ associated with $\sigma_{\mathbb H}$ is defined by
	\[
	H^2_\sigma(\mathbb H)
	=
	\overline{\mathbb C[a,x_1,x_2,x_3]}
	^{\,L^2(b\mathbb H,\sigma_{\mathbb H})}.
	\]
\end{definition}

In view of the last proposition, the norm on this space can also be expressed in terms of the measure $\mu_{\mathbb H}$. Indeed, for every polynomial $f\in\mathbb C[a,x_1,x_2,x_3]$,
\[
\|f\|_{H^2_\sigma(\mathbb H)}^2
=
\int_{b\mathbb H}|f(a,x_1,x_2,x_3)|^2
\,d\sigma_{\mathbb H}
=
\pi\int_{b\mathbb H}
|f(a,x_1,x_2,x_3)|^2
|a|\sqrt{1-|a|^2}\,
d\mu_{\mathbb H}.
\]
Thus, $H^2_\sigma(\mathbb H)$ can equivalently be viewed as the completion of the polynomial algebra with respect to the weighted norm determined by
\[
w(a,x_1,x_2,x_3)
=
\pi|a|\sqrt{1-|a|^2}.
\]
Note that $0\leq w(a,x_1,x_2,x_3)\leq \pi \slash 2$ for every $(a, x_1, x_2, x_3) \in b\mathbb H$. It follows that
\[
\|f\|_{H^2_\sigma(\mathbb H)}
\leq
\sqrt{\frac{\pi}{2}}\,
\|f\|_{H^2(\mathbb H)}.
\]
for every polynomial $f$. Consequently, the identity map on the polynomial algebra extends to a bounded linear map
\[
J:H^2(\mathbb H)\to H^2_\sigma(\mathbb H).
\]
Moreover, the map $J$ is injective. To see this, assume that
$f\in H^2(\mathbb H)$ satisfies $Jf=0$. Choose a sequence of
polynomials $(p_n)_{n\geq 1}$ in
$\mathbb C[a,x_1,x_2,x_3]$ such that
\[
p_n\to f
\quad\text{in }H^2(\mathbb H).
\]
Since $J$ is continuous and $Jf=0$ in $H^2_\sigma(\mathbb H)$, it follows that $Jp_n \to Jf$ in $H^2_\sigma(\mathbb H)$ and $Jp_n=p_n$ for every $n \geq 1$. Thus, we have
\[
\|p_n\|_{H^2_\sigma(\mathbb H)}^2=\|Jp_n-Jf\|_{H^2_\sigma(\mathbb H)}^2=\int_{b\mathbb H}|p_n|^2w\,d\mu_{\mathbb H}
\to 0.
\]
By passing to suitable subsequences, we may assume that
$p_{n_k}\to f$
$\mu_{\mathbb H}$-almost everywhere and $p_{n_k}\to 0$
$w\,d\mu_{\mathbb H}$-almost everywhere. Since $w>0$ $\mu_{\mathbb H}$-almost everywhere, the measures
$w\,d\mu_{\mathbb H}$ and $d\mu_{\mathbb H}$ have the same
null sets. Consequently, $p_{n_k}\to 0$
$\mu_{\mathbb H}$-almost everywhere. Comparing the two
almost-everywhere limits, we have that $f=0$
$\mu_{\mathbb H}$-almost everywhere and thus $f=0$ in
$H^2(\mathbb H)$.

\smallskip

Next, we show that the norms induced by the measures $\mu_{\mathbb H}$ and
$\sigma_{\mathbb H}$ are not equivalent by considering the sequence
of polynomials
\[
f_n(a,x_1,x_2,x_3)=a^n,\qquad n\geq 1.
\]
Recall that $d\sigma_{\mathbb H}(a,x_1,x_2,x_3)
=
\pi |a|\sqrt{1-|a|^2}\,
d\mu_{\mathbb H}(a,x_1,x_2,x_3)$. First, we calculate the norm of $f_n$ with respect to
$\mu_{\mathbb H}$. By the definition of $\mu_{\mathbb H}$, we have

\begin{align*}
	\|f_n\|_{L^2(b\mathbb H,\mu_{\mathbb H})}^2
	=\int_{b\mathbb H}|a|^{2n}\,
	d\mu_{\mathbb H}(a,x_1,x_2,x_3)
	=\int_{\mathbb T^4}
	\left|\frac{\lambda_1(\xi_1-1)}{2}\right|^{2n}
	\,dm_4(\xi_1,\xi_2,\lambda_1,\lambda_2).
\end{align*}
Since $|\lambda_1|=1$, the integrand depends only on $\xi_1$.
Writing $\xi_1=e^{it}$, we obtain
\begin{align*}
	\|f_n\|_{L^2(b\mathbb H,\mu_{\mathbb H})}^2
	=\frac{1}{2\pi}\int_0^{2\pi}
	\left|\frac{e^{it}-1}{2}\right|^{2n}\,dt
	=\frac{1}{2\pi}\int_0^{2\pi}
	\sin^{2n}\left(\frac{t}{2}\right)\,dt
	=\frac{(2n)!}{2^{2n}(n!)^2}.
\end{align*}
Here, we have used the standard integral formula
\[
\int_0^{\pi/2}\sin^{2n}\theta\,d\theta
=
\frac{\pi}{2}\frac{(2n)!}{2^{2n}(n!)^2}.
\]
Next, we calculate the norm of $f_n$ with respect to
$\sigma_{\mathbb H}$. Note that
\begin{align*}
	\|f_n\|_{L^2(b\mathbb H,\sigma_{\mathbb H})}^2
	&=\int_{b\mathbb H}|a|^{2n}\,
	d\sigma_{\mathbb H}(a,x_1,x_2,x_3)\\
	&=\pi\int_{b\mathbb H}
	|a|^{2n+1}\sqrt{1-|a|^2}\,
	d\mu_{\mathbb H}(a,x_1,x_2,x_3)\\
	&=\pi\int_{\mathbb T^4}
	\left|\frac{\lambda_1(\xi_1-1)}{2}\right|^{2n+1}
	\sqrt{1-\left|\frac{\lambda_1(\xi_1-1)}{2}\right|^2}
	\,dm_4(\xi_1,\xi_2,\lambda_1,\lambda_2).
\end{align*}
Again, since $|\lambda_1|=1$, writing $\xi_1=e^{it}$ gives
\begin{align*}
	\|f_n\|_{L^2(b\mathbb H,\sigma_{\mathbb H})}^2
	=\frac{\pi}{2\pi}\int_0^{2\pi}
	\left|\sin\frac{t}{2}\right|^{2n+1}
	\sqrt{1-\sin^2\frac{t}{2}}\,dt
	=\int_0^\pi
	\sin^{2n+1}\theta\,|\cos\theta|\,d\theta
	=\frac{1}{n+1}.
\end{align*}
Consequently, we have
\begin{align}\label{eqn_ratio}
	\frac{
		\|f_n\|_{L^2(b\mathbb H,\sigma_{\mathbb H})}^2
	}{
		\|f_n\|_{L^2(b\mathbb H,\mu_{\mathbb H})}^2
	}
	=\frac{1}{n+1}
	\frac{2^{2n}(n!)^2}{(2n)!} \to 0 \quad \text{as} \ n \to \infty.
\end{align}
Now, assume on the contrary that the two norms were equivalent. Then there
exists a constant $C>0$ such that
\[
\|f\|_{L^2(b\mathbb H,\mu_{\mathbb H})}
\leq C\|f\|_{L^2(b\mathbb H,\sigma_{\mathbb H})}
\]
for every $f\in\mathbb C[a,x_1,x_2,x_3]$.
In particular, 
\[
\frac{
	\|a^n\|_{L^2(b\mathbb H,\mu_{\mathbb H})}
}{
	\|a^n\|_{L^2(b\mathbb H,\sigma_{\mathbb H})}
}
\leq C
\]
for every $n\geq1$, which is a contradiction to \eqref{eqn_ratio}. Hence, the two norms induced by $\mu_{\mathbb H}$ and $\sigma_{\mathbb H}$ are not equivalent. Thus, the construction $H^2_\sigma(\mathbb H)$ above gives a weighted Hardy space on $\mathbb H$ that arises from the Haar measure on the unitary group $U(2)$. 

\smallskip 

Next, we present an equivalent description of the measure $\sigma_\mathbb H$. Let $\omega_2$ be the normalized rotationally invariant surface measure on the topological boundary  $\partial\mathbb B_2$ of $\mathbb B_2$. Recall from \cite{BiswasPalII} that the map $f:\mathbb T\times\partial\mathbb B_2
\to b\mathbb H$ given by
\[
f(\eta,z,w)=(-\eta z,w,\eta\overline w,\eta)
\]
is a homeomorphism. Consequently, we define a probability measure $\tau_{\mathbb H}$ on $b\mathbb H$ as
\[
\tau_{\mathbb H}=f_*(m_1\times\omega_2).
\]
Equivalently, for every $F\in C(b\mathbb H)$,
\[
\int_{b\mathbb H}F(a,x_1,x_2,x_3)\,
d\tau_{\mathbb H}(a,x_1,x_2,x_3)
=
\int_{\mathbb T}\int_{\partial\mathbb B_2}
F(-\eta z,w,\eta\overline w,\eta)
\,d\omega_2(z,w)\,dm_1(\eta).
\]
We claim that $\tau_{\mathbb H}=\sigma_{\mathbb H}$. Consider the map
\[
\Theta:\mathbb T\times\partial\mathbb B_2
\to U(2)
\quad 
\text{defined by}
\quad 
\Theta(\eta,z,w)=
\begin{pmatrix}
	w&\overline z\\
	-\eta z&\eta\overline w
\end{pmatrix}.
\]
Since $|z|^2+|w|^2=1$ and $|\eta|=1$, a direct
calculation shows that $\Theta(\eta,z,w)\in U(2)$. We now show that
$\Theta_*(m_1\times\omega_2)=m_{U(2)}$. To do so, note that
\[
\Theta(\eta,z,w)
=
\begin{pmatrix}
	1&0\\
	0&\eta
\end{pmatrix}
\begin{pmatrix}
	w&\overline z\\
	-z&\overline w
\end{pmatrix}.
\]
The map $(z,w)\mapsto
\begin{pmatrix}
	w&\overline z\\
	-z&\overline w
\end{pmatrix}$ identifies $\partial\mathbb B_2$ with the set $SU(2)$ of
unitary matrices having determinant one in a bijective manner. Under this bijection, the normalized surface measure
$\omega_2$ is the normalized Haar measure on this set. 
To verify the above measure identity, let $G\in C(U(2))$.
Using the definition of $\Theta$, we have
\begin{align*}
	\int_{\mathbb T\times\partial\mathbb B_2}
	G(\Theta(\eta,z,w))\,d(m_1\times\omega_2)(\eta,z,w)
	=
	\int_{\mathbb T}\int_{\partial\mathbb B_2}
	G\left(
	\begin{pmatrix}
		1&0\\
		0&\eta
	\end{pmatrix}
	\begin{pmatrix}
		w&\overline z\\
		-z&\overline w
	\end{pmatrix}
	\right)
	\,d\omega_2(z,w)\,dm_1(\eta).
\end{align*}
Let $m_{SU(2)}$ denote the normalized Haar measure on
$SU(2)$. The map
\[
\kappa:\partial\mathbb B_2\to SU(2),
\qquad
\kappa(z,w)=
\begin{pmatrix}
	w&\overline z\\
	-z&\overline w
\end{pmatrix},
\]
is a homeomorphism. Moreover, since $\kappa$ is an
isometry between the unit sphere $\mathbb B_2$ in $\mathbb C^2$ and
$SU(2)$, the normalized rotationally invariant surface
measure $\omega_2$ is carried by $\kappa$ to
$m_{SU(2)}$. Thus,
\[
\int_{\partial\mathbb B_2}
G(\Theta(\eta,z,w))\,d\omega_2(z,w)
=
\int_{SU(2)}G(D_\eta V)\,dm_{SU(2)}(V),
\quad 
\text{where} \ \ 
D_\eta=
\begin{pmatrix}
	1&0\\
	0&\eta
\end{pmatrix}.
\]
We now use the standard decomposition of Haar measure
on $U(2)$. Every $U\in U(2)$ can be uniquely written as $U=D_\eta V$, where 
$\eta=\det U\in\mathbb T$ and $V\in SU(2)$. Furthermore, the normalized Haar measure on $U(2)$ is the product measure under this parametrization.
Consequently, for every $G\in C(U(2))$,
\[
\int_{U(2)}G(U)\,dm_{U(2)}(U)
=
\int_{\mathbb T}\int_{SU(2)}
G(D_\eta V)\,dm_{SU(2)}(V)\,dm_1(\eta).
\]
Combining the above identities, we have that
\begin{align*}
	\int_{\mathbb T\times\partial\mathbb B_2}
	G(\Theta(\eta,z,w))\,
	d(m_1\times\omega_2)(\eta,z,w)
	&=\int_{\mathbb T}\left[\int_{\partial\mathbb B_2}
	G(\Theta(\eta,z,w))\,d\omega_2(z,w)\right]dm_1(\eta)\\
	&=
	\int_{\mathbb T}\int_{SU(2)}
	G(D_\eta V)\,dm_{SU(2)}(V)\,dm_1(\eta)\\
	&=
	\int_{U(2)}G(U)\,dm_{U(2)}(U).
\end{align*}
Therefore, $\Theta_*(m_1\times\omega_2)=m_{U(2)}$. A simple calculation shows that
\[
\pi_{\mathbb H}(\Theta(\eta,z,w))
=(-\eta z,w,\eta\overline w,
\det\Theta(\eta,z,w))
=(-\eta z,w,\eta\overline w,\eta)
=f(\eta,z,w)
\]
and thus, $f=\pi_{\mathbb H}\circ\Theta$. Consequently, it follows that
\begin{align*}
	\tau_{\mathbb H}
	=f_*(m_1\times\omega_2)
	=(\pi_{\mathbb H}\circ\Theta)_*
	(m_1\times\omega_2)
	=(\pi_{\mathbb H})_*
	\big(\Theta_*(m_1\times\omega_2)\big)
	=(\pi_{\mathbb H})_*m_{U(2)}
	=\sigma_{\mathbb H}.
\end{align*}
Putting everything together, we arrive at the following result.

\begin{theorem}
	Consider the spaces 
	\[
	H^2(\mathbb H)
	=
	\overline{\mathbb C[a,x_1,x_2,x_3]}
	^{\,L^2(b\mathbb H,\mu_{\mathbb H})}
	\quad {and} \quad H^2_\sigma(\mathbb H)
	=
	\overline{\mathbb C[a,x_1,x_2,x_3]}
	^{\,L^2(b\mathbb H,\sigma_{\mathbb H})},
	\]
	where $\mu_{\mathbb H}$ is the probability measure induced by
	$\Psi_{\mathbb H}$ and $\sigma_{\mathbb H}$ is the probability
	measure induced by the normalized Haar measure on $U(2)$.
	Then the following statements hold:
	
	\begin{enumerate}
		\item The measures $\mu_{\mathbb H}$ and $\sigma_{\mathbb H}$
		are mutually absolutely continuous with
		\[
		\frac{d\sigma_{\mathbb H}}{d\mu_{\mathbb H}}
		(a,x_1,x_2,x_3)
		=
		\pi |a|\sqrt{1-|a|^2}
		\]
		almost everywhere on $b\mathbb H$.
		
		\item The identity map on $\mathbb C[a,x_1,x_2,x_3]$
		extends to a bounded injective linear operator $J:H^2(\mathbb H)\to H^2_\sigma(\mathbb H)$ satisfying
		\[
		\|Jf\|_{H^2_\sigma(\mathbb H)}
		\leq
		\sqrt{\frac{\pi}{2}}\,
		\|f\|_{H^2(\mathbb H)},
		\qquad f\in H^2(\mathbb H).
		\]
		
		\item The norms induced by $\mu_{\mathbb H}$ and
		$\sigma_{\mathbb H}$ on $\mathbb C[a,x_1,x_2,x_3]$
		are not equivalent. 
		
		\item Let $\omega_2$ be the normalized rotationally invariant
		surface measure on
		\[
		\partial\mathbb B_2
		=
		\{(z,w)\in\mathbb C^2:|z|^2+|w|^2=1\},
		\]
		and let $m_1$ be the normalized Haar measure on $\mathbb T$.
		Define $f:\mathbb T\times\partial\mathbb B_2\to b\mathbb H$
		by $f(\eta,z,w)=(-\eta z,w,\eta\overline w,\eta)$ and set $\tau_{\mathbb H}=f_*(m_1\times\omega_2)$. Then
		$
		\tau_{\mathbb H}=\sigma_{\mathbb H}.
		$
	\end{enumerate}
\end{theorem}

\section{Operator theory on the hexablock via the map $\Psi_{\mathbb H}$}\label{sec007_hex}

\noindent In Section \ref{sec002_hex}, we studied the scalar mapping $\Psi_{\mathbb H}$ on the closed polydisc $\mathbb D^4$, which played a pivotal role in the construction of the Hardy space $H^2(\mathbb H)$ associated with hexablock $\mathbb H$. We now turn to the operator-theoretic aspects of this mapping. In particular, we explore how $\Psi_{\mathbb H}$ relates operator tuples on the polydisc $\mathbb D^4$ to operator tuples associated with the hexablock. We begin by recalling the relevant notions of $\mathbb H$-contractions and their dilations.

\begin{definition}
	A commuting quadruple $(A, X_1, X_2, X_3)$ of operators acting on a Hilbert space $\mathcal H$ is said to be  an \textit{$\mathbb{H}$-contraction} if 
	\[
	\|q(A, X_1, X_2, X_3)\| \leq \|q\|_{\infty, \overline{\mathbb H}}
	\]
	for every $q \in \mathbb{C}[a, x_1, x_2, x_3]$. We say that an $\mathbb H$-contraction $(A, X_1, X_2, X_3)$ acting on a Hilbert space $\mathcal H$ admits dilation to an $\mathbb H$-unitary if there exists an $\mathbb H$-unitary $(N_1, N_2, N_3, N_4)$ on a Hilbert space $\mathcal K \supseteq \mathcal H$ such that
	\[
	q(A, X_1, X_2, X_3)=P_{\mathcal H}q(N_1, N_2, N_3, N_4)|_{\mathcal H}
	\]
	for every $q \in \mathbb C[a, x_1, x_2, x_3]$. Here, $P_{\mathcal H}$ is the orthogonal projection of $\mathcal K$ onto $\mathcal H$.
\end{definition}

Several equivalent descriptions of $\mathbb H$-unitaries have been developed in \cite{PalTomarI}. In view of the scalar identity $b\mathbb H=\Psi_{\mathbb H}(\mathbb T^4)$ established in Section \ref{sec003_hex}, it is natural to seek a corresponding operator-theoretic formulation. The following theorem provides such a characterization by relating $\mathbb H$-unitaries to commuting unitary quadruples. Its proof follows the same approach employed in \cite{AglerYoung1999} for operator pairs associated with the symmetrized bidisc.

\begin{theorem}\label{thm_H_unitary}
	A commuting quadruple $(N_1, N_2, N_3, N_4)$ of operators acting on a Hilbert space $\mathcal H$ is an $\mathbb H$-unitary if and only if there exist commuting quadruple $(U_1, U_2, U_3, U_4)$ of unitaries on $\mathcal H$
	such that $\Psi_{\mathbb H}(U_1, U_2, U_3, U_4)=(N_1, N_2, N_3, N_4)$, i.e., 
	\[
	N_1=\frac{U_3(U_1-I)}{2},\ \
	N_2=\frac{U_2(U_1+I)}{2},\ \
	N_3=\frac{U_3U_4(U_1+I)}{2},\ \
	N_4=U_1U_2U_3U_4.
	\]
\end{theorem}

\begin{proof} Let $\underline{N}=(N_1,N_2,N_3, N_4)$ be an $\mathbb H$-unitary. Then there exists a unique spectral measure $E$ on $\sigma_T(\underline N)$ such that
	\[
	N_j=\int_{\sigma_T(\underline N)} z_j\,dE(z),\qquad j=1,2,3, 4.
	\]
	By Theorem \ref{prop:Psi-interior}, the map $\Psi_{\mathbb H}:\mathbb T^4 \to b\mathbb H$ is onto and continuous. Since $\mathbb T^4$ and $b\mathbb H$ are compact spaces, one can find a Borel measurable map $\rho_{\mathbb H}:b\mathbb H\to \mathbb T^4$ such that $\Psi_{\mathbb H}\circ\rho_{\mathbb H}=\operatorname{id}_{b\mathbb H}$. Write $\rho_{\mathbb H}=(\rho_{\mathbb H}^{(1)}, \rho_{\mathbb H}^{(2)}, \rho_{\mathbb H}^{(3)}, \rho_{\mathbb H}^{(4)})$ and $U_j=\rho_{\mathbb H}^{(j)}(N_1, N_2, N_3, N_4)$ for $1 \leq j \leq 4$, where $U_j$ is defined via the Borel functional calculus. Since each $|\rho_{\mathbb H}^{(j)}|=1$ on $\sigma_T(\underline N)$, each $U_j$ is a unitary operator. Also, $U_1, U_2, U_3, U_4$ commute since they lie in the commutative von Neumann algebra generated by $N_1, N_2, N_3, N_4$. It now follows from Borel functional calculus that $\underline{N}=\Psi_{\mathbb H}(U_1, U_2, U_3, U_4)$. 
	
	\smallskip 
	
	Conversely, suppose $U_1,U_2,U_3, U_4$ are commuting unitaries and let us write $(N_1, N_2, N_3, N_4)=\Psi_{\mathbb H}(U_1, U_2, U_3, U_4)$. Evidently, $N_1, N_2, N_3, N_4$ are commuting normal operators. We have by spectral mapping theorem that
	\[
	\sigma_T(N_1,N_2,N_3, N_4)
	=
	\Psi_{\mathbb H}\bigl(\sigma_T(U_1,U_2,U_3, U_4)\bigr)
	\subseteq
	\Psi_{\mathbb H}(\mathbb T^4)
	=
	b\mathbb H,
	\]
	which completes the proof. 
\end{proof}

The characterization of $\mathbb H$-unitaries in Theorem \ref{thm_H_unitary} gives the following result.

\begin{corollary}
	Let $(V_1,V_2,V_3, V_4)$ be a commuting quadruple of isometries on a Hilbert space $\mathcal H$. Then $\Psi_{\mathbb H}(V_1,V_2,V_3, V_4)$ is an $\mathbb H$-isometry.
\end{corollary}

\begin{proof}
	By the dilation theorem for commuting isometries, there exist a Hilbert space $\mathcal K$ containing $\mathcal H$ and a commuting quadruple of unitaries $(U_1,U_2,U_3, U_4)$ on $\mathcal K$ such that $\mathcal H$ is invariant under each $U_j$ and $V_j=U_j|_{\mathcal H}$ for $1\leq j \leq 4$.
	By Theorem \ref{thm_H_unitary}, $\Psi_{\mathbb H}(U_1,U_2,U_3, U_4)$ is an $\mathbb H$-unitary and
	\[
	\Psi_{\mathbb H}(V_1,V_2,V_3, V_4)
	=
	\Psi_{\mathbb H}(U_1,U_2,U_3, U_4)|_{\mathcal H}.
	\]
	Thus, $\Psi_{\mathbb H}(V_1,V_2,V_3, V_4)$ is the restriction of an $\mathbb H$-unitary to a joint invariant subspace and hence, it is an $\mathbb H$-isometry.
\end{proof}

Constructing examples of $\mathbb H$-contractions directly can be a rather delicate endeavour. The mapping $\Psi_{\mathbb H}$ provides a systematic approach to this problem by generating $\mathbb H$-contractions from suitable commuting operator quadruples of contractions. Beyond this, the following result identifies a family of $\mathbb H$-contractions that admit dilations to $\mathbb H$-unitaries.

\begin{proposition}
	Let $(T_1, T_2, T_3, T_4)$ be a commuting quadruple of contractions on a Hilbert space $\mathcal H$ for which $\overline{\mathbb D}^4$ is a spectral set for $(T_1, T_2, T_3, T_4)$. Then $\Psi_{\mathbb H}(T_1, T_2, T_3, T_4)$ is an $\mathbb H$-contraction. In particular, if $(T_1, T_2, T_3, T_4)$ dilates to a commuting quadruples of unitaries, then $\Psi_{\mathbb H}(T_1, T_2, T_3, T_4)$ dilates to an $\mathbb H$-unitary.
\end{proposition}

\begin{proof}
	Let $q$ be a holomorphic polynomial in four variables. Then $g=q\circ\Psi_{\mathbb H}$ is holomorphic on $\overline{\mathbb D}^{\,4}$. Since $\overline{\mathbb D}^{\,4}$ is a spectral set for $(T_1,T_2,T_3,T_4)$, we have 
	\begin{align*}
		\|q(\Psi_{\mathbb H}(T_1,T_2,T_3,T_4))\|
		=\|g(T_1,T_2,T_3,T_4)\|
		\leq \|g\|_{\infty,\overline{\mathbb D}^{\,4}}
		\leq \|q\|_{\infty,\overline{\mathbb H}},
	\end{align*}
	where the last inequality follows from the fact that $\Psi_{\mathbb H}(\overline{\mathbb D}^{\,4})
	\subseteq\overline{\mathbb H}$ (see Theorem \ref{prop:Psi-interior}). Thus, $\overline{\mathbb H}$ is a spectral set for
	$\Psi_{\mathbb H}(T_1,T_2,T_3,T_4)$ and so, it is an $\mathbb H$-contraction. Next, assume that $(T_1,T_2,T_3,T_4)$ dilates to a commuting quadruple of unitaries $(U_1,U_2,U_3,U_4)$ acting on a Hilbert space $\mathcal K$ containing $\mathcal H$. Then
	\begin{align*}
		q(\Psi_{\mathbb H}(T_1,T_2,T_3,T_4))
		=g(T_1,T_2,T_3,T_4)
		=P_{\mathcal H}g(U_1,U_2,U_3,U_4)|_{\mathcal H}
		=P_{\mathcal H}
		q(\Psi_{\mathbb H}(U_1,U_2,U_3,U_4))
		|_{\mathcal H}.
	\end{align*}
	Since $(U_1,U_2,U_3,U_4)$ is a commuting quadruple of unitaries, we have by Theorem \ref{thm_H_unitary} that the quadruple $\Psi_{\mathbb H}(U_1,U_2,U_3,U_4)$ is an $\mathbb H$-unitary. Consequently, $\Psi_{\mathbb H}(U_1,U_2,U_3,U_4)$ is an $\mathbb H$-unitary dilation of $\Psi_{\mathbb H}(T_1,T_2,T_3,T_4)$. This completes the proof.
\end{proof}

\section{Operators associated with the tetrablock via the map $\Psi_{\mathbb E}$}\label{sec008_hex}

\noindent We studied in Section \ref{sec007_hex}, the operator-theoretic aspects of the map $\Psi_{\mathbb H}$, which first appeared in Section \ref{sec003_hex}. In the same section, we discussed another parametrization map $\Psi_{\mathbb E}$ associated with the tetrablock.  Recall that the polynomial map $\Psi_{\mathbb E}:\mathbb C^3\to\mathbb C^3$ is defined as
\[
\Psi_{\mathbb E}(\lambda_1, \lambda_2, \xi)=\left(\frac{\xi(\lambda_1+\lambda_2)}{2}, \frac{\lambda_1+\lambda_2}{2}, \xi \lambda_1 \lambda_2\right) 
\]
and the above map satisfies $\Psi_{\mathbb E}(\mathbb D^2 \times \overline{\mathbb D})\subseteq\mathbb E$ and $\Psi_{\mathbb E}(\mathbb T^3)=b\mathbb E$. We begin by recalling from \cite{Bhatt2014} the notions of $\mathbb E$-contractions and the dilations to $\mathbb E$-unitaries.

\begin{definition}
	A commuting triple $(X_1, X_2, X_3)$ of operators acting on a Hilbert space $\mathcal H$ is said to be an \textit{$\mathbb{E}$-contraction} if 
	\[
	\|q(X_1, X_2, X_3)\| \leq \|q\|_{\infty, \overline{\mathbb E}}
	\]
	for every $q \in \mathbb{C}[x_1, x_2, x_3]$. A commuting triple $(X_1, X_2, X_3)$ of normal Hilbert space operators is said to be an $\mathbb E$-unitary if $\sigma_T(X_1, X_2, X_3) \subseteq b\mathbb E$. An $\mathbb E$-isometry is a commuting triple of operators which is the restriction of an $\mathbb E$-unitary to a joint invariant subspace.

	\smallskip 
	
	An $\mathbb E$-contraction $(X_1, X_2, X_3)$ on a Hilbert space $\mathcal H$ is said to dilate to an $\mathbb E$-unitary if there exists an $\mathbb E$-unitary $(N_1, N_2, N_3)$ on a Hilbert space $\mathcal K \supseteq \mathcal H$ such that
	\[
	q(X_1, X_2, X_3)=P_{\mathcal H}q(N_1, N_2, N_3)|_{\mathcal H}
	\]
	for every $q \in \mathbb C[x_1, x_2, x_3]$.
\end{definition}

The author of \cite{Bhatt2014} obtained various characterizations of $\mathbb E$-unitaries. Next, we present an alternative characterization of $\mathbb E$-unitaries via the map $\Psi_{\mathbb E}$.

\begin{theorem}\label{thm_E_unitary}
	A commuting quadruple $(N_1, N_2, N_3)$ of operators acting on a Hilbert space $\mathcal H$ is an $\mathbb E$-unitary if and only if there exist commuting triple $(U_1, U_2, U_3)$ of unitaries on $\mathcal H$
	such that $\Psi_{\mathbb E}(U_1, U_2, U_3)=(N_1, N_2, N_3)$, i.e., 
	\[
	N_1=\frac{U_3(U_1+U_2)}{2},\quad
	N_2=\frac{U_1+U_2}{2},\quad 
	N_3=U_1U_2U_3.
	\]
\end{theorem}

\begin{proof} Let $\underline{N}=(N_1,N_2,N_3)$ be an $\mathbb E$-unitary. Then there exists a unique spectral measure $E$ on $\sigma_T(\underline N)$ such that
	\[
	N_j=\int_{\sigma_T(\underline N)} z_j\,dE(z),\qquad j=1,2,3.
	\]
	By Theorem \ref{prop:Psi-interior}, the map $\Psi_{\mathbb E}:\mathbb T^3 \to b\mathbb E$ is onto and continuous. Since $\mathbb T^3$ and $b\mathbb E$ are compact spaces, there exists a Borel measurable map $\rho_{\mathbb E}:b\mathbb E\to \mathbb T^3$ such that $\Psi_{\mathbb E}\circ\rho_{\mathbb E}=\operatorname{id}_{b\mathbb E}$. Write $\rho_{\mathbb E}=(\rho_{\mathbb E}^{(1)}, \rho_{\mathbb E}^{(2)}, \rho_{\mathbb E}^{(3)})$ and $U_j=\rho_{\mathbb E}^{(j)}(N_1, N_2, N_3)$ for $1 \leq j \leq 3$, where $U_j$ is defined via the Borel functional calculus. Since each $|\rho_{\mathbb E}^{(j)}|=1$ on $\sigma_T(\underline N)$, each $U_j$ is a unitary operator. Also, $U_1, U_2, U_3$ commute as each of them belongs to the commutative von Neumann algebra generated by $N_1, N_2, N_3$. It now follows from Borel functional calculus that $\underline{N}=\Psi_{\mathbb E}(U_1, U_2, U_3)$. Conversely, let $U_1,U_2,U_3$ be commuting unitaries and let $(N_1, N_2, N_3)=\Psi_{\mathbb E}(U_1, U_2, U_3)$. Evidently, $N_1, N_2, N_3$ are commuting normal operators. By spectral mapping theorem, $\sigma_T(N_1,N_2,N_3)=
	\Psi_{\mathbb E}\bigl(\sigma_T(U_1,U_2,U_3)\bigr)
	\subseteq
	\Psi_{\mathbb E}(\mathbb T^3)
	=
	b\mathbb E$. The proof is now complete. 
\end{proof}

Next, we have the following result providing a way to construct $\mathbb E$-isometries.

\begin{corollary}
	Let $(V_1,V_2,V_3)$ be a commuting triple of isometries on a Hilbert space $\mathcal H$. Then $\Psi_{\mathbb E}(V_1,V_2,V_3)$ is an $\mathbb E$-isometry.
\end{corollary}

\begin{proof}
	By the dilation theorem for commuting isometries, there exist a Hilbert space $\mathcal K$ containing $\mathcal H$ and a commuting triple of unitaries $(U_1,U_2,U_3)$ on $\mathcal K$ such that $\mathcal H$ is invariant under each $U_j$ and $V_j=U_j|_{\mathcal H}$ for $1\leq j \leq 3$.
	By Theorem \ref{thm_E_unitary}, $\Psi_{\mathbb E}(U_1,U_2,U_3)$ is an $\mathbb E$-unitary and
	$
	\Psi_{\mathbb E}(V_1,V_2,V_3)
	=
	\Psi_{\mathbb E}(U_1,U_2,U_3)|_{\mathcal H}.
	$
	Thus, $\Psi_{\mathbb E}(V_1,V_2,V_3)$ is an $\mathbb E$-isometry.
\end{proof}

Next, we show that the map $\Psi_{\mathbb E}$ also provides a systematic approach to construct examples of $\mathbb E$-contractions.

\begin{proposition}
	Let $(T_1, T_2, T_3)$ be a commuting triple of contractions on a Hilbert space $\mathcal H$ for which $\overline{\mathbb D}^3$ is a spectral set for $(T_1, T_2, T_3)$. Then $\Psi_{\mathbb E}(T_1, T_2, T_3)$ is an $\mathbb E$-contraction. In particular, if $(T_1, T_2, T_3)$ dilates to a commuting triples of unitaries, then $\Psi_{\mathbb E}(T_1, T_2, T_3)$ dilates to an $\mathbb E$-unitary.
\end{proposition}

\begin{proof}
	Let $q$ be a holomorphic polynomial in three variables. Then $g=q\circ\Psi_{\mathbb E}$ is holomorphic on $\overline{\mathbb D}^{3}$. Since $\overline{\mathbb D}^{3}$ is a spectral set for $(T_1,T_2,T_3)$, we have 
	\begin{align*}
		\|q(\Psi_{\mathbb E}(T_1,T_2,T_3))\|
		=\|g(T_1,T_2,T_3)\|
		\leq \|g\|_{\infty,\overline{\mathbb D}^{3}}
		\leq \|q\|_{\infty,\overline{\mathbb E}},
	\end{align*}
	where the last inequality follows from the fact that $\Psi_{\mathbb E}(\overline{\mathbb D}^{3})
	\subseteq\overline{\mathbb E}$ (see Theorem \ref{prop:Psi-interior}). Thus, $\overline{\mathbb E}$ is a spectral set for
	$\Psi_{\mathbb E}(T_1,T_2,T_3)$ and so, it is an $\mathbb E$-contraction. Next, assume that $(T_1,T_2,T_3)$ dilates to a commuting triple of unitaries $(U_1,U_2,U_3)$ acting on a Hilbert space $\mathcal K$ containing $\mathcal H$. Then
	\begin{align*}
		q(\Psi_{\mathbb E}(T_1,T_2,T_3))
		=g(T_1,T_2,T_3)
		=P_{\mathcal H}g(U_1,U_2,U_3)|_{\mathcal H}
		=P_{\mathcal H}
		q(\Psi_{\mathbb E}(U_1,U_2,U_3))
		|_{\mathcal H}.
	\end{align*}
	Since $(U_1,U_2,U_3)$ is a commuting triple of unitaries, we have by Theorem \ref{thm_E_unitary} that $\Psi_{\mathbb E}(U_1,U_2,U_3)$ is an $\mathbb E$-unitary. Consequently, $\Psi_{\mathbb E}(U_1,U_2,U_3)$ is an $\mathbb E$-unitary dilation of $\Psi_{\mathbb E}(T_1,T_2,T_3)$. 
\end{proof}

%-----------------------------------------------------------------------
\smallskip 

\noindent \textbf{Acknowledgments.} 
The first named author is supported by the `Core Research Grant (CRG)' with Award No. CRG/2023/005223 of Anusandhan National Research Foundation (ANRF), Govt. of India. The second named author is partially supported by the National Natural Science Foundation of China (Nos.\ 1231101308, \ 12671102) and the Natural Science Foundation of Shandong Province (No.\ ZR2025MS77). The third named author is supported via the Postdoctoral Fellowship (Ref. No. 0204/9(10)/2026-R\&D-II/4805) awarded by the National Board for Higher Mathematics (NBHM), Government of India.

\bibliographystyle{plain}   
\bibliography{references} 
\end{document}